\documentclass[10pt, english]{amsart}

\usepackage{Style}

\usepackage{geometry}

\usepackage[utf8]{inputenc}
\usepackage{amsmath}
\usepackage{amsfonts}
\usepackage{amssymb}
\usepackage{comment}
\usepackage{hyperref}
\usepackage{enumitem,pstricks,xy,pst-node}
\xyoption{all}
\usepackage{amsthm}
\usepackage{booktabs} 
\usepackage[color]{changebar} 
\cbcolor{red} 
\usepackage{mathrsfs} 
\usepackage{verbatim} 
\usepackage[disable] 
{todonotes}
\usepackage[normalem]{ulem} 
\usepackage{caption} 
\usepackage{nicefrac} 
\usepackage{enumitem} 
\usepackage{pdfpages}
\usepackage[dvipsnames]{xcolor}

\newcommand{\tG}{{{\widetilde{\Gamma}}}}
\newcommand{\PP}{\mathbb{P}}

\newcommand{\fX}{{{\mathfrak{X}}}}

\newcommand\formal[1]{{[\![{#1}]\!]}}
\newcommand{\uk}{{\underline{k}}}

\def\git{/\!\!/}

\def\cF{{\mathcal F}}

\def\cU{{\mathcal U}}
\def\cM{{\mathcal M}}
\def\cS{{\mathcal S}}
\def\cY{{\mathcal Y}}

\def\bbP{\mathbb P}

\def\CC{\mathbb{C}}

\def\QQ{\mathbb{Q}}
\def\QQQ{\mathbb{Q}}
\def\ZZ{\mathbb{Z}}

\def\cO{{\mathcal O}}
\def\cF{{\mathcal F}}

\def\cM{{\mathcal M}}

\def\fg{\mathfrak{g}}

\def\af1{\mathbf{aff}_1}

\def\cD{{\mathcal D}}

\DeclareMathOperator{\eu}{Eu}

\DeclareMathOperator{\im}{Im}

\DeclareMathOperator{\Sym}{Sym}

\DeclareMathOperator{\HHH}{H}

\DeclareMathOperator{\Gr}{Gr}
\DeclareMathOperator{\IGr}{IGr}

\DeclareMathOperator{\pt}{pt}

\DeclareMathOperator{\Stab}{Stab}

\DeclareMathOperator{\amb}{amb}

\DeclareMathOperator{\van}{van}
\DeclareMathOperator{\QH}{QH}
\DeclareMathOperator{\QHo}{{\overline{QH}}}
\DeclareMathOperator{\BQH}{BQH}
\DeclareMathOperator{\BQHo}{{\overline{BQH}}}
\DeclareMathOperator{\QHoa}{{\overline{QH}}^{amb}}

\DeclareMathOperator{\QHa}{{QH}^{amb}}
\DeclareMathOperator{\HHHa}{{H}^{amb}}
\DeclareMathOperator{\HHHna}{{H}^{van}}
\DeclareMathOperator{\BQHa}{{ BQH}^{amb}}

\DeclareMathOperator{\Ham}{H^{amb}}

\newcommand{\tone}{{t_{1}}}
\newcommand{\ttwo}{{t_{2}}}

\DeclareMathOperator{\OGP}{O}

\newcommand{\ladi}{\begin{lastadd}}
\newcommand{\ladf}{\end{lastadd}}
\newcommand{\lrei}{\begin{lastrem}}
\newcommand{\lref}{\end{lastrem}}
\newenvironment{lastadd}
{\cbstart\color{red}}
{\todo{red to remove}\cbend}
\newenvironment{lastrem}
{\cbstart\color{yellow}}
{\cbend}

\author{Vladimiro Benedetti, Simone Billi, Nicolas Perrin}

\address{
Universit\'e C\^ote d'Azur, CNRS, Laboratoire J.-A. Dieudonn\'e, Parc Valrose, F-06108 Nice Cedex 2, France} 
\email{vladimiro.benedetti@univ-cotedazur.fr}

\address{Fakultät für Mathematik und Informatik, Universität des Saarlandes, Campus E2.4, 66123 Saarbrücken, Germany} 
\email{billi@math.uni-sb.de}

\address{
Centre de Math\'ematiques Laurent Schwartz (CMLS), CNRS, \'Ecole polytechnique,
Institut Polytechnique de Paris, F-91120 Palaiseau, France}
\email{nicolas.perrin.cmls@polytechnique.edu}

\title{Monodromy and big quantum cohomology of hyperplane sections of (co)adjoint varieties}
\begin{document}
\maketitle

\begin{abstract}
    We study the monodromy representation of general hyperplane sections of (co)adjoint varieties. We apply this result to reduce the generic semisemplicity of the big quantum cohomology to its ambient part, and to prove that the big quantum cohomology of general hyperplane sections of \(\IGr(2,6)\) is generically semisimple and satisfies Dubrovin's conjecture.
\end{abstract}

\section{Introduction}

Dubrovin's conjecture \cite{dubrovin1998geometry} relates the semisimplicity of the big quantum cohomology \(\BQH(X)\) of a smooth Fano variety \(X\) to the existence of a full exceptional collection in \(\text{D}^b(X)\) the bounded derived category of coherent sheaves on \(X\). This conjecture has many interesting refinements giving more precise information on both sides of the conjecture (see for example \cite{CDG} or \cite{KS}). Understanding for which Fano varieties \(X\) the big quantum cohomology \(\BQH(X)\) is semisimple is therefore significant. A stronger condition is the semisimplicity of \(\QH(X)\) the small quantum cohomology which is a quotient of \(\BQH(X)\). In \cite{KS}, under the condition that \(\BQH(X)\) is semisimple, the algebraic structure of \(\QH(X)\) and especially of its radical gives conjectural strong conditions on the structure of full exceptional collections in \(\text{D}^b(X)\).

\medskip

Rational projective homogeneous spaces are an important source of Fano varieties. It is therefore natural to search for examples among them. The semisimplicity of \(\QH(X)\) and \(\BQH(X)\) for \(X\) a rational projective homogeneous space has been studied in several works (see for example \cite{iritani-koto} and \cite{semisimplicity}). In \cite{ChPe} the first example of rational projective homogeneous spaces \(X\) with \(\QH(X)\) not semisimple appeared, while in \cite{galkin2015dubrovin} and \cite{perrin2014semisimple} the first examples for which \(\BQH(X)\) is semisimple but \(\QH(X)\) is not were presented. Adjoint or coadjoint varieties (see \autoref{def:adj-coadj}) are among the simplest examples of such homogeneous spaces. This case is completely settled in \cite{maxim}, where the semisimplicity of \(\BQH(X)\) is proved and the structure of the radical in \(\QH(X)\) is described. Several conjectural result predicted by refinements of Dubrovin's conjecture were proved in further works about (co)adjoint varieties (see for example \cite{KS1} and \cite{KS}).

\medskip

Another important source of Fano varieties is provided by zero loci \(Y\) of general sections of equivariant vector bundles on projective rational homogeneous spaces \(X\). Among those, complete intersections and linear sections are the easiest to understand. The study of their quantum cohomology attracted interest for a long time (see for example \cite{beauville-qh}) but it was recently applied to important rationality questions (see \cite{kkpy}, \cite{guere}, \cite{GM} or \cite{verra}). A systematic way to compute the ambient part of \(\BQH(Y)\) along deformations of ambient classes (see Subsection \ref{subsection:qh} for more details) has recently been proposed in \cite{atomic.ci.flag}.

\medskip

This paper focuses on quantum cohomology of a hyperplane section \(Y\) of a (co)adjoint variety \(X\) (see \autoref{def:adj-coadj}). In a previous paper \cite{BP}, a description of the primitive cohomology was obtained as well as partial results on semisimplicity and non-semisimplicity of \(\QH(Y)\). These results were complete in \cite{belmans2025ade}. In particular \(\QH(Y)\) is not semisimple for \(Y\) a smooth hyperplane section of a coadjoint variety. The main purpose of this paper is to perform computations in \(\BQH(Y)\) with a view toward proving its semisimplicity. A key ingredient, among others, is the study of the monodromy representation for such hyperplane sections. From \cite{BP}, if \(Y\) is such an hyperplane section, then the primitive cohomology \(\HHHna(Y,\mathbb{Z})\) lies in the middle cohomology and it is isometric to an associated root lattice \(\Lambda_\aleph\) with Weyl group \(W_\aleph\). 

\begin{theorem}[{\autoref{thm:main_theorem_monodromy}}]\label{thm:main_theorem_monodromy_intro}
Let \(Y\subset X\) be a smooth hyperplane section of a (co)adjoint variety \(X\).
    There is only one non-trivial subrepresentation of the monodromy representation on \(\HHH^*(Y,\ZZ)\), which factors through the natural representation of \(W_\aleph\) on \(\Lambda_\aleph\).
\end{theorem}
We refer to \autoref{sec:monodromy} for precise definitions and statements. The Weyl group coincides in fact with the entire monodromy group when the associated Dynking diagram has no exterior automorphisms (\autoref{cor:monodromy_no_exterior}) or when it is simply-laced in the adjoint case (\autoref{cor:monodromy_simplylaced} and \cite{belmans2025ade}).

\medskip

From the knowledge of the monodromy representation and using the invariance of Gromov-Witten under deformation, we prove the following criterion for semisimplicity.

 \begin{theorem}[{\autoref{thm_QBH_nonAn}}]
 \label{thm-intro-reduction}
     Let \(Y\subset X\) be a smooth hyperplane section of a (co)adjoint variety \(X\) for a simply laced group. Let \(\eta\in \HHHa(Y)\) be an ambient class and consider deformations of the quantum cohomology along \(\eta\). If the ambient part of the quantum cohomology \(\BQH_\eta^{\mathrm{amb}}(Y)\) is generically semisimple, then so is \(\BQH(Y)\).
 \end{theorem}
 
 This result together with the recent computational software developed in the paper \cite{atomic.ci.flag} offer many perspectives at least for simply laced groups (see \autoref{prop:QBH_nonAn} for the general case).
 
 \medskip
 
We then focus on cases not covered by \autoref{thm-intro-reduction}: a coadjoint variety \(X\) for a non-simply laced group \(G\) and \(Y \subset X\) a smooth hyperplane section. In types \(\text{B}_n\) and \(\text{G}_2\) both \(X\) and \(Y\) are smooth quadrics for which the semisimplicity of \(\QH(Y)\) is well known. We thus consider the remaining cases: type \(\text{C}_n\) and type \(\text{F}_4\). We first give a presentation of \(\QH(Y)\) recovering the fact that this algebra is not semisimple (see \cite{BP} and \cite{belmans2025ade}). 

To state this result we introduce some notations. Let \(X\) be coadjoint not adjoint for the simple group \(G\) and \(Y \subset X\) a smooth hyperplane section. Let \(s\) be the number of short roots in the Dynkin diagram of \(G\), let \(i_Y\) be the index of \(Y\), let \(h \in \HHH^2(Y)\) be the hyperplane class and let \((\Gamma_\alpha)_{\alpha \in \aleph}\) be the basis of \(\HHHna(Y)\) obtained from the Bia{\l}ynicki-Birula decomposition (see Subsection \ref{subsection:BB-decomposition-and-coho-basis}). Define \(\QHo(Y) := \QH(Y)/(q - 1)\) (see Subsection \ref{subsection:qh}).

\begin{theorem}[\autoref{thm_pres_smallQHbar}]
There exists an ambient cohomology class \(\tau \in \HHHa(Y)\) and a polynomial \(P\) with only simple non-zero roots such that
\begin{align*}
    \QHo(Y)\simeq \QQ[h, \tau, (\Gamma_\alpha)_{\alpha\in \Phi_\aleph}]/\!\left(hP(h^{i_Y}),\tau^s, h\tau, \{h\Gamma_\alpha\}_{\alpha\in \Phi_\aleph},\{\Gamma^2_\alpha - \Gamma_\beta^2\}_{\alpha\neq \beta\in \Phi_\aleph}\right.,\\ 
    \nonumber\left.\{\Gamma^2_\alpha - 2\Gamma_\alpha \Gamma_\beta\}_{\alpha\neq \beta\in \Phi_\aleph}, \{\tau \Gamma_\alpha\}_{\alpha\in \Phi_\aleph}, \{\Gamma^2_\alpha-\tau^{s-1}\}_{\alpha\in \Phi_\aleph}\right).
\end{align*}
\end{theorem}

This result together with the refinement of Dubrovin's conjecture in \cite{KS} give strong contrains on the expected exceptional collections in \(\text{D}^b(Y)\). Note also that the above result applies in types \(\text{B}_n\) and \(\text{G}_2\) in which case \(s =1\) and all classes \(\Gamma_\alpha\) vanish so that we indeed recover Beauville's presentation \cite{beauville-qh} of the small quantum cohomology of quadrics and its semisimplicity.

\medskip

We finish the paper with a new semisimplicity result for the big quantum cohomology. We focus on the simplest case where \(Y\) is a hyperplane section of \(\IGr(2,6)\)\footnote{Such hyperplane sections are also called \textit{bisymplectic Grassmannians}} for which we consider deformations along primitive classes. We prove the following result
 
 \begin{theorem}[\autoref{cor:semi-semplicity_for_IG(2,6)}]\label{thm:main_thm_intro}
    Let \(Y\subset \IGr(2,6)\) be a smooth hyperplane section, then the big quantum cohomology ring \(\BQH(Y)\) is generically semisimple.
\end{theorem}

Dubrovin's conjecture \cite{dubrovin1998geometry} predicts that the semisimplicity of \(\BQH(Y)\) should be equivalent to the existence of a full exceptional collections for the bounded derived category \(D^b(Y)\). It was proved by Kuznetsov \cite{kuznetsov2008exceptional} that for a hyperplane section \(Y\subset \IGr(2,6)\) such full exceptional collection exists, which together with our result implies Dubrovin's conjecture in this case.
 
 \begin{corollary}[\autoref{cor:Dubrovin_for_IG(2,6)}]
     Dubrovin's conjecture holds for smooth hyperplane sections \(Y\subset \IGr(2,6)\).
 \end{corollary}
 

\bigskip\noindent {\it Acknowledgements}. V.B. and N.P. were supported by the project FanoHK ANR-20-CE40-0023. S.B. would like to thank Stevell Muller for useful discussions.

\setcounter{tocdepth}{2}
\tableofcontents

\section{(Co)adjoint varieties and their hyperplane sections}
All varieties are defined over \(\mathbb{C}\). For \(Z\) a smooth variety, we write \(\HHH^*(Z,\ZZ)\) for its integral cohomology and \(\HHH^*(Z) := \HHH^*(Z,\ZZ) \otimes_\ZZ\QQ\) for its rational cohomology.

\subsection{(Co)adjoint varieties}
\label{subsection:coadjoint}

Let \(G\) be a connected and simply connected simple Lie group and denote by \(\fg\) its Lie algebra.  Fix \(T \subset B \subset G\) a maximal torus and a Borel subgroup of \(G\) and denote by \(W = N_G(T)/T\) the Weyl group of \(G\) associated to \(T\). Denote by \(\Delta\) the set of roots of \(G\), by \(\Delta^+\) the set of positive roots and set \(\Delta^- = - \Delta^+\). Let \(\Phi \subset \Delta^+\) be the set of simple roots. Denote by \(\Delta_l\) and \(\Delta_s\) the set of long roots and short roots respectively and set \(\Delta_l^\pm = \Delta_l \cap \Delta^\pm\) and  \(\Delta_s^\pm = \Delta_s \cap \Delta^\pm\). Denote by \(\Theta \in \Delta_l^+\) and \(\theta \in \Delta_s^+\) the highest root and the highest short root. If \(G\) is simply laced, we have \(\Delta = \Delta_l = \Delta_s\) and \(\Theta = \theta\). For a dominant weight \(\lambda\) of \(G\), denote by \(\nabla_\lambda\) the highest weight representation of highest weight \(\lambda\). Note that we have \(\nabla_\Theta = \fg\). The representations \(\nabla_\Theta\) and \(\nabla_\theta\) will play an important role in this paper.

\begin{definition}[Adjoint and coadjoint varieties]\label{def:adj-coadj} Let \(G\) be as above.
\begin{enumerate}
\item The \(G\)-\textit{adjoint variety} is the unique closed \(G\)-orbit inside \(\bbP(\nabla_\Theta)\). 
\item The \(G\)-\textit{coadjoint variety} is the unique closed \(G\)-orbit inside \(\bbP(\nabla_\theta)\). 
\end{enumerate}
\end{definition}

From now on \(X\) will be either the adjoint or coadjoint variety (we will specify when a distinction needs to be made). Set \(V = \nabla_\Theta\) in the former case and \(V = \nabla_\theta\) in the latter case. Note that in both cases, we have an isomorphism \(V \simeq V^\vee\) of \(G\)-representations. The embedding \(X \subset \bbP(V)\) is the Pl\"ucker embedding of \(X\) \emph{i.e.} the minimal embedding of \(X\) in a projective space. Denote by $\aleph$ the \(W\)-orbit of the highest weight of \(V\). If \(X\) is \(G\)-adjoint, we have \(\aleph = W.\Theta = \Delta_l\) while if \(X\) is \(G\)-coadjoint, we have \(\aleph = W.\theta = \Delta_s\). We set $\Phi_\aleph = \Phi \cap \aleph$.

\medskip 

Let \(Y\subset X\) be a general smooth hyperplane section. By \cite[Theorem 1.2]{BP}, the stabilizer \(\Stab_G(Y)\) of \(Y\) in \(G\) contains a maximal torus. Up to conjugation, we may (and will) assume that \(Y\) is stable by the maximal torus \(T \subset G\). Since \(X\) has only finitely many \(T\)-stable curves and \(T\)-fixed points, the same holds for \(Y\). More precisely one can show (see \cite{BP}) that \(X^T=Y^T\). The Bia{\l}ynicki-Birula decomposition therefore induces a natural cellular decompositions in \(X\) and \(Y\) and cohomology basis. We discuss these facts in more details in the next subsection.

\subsection{Bia{\l}ynicki-Birula decomposition and cohomology basis}
\label{subsection:BB-decomposition-and-coho-basis}


Let \(Z\) be a smooth projective variety acted on by a one dimensional torus \(\top\) such that the set of fixed points \(Z^\top\) is finite. Identify \(\top\) with \(\CC^*\) and embed it in \(\mathbb{P}^1\). For \(z_0 \in Z^\top\), set \(\Omega_{z_0,Z}^+ = \{ z \in Z \ | \ \lim_{t \to \infty}t \cdot z = z_0\}\) and \(\Omega_{z_0,Z}^- = \{ z \in Z \ | \ \lim_{t \to 0}t \cdot z = z_0\}\). Let \(Z_{z_0}^\pm\) be the Zariski closure in \(Z\) of  \(\Omega_{z_0,Z}^\pm\) and set \(\sigma_{z_0,Z}^\pm = [Z_{z_0}^\pm] \in \HHH^*(Z,\ZZ)\). The classes \((\sigma_{z_0,Z}^+)_{z_0 \in Z^\top}\) and \((\sigma_{z_0,Z}^-)_{z_0 \in Z^\top}\) form two basis of \(\HHH^*(Z,\ZZ)\). We will refer to \(\Omega_{z_0,Z}\) and \(\Omega_{z_0,Z}^-\) as the \emph{Schubert} cells and \emph{opposite Schubert} cells of \(Z\) and to \(Z_{z_0}^+\) and \(Z_{z_0}^-\) as the \emph{Schubert} varieties and \emph{opposite Schubert} varieties of \(Y\), in analogy with the homogeneous case.

We may apply the above setting to \(Y \subset X \subset \bbP(V)\). Indeed, if \(\lambda\) is a non-zero weight of \(V\) and \(V_\lambda = \{ v \in V \ | \ t.v = \lambda(t) v, \forall t \in T\}\) is the corresponding weight space, we have \(\dim(V_\lambda) = 1\) and \(X^T = Y^T = \{ x_\alpha = [V_\alpha] \in \bbP(V) \ | \ \alpha \in \aleph\}\). For \(X\), this leads to the usual Schubert (and opposite Schubert) cells, varieties and classes. We set \(\Omega^\pm_{\alpha,X} := \Omega^{\pm}_{x_\alpha,X}\) , \(X_{\alpha}^\pm := X^{\pm}_{x_\alpha}\)  and \(\sigma^{\pm}_{\alpha,X} := \sigma_{x_\alpha,X}\). For \(Y\), we set \(\Omega_{\alpha}^\pm := \Omega^\pm_{x_\alpha,Y}\) , \(Y_{\alpha}^\pm := Y_{x_\alpha}^\pm\)  and \(\sigma_{\alpha}^\pm := \sigma_{x_\alpha,Y}^\pm\). We shall often drop the exponent for the positive basis and write $\sigma_{\alpha,X}$ for $\sigma_{\alpha, X}^+$ and $\sigma_\alpha$ for $\sigma_\alpha^+$. In particular \((\sigma_{\alpha,X})_{\alpha \in \aleph}\) is a basis of \(\HHH^*(X,\ZZ)\) and \((\sigma_\alpha)_{\alpha\in \aleph}\) is a basis of \(\HHH^*(Y,\ZZ)\).

\medskip 

We now recall a few facts from \cite{BP} on middle cohomology of $Y$ that are useful to study the monodromy action in Section \ref{sec:monodromy}. We shall give more results such as cup product and pull-back formulae in Section \ref{sec:quantum}. 

For \(\alpha \in \aleph\), we have (\(\sigma_{\alpha} \in \HHH^{\dim(Y)}(Y)\))  \(\Leftrightarrow\) (\(\alpha\) or \(-\alpha\) is simple) \(\Leftrightarrow\) (\(\pm \alpha \in \Phi_\aleph\)).

\begin{definition} 
\label{def:vanish-coho}
Denote by \(j\colon Y\hookrightarrow X\) the embedding of \(Y\) in \(X\).
\begin{enumerate}
\item The \emph{ambient cohomology} \(\HHHa(Y,\ZZ)\) of \(Y\) is defined as \(\HHHa(Y,\ZZ)=\im(j^*\mid_{\HHH^*(X,\ZZ)})\).
\item  The \emph{vanishing cohomology} \(\HHHna(Y,\ZZ)\) of \(Y\) is defined as \(\HHHna(Y,\ZZ)=\Ker(j_*\mid_{\HHH^*(Y,\ZZ)})\).
\end{enumerate}
\end{definition}

Following our convention, we set \(\HHHa(Y) := \HHHa(Y,\ZZ) \otimes_\ZZ \QQ\) and  \(\HHHna(Y) := \HHHna(Y,\ZZ) \otimes_\ZZ \QQ\)

\begin{remark}
A few comments about the previous definitions are in order.
\begin{enumerate}
\item Whenever the cohomological degree of the ambient cohomology needs to be specified, we will denote the degree \(k\) part of \(\HHHa(Y,\ZZ)\)  and \(\HHHa(Y)\) as \(\HHH^k(Y,\ZZ)^{\amb}\) and \(\HHH^k(Y)^{\amb}\).  
\item Classical results imply that \(\HHHna(Y,\mathbb{Z})\subset \HHH^{\dim Y}(Y,\mathbb{Z})\) (see \cite[Corollary II.2.25]{voisin2002theorie}) so we will not need to specify the degree when working with vanishing cohomology. 
\item Since \(X\) has no primitive cohomology in degree \(\dim Y\), Proposition II.2.27 in \cite{voisin2002theorie} implies that \(\HHHna(Y,\ZZ)\) coincides with primitive classes in \(\HHH^{\dim Y}(Y,\ZZ)\). 
\end{enumerate}
\end{remark}

\begin{definition}
For \(\alpha \in \Phi_\aleph\), set \(\Gamma_\alpha = \sigma_\alpha - \sigma_{-\alpha}\) and \(\Gamma_\alpha^- = \sigma_{-\alpha}^- - \sigma_{\alpha}^-\).
\end{definition}

The following was proved in \cite[Lemma 4.24]{BP}.

\begin{lemma}
\label{lemma:gamma-gen-van}
The classes \((\Gamma_\alpha)_{\alpha \in \Phi_\aleph}\) and \((\Gamma_\alpha^-)_{\alpha \in \Phi_\aleph}\) form two basis of \(\HHHna(Y)\).
\end{lemma}

 Let \(I_\aleph\) be the identity matrix on the set \(\aleph\) \emph{i.e.} \(I_\aleph = (\delta_{\alpha,\beta})_{\alpha,\beta \in \Phi_\aleph}\) and let \(C_\aleph\) be the Cartan matrix associated to \(\aleph\) and defined by \(C_\aleph = (\langle \alpha^\vee,\beta \rangle)_{\alpha,\beta \in \Phi_\aleph}\). Note that \(C_\aleph\) is symmetric since all roots have the same length in \(\aleph\).

\begin{theorem}[{\cite[Theorem 4.36]{BP}}]
\label{thm_inter_product_middle_cohom}
  We have the following intersection matrices
  \begin{equation*}
      (\Gamma_\alpha \cup \Gamma_\beta)_{\alpha,\beta \in \Phi_\aleph} = (\Gamma^-_\alpha \cup \Gamma^-_\beta)_{\alpha,\beta \in \Phi_\aleph} =
  (\sqrt{-1})^{\dim Y}(4I_\aleph - C_\aleph).
  \end{equation*}
\end{theorem}

We recall that here and in \cite{BP}, we have that \(\dim Y\) is always an even number.
\section{The monodromy representation}\label{sec:monodromy}

In this section, we study the action of monodromy on families of smooth hyperplane sections of adjoint and coadjoint varieties and prove \autoref{thm:main_theorem_monodromy_intro}. Hyperplane sections are parametrized by \(V\simeq V^\vee=\HHH^0(X,\mathcal{O}(1))^\vee\). For \(h \in V\), set 
\begin{equation*}
Y_h = \{x \in X \ | \ h(x) = 0\} \subset X \subset \bbP(V). 
\end{equation*}

Smooth hyperplane sections are parametrized by an open subset \(U\subset V\). 
Let \(V_0\subset V\) be the 0-eigenspace for the action of the Lie algebra of \(T\) on \(V\) (in the adjoint case \(V_0\) is a Cartan subalgebra of \(\fg\)). For \(h\in V_0\), a section \(Y_h \subset X\)  is smooth if and only if it lies in the complement \(U_0\) in \(V_0\) of the hyperplane arrangement of orthogonal spaces to roots in \(\aleph\) (for these results see \cite[Theorem 2.5]{tevelev} and \cite[Proposition 7.2]{manben}). Note that since \(G.V_0\) is dense in \(V\) (see \cite{BP}), the \(G\)-translate of \(U_0\) is a dense open subset of \(V\) such that \(G.U_0\subseteq U\). The adjoint case has been studied in great details. We review a few features of this case in the following example and refer to Tevelev's book \cite{tevelev} for more details.

\begin{example}[The adjoint case]
\label{exam:reg-ss-adjoint}
Assume that we are in the adjoint case \emph{i.e.} \(V = \fg\). In that case, all elements in \(U_0\) and therefore in \(G.U_0\) are regular and semisimple. The converse does not always hold but it holds for simply laced groups. 

The dual projective variety \(\cD_X := \{ [h] \in \bbP(V) \ | \ Y_h \text{ is singular}\}\) of \(X\) is the locus of singular hyperplane sections of \(X\). In particular, \(\cD_X\) is an hypersurface whose equation \(D_l\) is such that the restriction of \(D_l\) to any Cartan subalgebra is the product of long roots (see \cite[Theorem 2.5]{tevelev}). Let \(D\) be the equation whose restriction to any Cartan subalgebra is the product of all roots and \(D_s\) be the equation whose restriction to any Cartan subalgebra is the product of short roots (these equations exist by Chevalley's Theorem \(\fg /\!\!/G \cong  V_0/W\)). In particular \(D = D_lD_s\). For \(G\) simply laced \(D = D_l\) and \(D_s = 1\). By \cite[Lemma 2.1.9]{coll-mcgov}, an element \(h \in \fg = V\) is regular and semisimple if and only if \(D(h) \neq 0\). In particular, for \(G\) simply laced we have that \(Y_h\) is smooth if and only if \(h\) is regular and semisimple. For \(G\) non-simply laced, there are smooth hyperplane sections \(Y_h\) for which \(h\) is not semisimple.
\end{example}

Denote by \(\chi\colon \mathcal{F}\to F\) any family of smooth hyperplane sections in \(X\) with \(\cF\subset X\times F\). The action of parallel transport operators from a smooth hypersurface \(Y_h\) to itself, where \(h\in F\) and \(\chi(Y_h)=h\), is encoded by the monodromy representation of the family $F$: 
\begin{equation*}
    \pi_1(F,h)\to\Aut(\HHH^*(Y_h,\mathbb{Z})),
\end{equation*}
whose image is called the \textit{monodromy group} \(\Mon(Y_h,F)\) associated to the family. In this case, we actually have an equality \(\Mon(Y_h,F) \simeq \Mon^{\van}(Y_h,F)\), where the latter is the image of the map 
\begin{equation*}
    \rho_F\colon\pi_1(F,h)\to \OGP(\HHHna(Y_h,\mathbb{Z})),
\end{equation*}
where \(\OGP(\HHHna(Y_h,\mathbb{Z}))\) is the isometry group (for Poincaré pairing) of the vanishing cohomology.

\begin{definition}
In this setting, we define the \textit{monodromy group} \(\Mon(Y_h)\) (and its avatar \(\Mon^{\van}(Y_h)\)) as the group generated by \(\Mon(Y_h,F)\) (resp. \(\Mon^{\van}(Y_h,F)\)) for all families $\cF \to F$ of smooth hypersurfaces such that \(\cF\subset X\times F\).\footnote{This is an ad-hoc definition for this context, more generally the monodromy group is bigger since families \(\mathcal{F}\to F\) are not required to satisfy the property \(\mathcal{F}\subset X\times F\). See the discussion in \autoref{rmk:definition_of_monodromy}.} 
\end{definition}

The universal family is \(\cF_U := \{ (x,h) \in X \times U) \ | \ h(x) = 0 \} \to U\), where \(U\subset V\) is the set of \emph{all} smooth hyperplane sections in \(X\). Any family \(\cF \subset X \times F\) can be obtained from the universal family \(\cF_U \to U\) by pullback. We thus have an identification \(\Mon(Y_h)\simeq \Mon(Y_h, U)\) (see \cite{beauville-monodromie}). In particular, we set \(\rho:=\rho_U\) for the monodromy representation of the universal family. 

Define the family \(S_0 \subset V_0\) of hypersurfaces as follows. In the adjoint case, the space \(V_0\) is a Cartan subalgebra; let \(H_\alpha\in V_0\) be the coroot associated to the root \(\alpha\) and set \(S_0 :=U_0 \cap \sum_{\alpha \in \Delta_l} \bbC H_\alpha\).  In the coadjoint case, set \(S_0 := U_0\). The Weyl group \(W_\aleph\) of the root system \(\aleph\) acts on \(S_0\), and we have an isomorphism \(Y_h \simeq Y_{w(h)}\) for any \(w\in W_\aleph\). Set \(\cS_0 := \cF_U \times_U S_0\). Then \(\cS_0\) has a \(W_\aleph\)-action such that the map \(\cS_0 \to S_0\) is \(W_\aleph\)-equivariant. Consider the family of hyperplane sections \(S\) obtained by taking the quotient under this action. More precisely, setting \(S := S_0/W_\aleph\) and  \(\cS := \cS_0/W_\aleph\), we consider the family 
\begin{equation*}
\cS \to S.
\end{equation*}
Notice that \(W_\aleph\) acts on \(S_0\) freely because any isotropy group in a complex reflection group is generated by reflections \cite[Section 1.2]{humphreys1992reflection}, and the regularity condition, which is there to ensure that the hyperplane sections are smooth, implies that \(h \in S_0\) is fixed by no reflection. 

\begin{remark}
Recall from \autoref{exam:reg-ss-adjoint} that in the adjoint simply-laced case, \(S_0 = U_0\), \(G.S_0 = U\), \(W = W_\aleph\) and by the Chevalley Theorem \(S \simeq U \git W \simeq S_0 / W_\aleph\). The last isomorphism is a consequence of the fact that elements in \(G.V_0\) have closed orbit (because they are semisimple).
\end{remark}

Notice that \(\cS \subset (X \times S_0)/W_\aleph\) is not contained in \(X \times S\). However, for \(h \in S\), the following lemma shows that we can still define a monodromy action and a group \(\Mon(Y_h,S)\leq \Mon(Y_h)\). 

\begin{lemma}\label{lem:rhoB_exists}
    There is a well-defined map 
\(\varrho_S:\pi_1(S,h) \to \Aut(\HHH^*(Y_h,\ZZ))\) that factors through $\mathrm{Mon}(Y_h)\hookrightarrow \Aut(\HHH^*(Y_h,\ZZ))$.
\end{lemma}
\begin{proof}
Given a loop \(\gamma:[0,1]\to S\), we can lift this loop to \(\widetilde{\gamma}:[0,1]\to S_0\) such that \(\widetilde{\gamma}(0)=h\) and \(\widetilde{\gamma}(1)=w(h)\) for a certain \(w\in W_\aleph\). Let \(\widetilde{w}\) be a lift of \(w\in W_\aleph\leq W=N_G(T)/T\) to \(N_G(T)\), and let \(\delta:[0,1]\to G\) be a path such that \(\delta(0)=\widetilde{w}\) and \(\delta(1)=\id\) (this is possible since \(G\) is connected). We can extend the path \(\widetilde{\gamma}\) to a path \(\gamma':[0,1]\to G.S_0\) by concatenating \(\widetilde{\gamma}\) with \(\widetilde{\delta}(t):= \delta(t)\cdot (h)\), \(t\in [0,1]\); this is possible since \(w(h)=\widetilde{w}(h)\in S_0\subset V_0\) as \(T\) acts trivially on \(V_0\). By the simple connectedness of \(G\), the homotopy class of \(\gamma'\) does not depend on the choice of \(\delta\). Moreover \(\gamma'\) is a loop from \(h\) to itself. Let \([\gamma]\in \pi_1(S,h)\) (we denote by \(h\) also the point in \(S\) by abuse of notation). Then we define \(\varrho_S:\pi_1(S,h) \to \Aut(\HHH^*(Y_h,\ZZ))\) as:
\begin{equation*}
    \varrho_S([\gamma]):=\rho([\gamma']).
\end{equation*}
The definition does not depend on the lift \(\widetilde{w}\) of \(w\) since \(T\) fixes \(Y\) and is connected, so it acts trivially on cohomology.
\end{proof}

\begin{proposition}\label{prop:rhoB_morphism}
The map \(\varrho_S\) is a morphism of groups.
\end{proposition}

\begin{proof}
Let \(\circ\) denote concatenation of paths (the first loop being on the left) and, for \(l:[0,1]\to V\) a path, denote by \(-l\) the path \((-l)(t)=l(1-t)\) obtained by following \(l\) backwards. Then in order to show that \(\varrho_S\) is a morphism of groups it is sufficient to show that \((\gamma\circ\gamma_1)'\) is homotopic to \(\gamma' \circ\gamma'_1\), and then use the fact that \(\rho\) is a morphism of groups. Let us use the same notation as before; moreover let us suppose that \(\widetilde{\gamma_1}(0)=h\) and \(\widetilde{\gamma_1}(1)=\widetilde{w_1}(h)\). Let \(\delta_1:[0,1]\to G\) be a path such that \(\delta_1(0)=\widetilde{w_1}\) and \(\delta(1)=\id\), and define as before \(\widetilde{\delta_1}\). Then \(\gamma'=\widetilde{\gamma} \circ \widetilde{\delta}\) and \(\gamma'_1=\widetilde{\gamma_1} \circ \widetilde{\delta_1}\), while 
\[(\gamma\circ \gamma_1)'=\widetilde{\gamma} \circ (\widetilde{w}\cdot(\widetilde{\gamma_1})) \circ (\widetilde{\omega}\cdot(\widetilde{\delta_1})) \circ \widetilde{\delta}.\]
Let \(\eta:[0,1]\to G\) be a path such that \(\eta(0)=\id\) and \(\eta(1)=\widetilde{\omega}\), and \(\widetilde{\eta}(t):=\eta(t)\cdot(h)\). Moreover, for \(s\in [0,1]\) define by \(\widetilde{\eta}|_{[0,s]}\) the path obtained from \(\widetilde{\eta}\) by restricting to \([0,s]\). Then define the homotopy
\[\Gamma:[0,1]\times [0,1] \to G.S_0 \]
by
\[\Gamma(\bullet, s)=\widetilde{\gamma} \circ \widetilde{\delta} \circ \widetilde{\eta}|_{[0,s]} \circ (\eta(s)\cdot (\widetilde{\gamma_1})) \circ (-(\widetilde{(\eta\cdot \widetilde{\omega_1})}|_{[0,s]})) \circ \widetilde{\delta_1}.\]
Since \(G\) is simply connected, the paths \(\widetilde{\delta}\circ \widetilde{\eta}\) and \((-\widetilde{(\eta\cdot \widetilde{\omega_1})})\circ \widetilde{\delta_1} \circ (-\widetilde{\delta})\circ (-\widetilde{\omega}\cdot (\widetilde{\delta_1}))\) are homotopically trivial, so we deduce that \(\Gamma(\bullet, 1)\) is homotopically equivalent to \((\gamma\circ \gamma_1)'\); since \(\Gamma(\bullet, 0)=\gamma' \circ \gamma'_1\), we deduce that \(\gamma' \circ \gamma'_1\) and \((\gamma\circ \gamma_1)'\) are homotopically equivalent.
\end{proof}

\begin{remark}
\label{rem_parallel_and_action}
Let \(\Sigma\subset Y_{\tilde{w}(h)}\) be a subvariety and \(\sigma\) its class in cohomology. Parallel transport of \(\sigma\) along \(\tilde{\delta}\) exactly gives the class of \(\tilde{w}^{-1}(\Sigma)\subset Y_h\) in the cohomology of \(Y_h\). Thus in practice, to compute the action of \(\varrho_S\), we will do parallel transport along \(\tilde{\gamma}\) and then we will apply \(\tilde{w}^{-1}\).
\end{remark}

\begin{remark}\label{rmk:definition_of_monodromy}
For any smooth projective fibration \(p : \cF\to F\) if \(f\in F\), denote by \(\cF_f:=p^{-1}(f)\) the fiber over \(f\). In fact, in this more general case, following \cite[Section II.3.1.2]{voisin2002theorie}  one has a monodromy action \(\rho_F:\pi_1(F,f)\to \Aut(\HHH^*(\cF_f,\ZZ))\) with image \(\Mon(\cF_f,F)\) being the monodromy group of the family. When the family \(\cF\) satisfies \(\cF\subset X\times F\), then it is a subfamily of the universal family \(\cF_U \to U\) of smooth hyperplane sections, and so the monodromy group will be a subgroup of \(\Mon(Y_h)\) (where \(Y_h\) is a hyperplane section for \(h \in F\subset U\)). With this point of view, the above definition of \(\varrho_S\) together with \autoref{lem:rhoB_exists} and \autoref{prop:rhoB_morphism} actually show that 
\[\varrho_S = \rho_S \quad \mbox{ and }\quad \Mon(S,h)\leq \Mon(Y_h).\]
\end{remark}


\subsection{Upper bound on the monodromy group}

By definiton \(\Mon(Y_h)\) is a subgroup of \(\OGP(\HHHna(Y_h,\ZZ))\). We compute the latter. By \autoref{lemma:gamma-gen-van},
the classes \((\Gamma_\alpha)_{\alpha\in \Phi_\aleph}\) form a basis for \(\HHHna(Y_h,\mathbb{Z})\). We need a slightly different basis for the vanishing cohomology.
Assign a number \(\epsilon_\alpha\in\{\pm 1\}\) to each \(\alpha\in\Phi_\aleph\) so that when two distinct roots \(\alpha\) and \(\beta\) are not orthogonal, we have \(\epsilon_\alpha \epsilon_\beta=-1\) (this is possible since there is no loop in the Dynkin diagram). Define the alternative basis \((\tG_\alpha)_{\alpha \in\aleph}\) by \(\tG_\alpha \coloneqq \epsilon_\alpha \Gamma_\alpha\) for all \(\alpha \in \Phi_{\aleph}\).


\begin{lemma}\label{lemma:gram_matrix_new_basis}
We have the matrix equality 
\(({\tG_\alpha}\cup {\tG_\beta})_{\alpha,\beta \in \aleph}=(\sqrt{-1})^{\dim Y_h}C_\aleph\). 
\end{lemma}
\begin{proof} Set \(2k=\dim Y_h\). To simplify notation, we drop the subscript \(\aleph\) in the matrices \(I_\aleph\) and \(C_\aleph\). By \autoref{thm_inter_product_middle_cohom}, we have \((\Gamma_\alpha\cup \Gamma_\beta)_{\alpha,\beta}=(-1)^k(4I-C)\). Setting \(M:=4I-C\), we have:
\begin{itemize}
    \item \(M_{\alpha,\beta}=0\) if \(\alpha\perp \beta,\)
    \item \(M_{\alpha,\beta}=1\) if \(\alpha\not\perp \beta\) and \(\alpha\not=\beta\),
    \item \(M_{\alpha,\alpha}=2\).
\end{itemize}
In particular, notice that \(C_{\alpha,\beta}=-M_{\alpha,\beta}\) if \(\alpha\not=\beta\) and \(C_{\alpha,\alpha}=M_{\alpha,\alpha}\). We have:
\begin{itemize}
    \item If \(\alpha\perp \beta,\) then  \({\tG_\alpha}\cup {\tG_\beta}= \epsilon_\alpha\epsilon_\beta (\Gamma_\alpha\cup \Gamma_\beta)=0=(-1)^kC_{\alpha,\beta}\).
    \item If \(\alpha\not\perp \beta\) and \(\alpha\not=\beta\), then \(\epsilon_\alpha=-\epsilon_\beta\) and  
    \({\tG_\alpha}\cup {\tG_\beta}= -
    \Gamma_\alpha\cup \Gamma_\beta 
    =(-1)^{k+1}=(-1)^k C_{\alpha,\beta}\).
    \item Finally, \({\tG_\alpha}\cup {\tG_\beta}=\epsilon_\alpha^2(\Gamma_\alpha\cup\Gamma_\alpha)=2(-1)^k=(-1)^kC_{\alpha,\alpha}\).
\end{itemize}
The claim follows.
\end{proof}

Denote by \(\Lambda_\aleph\) the root lattice associated to the root system \(\Phi_\aleph\), and by \(D_\aleph\) the associated Dynkin diagram. The above \autoref{lemma:gram_matrix_new_basis} together with the fact that \(C_\aleph\) is precisely the Gram matrix of \(\Lambda_{\aleph}\), shows that we have a lattice isomorphism \(\Lambda_\aleph \simeq \HHHna(Y_h,\mathbb{Z})\).
\begin{lemma}\label{lemma:isometry_group_root_lattice}
    We have an isomorphism of groups 
    \[\OGP(\HHHna(Y_h,\mathbb{Z})) \simeq \OGP(\Lambda_\aleph) \simeq W_\aleph \rtimes\Aut(D_\aleph),\]
    where \(\Aut(D_\aleph)\) denotes the group of automorphisms of the diagram \(D_\aleph\).
\end{lemma}
\begin{proof}
Let \(\mathcal{C}\) be a chamber for the root system \(\Phi_\aleph\). 
The action of the subgroup \(W_\aleph\leq \OGP(\Lambda_\aleph)\) on the set of chambers of \(\Lambda_\aleph\) is free and transitive. We get an isomorphism \(\OGP(\Lambda_\aleph) = W_\aleph\rtimes \OGP_{\mathcal{C}}(\Lambda_\aleph)\) where \(\OGP_{\mathcal{C}}(\Lambda_\aleph)\) denotes the stabilizer of the chamber \(\mathcal{C}\). We conclude by observing that \(\OGP_{\mathcal{C}}(\Lambda_\aleph)\) can be identified with the automorphism group of the diagram \(D_\aleph\), after identifying the walls of the fundamental chamber with the vertices of the diagram. 
\end{proof}

\begin{corollary}\label{cor:isometry_group_middle_coh}
We have \(\Mon(Y_h)\leq  W_\aleph\rtimes\Aut(D_\aleph).\) 
\end{corollary}

\begin{proof}
We have seen that \(\Mon(Y_h)\leq \OGP(\HHHna(Y_h,\mathbb{Z})) \simeq W_\aleph\rtimes\Aut(D_\aleph)\).
\end{proof}

\subsection{The Weyl group acts via monodromy operators}
\label{subsection:weyl-gp-action}

We now use the family \(S\) to prove that \(W_\aleph \leq \Mon(Y_h)\). We have a natural action of \(T\) on both the tautological family \(\mathcal{U}_0=\{(x,h)\in X\times U_0\mid h(x)=0\}\) and the space \(X\times U_0\). Since \(T\) acts trivially on \(V_0 \supset U_0\), the \(T\)-fixed locus inside \(X\times U_0\) is given by \(    (X\times U_0)^T=X^T\times  U_0\).

We consider the Bia\l{}ynicki-Birula decompositions (BB-decomposition) for \(\cU_0\) and \(X\times U_0\). The restriction of the former to the fiber over \(h\in U_0\) give the 
BB-decomposition of \(Y_h\). Recall that \(X^T=\{x_\alpha\}_{\alpha\in\aleph}=Y_h^T\) for any \(h\in U_0\). Recall from Subsection \ref{subsection:BB-decomposition-and-coho-basis} the definition of \(\Omega_{\alpha,X}^+\) and set \(\Omega_{\alpha,X} = \Omega_{\alpha,X}^+\). The BB-decomposition for \(X\) is given by \(X = \bigsqcup_{\alpha\in \aleph} \Omega_{\alpha,X}\).
This decomposition trivially extends to a BB-decomposition \(X\times U_0=\bigsqcup_{\alpha\in \aleph} (\Omega_{\alpha,X} \times U_0)\).

We now apply the same BB-decomposition to \(Y_h\) and to the family \(\cS_0\subset X\times S_0\). We get a stratification
\(\cS_0 =\bigsqcup_{\alpha\in \aleph} \Omega_{\alpha,\cS_0}\)
where \(\Omega_{\alpha,\cS_0}=\{y \in \cS_0\mid \lim_{t\to\infty} t\cdot y\in x_\alpha \times S_0\}\)
 and set \((\cS_0)_\alpha := \overline{\Omega_{\alpha,\cS_0}}\). The fiber over \(h \in S_0\) of the fibration \(\phi\colon\cS_0\to S_0\) is given by \(Y_h = X\cap H\) where \(H=\{h=0\}\). In particular, this gives the BB-decomposition \(Y_h=\bigsqcup_{\alpha\in \aleph} \Omega_{\alpha,h}\)
with \(\Omega_{\alpha,h}=\{y\in Y_h\mid \lim_{t\to\infty} t\cdot y = x_\alpha\}\). Set \(Y_{\alpha,h} := \overline{\Omega_{\alpha,h}}\). The subvarieties $Y_{\alpha,h}\subset Y_h$ define classes in cohomology that we denote by $\sigma_{\alpha,h}$. Note that for fixed \(h\), we recover the classes \(\sigma_\alpha\) from Subsection \ref{subsection:BB-decomposition-and-coho-basis}.

\medskip

Let \(\gamma\colon [0,1]\to S_0\) be a continuous path. We also write \(\gamma\) for its image. It defines families
of hyperplane sections with fiber \(Y_h\) over \(h\in \gamma\) and of BB-cells in the fiber \(Y_h\):
\begin{equation*}
    \cY_\gamma = \{ (y,h) \in X \times S_0 \mid h(y)=0, h \in \gamma \} \subset \cS_0 \text{ and } \Omega_{\alpha,\gamma}=\{y\in \cY_\gamma\mid \lim_{t\to\infty} t\cdot y\in x_\alpha\times \gamma\}\subset \Omega_{\alpha,\cS_0}.
\end{equation*} 
Set \(\cY_{\alpha,\gamma} := \overline{\Omega_{\alpha,\gamma}} \subset (\cS_0)_\alpha\). 


\begin{lemma}\label{lemma:action_on_paths}
Let \(\gamma\colon [0,1]\to S_0\) be a continuous path. The induced parallel transport action \(\gamma^*\colon \HHHna(Y_{\gamma(1)},\ZZ)\to\HHHna(Y_{\gamma(0)},\ZZ)\) satisfies: 
\[\gamma^*(\sigma_{\alpha,\gamma(1)})=\sigma_{\alpha,\gamma(0)}.\]
\end{lemma}

\begin{proof}
The fibration \(\phi\colon \cS_0\to S_0\) determines the trivial local system \(R^{2k}\phi_*\mathbb{Z}\) with fiber \(\HHH^{2k}(Y_h,\mathbb{Z})=(R^{2k}\phi_*\mathbb{Z})_h\) over \(h\in S_0\). The path \(\gamma\) determines the family \(\phi_{\gamma}\colon \cY_\gamma\to[0,1]\); we get a (differentiable) trivialization by Ehresmann's Theorem: \[ \cY_\gamma\simeq \phi_\gamma^{-1}(0)\times[0,1]=Y_{\gamma(0)}\times [0,1]\]
which gives an homeomorphism \(\psi\colon Y_{\gamma(0)}\simeq Y_{\gamma(1)}\).
We have 
\(\gamma^*(\eta)=\psi^*(\eta)\)
for \(\eta\in \HHH^{2k}(Y_{\gamma(1)},\mathbb{Z})\). 

By definition, \(\psi^*(\sigma_{\alpha,\gamma(1)})=[\psi^{-1}(Y_{\alpha,\gamma(1)})]\). 
Identifying \(\cY_\gamma\) and \(Y_{\gamma(0)}\times [0,1]\), the subvariety \(\cY_{\alpha,\gamma} \subset Y_{\gamma(0)}\times [0,1]\) can be seen as a homotopy between \(Y_{\alpha,\gamma(0)}\) and \(\psi^{-1}(Y_{\alpha,\gamma(1)})\) inside \(Y_{\gamma(0)}\); by homotopy invariance of cohomology classes, we deduce that
\[\psi^*(\sigma_{\alpha,\gamma(1)})=[\psi^{-1}(Y_{\alpha,\gamma(1)})] = \sigma_{\alpha,\gamma(0)},\]
proving the claim.
\end{proof}

\begin{corollary}
\label{cor:path-invariance}
Closed paths in $S_0$ induce a trivial action on cohomology; thus the parallel transport action does not depend on the chosen path \(\gamma : [0,1] \to S_0\) but only on \(\gamma(0)\) and \(\gamma(1)\).
\end{corollary}

\begin{remark}
The conclusion of \autoref{lemma:action_on_paths} also holds for a continuous path \(\gamma : [0,1] \to U_0\).
\end{remark}


\begin{proposition}\label{prop_mon_induces_action_W}
    The monodromy action on \(\HHH^*(Y_h,\mathbb{Z})\) (and on \(\HHHna(Y_h,\mathbb{Z})\)) induces a representation that we call again \(\varrho_S\colon  W_\aleph\to \Aut(\HHH^*(Y_h,\mathbb{Z}))\) (respectively \(\OGP(\HHHna(Y_h,\mathbb{Z}))\) by restriction).
\end{proposition}

\begin{proof}
    As already noticed, we have a \(W_\aleph\)-covering \(S_0\to S\) since the action of \(W_\aleph\) on \(S_0\) is free. We know by \cite{van1983homotopy} (or \cite[Theorem 1.2]{paris2014k}) that \(\pi_1(S,h) \simeq A_{W_\aleph}\), the Artin group of \(W_\aleph\). Moreover, the short exact sequence of homotopy groups associated to the covering \(S_0\to S\) is exactly the short exact sequence
    \[1\to PA_{W_\aleph}\to A_{W_\aleph}\to W_\aleph \to 1 .\]
   Thus we deduce that \(PA_{W_\aleph} = \pi_1(S_0,h)\). However, \autoref{lemma:action_on_paths} implies that \(\rho_{S_0}(\pi_1(S_0,h)) = \varrho_S(PA_{W_\aleph})\) acts trivially on cohomology. Factorising through \(W_\aleph \simeq A_{W_\aleph}/PA_{W_\aleph}\), we get an induced action \(\varrho_S : W_\aleph \to \Aut(\HHH^*(Y_h,\mathbb{Z}))\).
\end{proof}

Let \(w\in W_\aleph\) and let \(\widetilde{w}\in N_G(T)\subset G\) be a representative of \(w \in W_\aleph = N_G(T)/T\). The torus \(T\) acts trivially on \(S_0\) and because it is connected, it also acts trivially on cohomology by homotopy invariance; thus setting 
\begin{equation}
\label{eq_action_morphism_mon}
    w(\sigma_{\alpha,h}):= [\widetilde{w}(Y_{\alpha,h})]\in \HHH^*(Y_{w(h)},\ZZ)
\end{equation} 
is well-posed since in fact it does not depend on the choice of the representative \(\widetilde{w}\) of \(w\). We now describe the monodromy representation \(\varrho_S\colon  W_\aleph \to \Aut(\HHH^*(Y_h,\mathbb{Z}))\) from \autoref{prop_mon_induces_action_W}. 

\begin{lemma}\label{lemma:action_of_elementary_loops}
Let \(\gamma\colon [0,1]\to  S\) be a closed path with lift \(\widetilde{\gamma}\) to \(S_0\) such that \(\widetilde{\gamma}(1) = w\cdot \widetilde{\gamma}(0)\) for \(w\in W_\aleph\). Then the associated homotopy class \([\gamma]\in \pi_1(S,\gamma(0))\) acts on \(\HHH^{*}(Y_{\gamma(0)},\ZZ)\) by
\[\varrho_S([\gamma])(\sigma_{\alpha,\gamma(0)}) = (w^{-1})(\sigma_{\alpha,\gamma(0)}).\]
\end{lemma}

\begin{proof}
Fix \(w \in W_\aleph\) and let \(\widetilde{w} \in N_G(T)\subset G\) a representative of \(w \in W_\aleph= N_G(T)/T\). Then \(\widetilde{w}\) acts on \(X\) and induces an isomorphism \(Y_{\widetilde{\gamma}(1)} \simeq Y_{w^{-1}\cdot \widetilde{\gamma}(1)} =Y_{\widetilde{\gamma}(0)}\). \autoref{lemma:action_on_paths} implies that if we start with a class \(\sigma_{\alpha,\widetilde{\gamma}(1)}\) and follow the path \(\widetilde{\gamma}\subset  U\), the mentioned class will be \emph{transported} to the class \(\sigma_{\alpha,\widetilde{\gamma}(0)}\). However, in order to get the monodromy action for the base \(S = S_0/W_\aleph\) of our family, we need to identify \(Y_{\widetilde{\gamma}(0)}\) with \(Y_{\widetilde{\gamma}(1)}\) through the action of \(W_\aleph\). Recall that \(\gamma\) is closed so that \(\gamma(0)=\gamma(1)\) and then obviously \(Y_{\gamma(0)}=Y_{\gamma(1)}\) As a result the monodromy action is given by
\[\varrho_S([\gamma])(\sigma_{\alpha,\gamma(0)})=\varrho_S([\gamma])(\sigma_{\alpha,\gamma(1)})=[\gamma]^* (\sigma_{\alpha,\gamma(1)})=[\widetilde{w}^{-1}(Y_{\alpha,\widetilde{\gamma}(0)})] = (w^{-1})(\sigma_{\alpha,\gamma(0)}),\]
where the second last equality follows from \autoref{rem_parallel_and_action} and the last equality from Equation \eqref{eq_action_morphism_mon}.
\end{proof}

It remains then to explicitly compute the action of \(W_\aleph\) on \(\HHHna(Y_h,\mathbb{Z})\). This is given in the alternative basis \(\{\widetilde{\Gamma_\alpha}\coloneqq \epsilon_\alpha \Gamma_\alpha\}_{\alpha\in \Phi_{\aleph}}\) and uses equivariant cohomology in the next subsection.

\subsection{Explicit action of the Weil group}

We compute explicitly the action of \(W_\aleph\) defined above on vanishing cohomology \(\HHHna(Y_h,\mathbb{Z})\), for \(h\in S_0\). We will omit the subscript \(h\) in the following, and set \(Y=Y_h\). Let \(w = s_\alpha \in W_\aleph\) be the reflection around a simple root \(\alpha \in \Phi_\aleph\). Since an additive basis of \(\HHH^{\dim Y}(Y,\mathbb{Z})\) is given by \((\sigma_\gamma)_{\pm\gamma\in \Phi_\aleph}\), our objective is to compute the coefficients \(a_{\beta,\gamma}\) appearing in
\begin{equation}\label{eq:coefficients_a}
    w(\sigma_{-\beta})=\sum_{\pm \gamma \in \Phi_\aleph} a_{\beta,\gamma} \sigma_\gamma,
\end{equation}
for \(\beta\in \Phi_\aleph\) (this will be sufficient to get the whole action). We use \(T\)-equivariant cohomology and recall the few of its basic properties that we need. We refer to \cite{Brion} for more details

For a \(T\)-variety \(Z\), define its  \(T\)-equivariant cohomology by \(\HHH_T^*(Z) := \HHH^*(Z_T)\) where \(Z_T := (Z \times E_T)/T\) with \(E_T \to B_T\) a principal \(T\)-bundle such that \(E_T\) is contractible and \(B_T\) is the classifying space of \(T\). Equivariant cohomology \(\HHH_T^*(Z)\) is a \(\HHH_T^*(\mathrm{pt})\)-algebra and \(\HHH_T^*(\mathrm{pt})\) can be identified with  the polynomial ring \(\QQ[\fX(T)]\) over the character group \(\fX(T)\) of \(T\). 

From now on, we assume that \(Z\) has a finite number of \(T\)-fixed points and \(T\)-stable curves. As explained in Subsection \ref{subsection:BB-decomposition-and-coho-basis}, we have a BB-decomposition for \(Z\) with classes \((\sigma_{z,Z})_{z \in Z^T}\) (here \(\sigma_{z,Z} = \sigma_{z,Z}^+\) with notations as in Subsection \ref{subsection:BB-decomposition-and-coho-basis}) generating \(\HHH(Z)\) as \(\QQ\)-vector space. The \(\HHH_T^*(\pt)\)-module structure of \(\HHH_T^*(Z)\) is very simple and given by the decomposition
\begin{equation*}
     \HHH_T^*(Z) \simeq \bigoplus_{z \in Z^T} \QQ[\fX(T)] \sigma_{z,Z},
\end{equation*}
where here, by abuse of notation, we also denote by \(\sigma_{z,Z}\) the \(T\)-equivariant cohomology class of the closure of positive BB-cells. Moreover, we recover the classical cohomology by:
\begin{equation*}
    \HHH^*(Z) \simeq \HHH_T^*(Z)  / \QQ[\fX(T)]_+,
\end{equation*}
where \(\QQ[\fX(T)]\) is the ideal of polynomial with trivial constant term. We use this property to recover the coefficients \(a_{\beta,\gamma}\) in \eqref{eq:coefficients_a} from the corresponding coefficients in \(T\)-equivariant cohomology. 

The multiplicative structure of \(\HHH_T^*(Z)\) is determined by the so-called \emph{localization at fixed points}, i.e. the pullback morphism
\begin{equation*}
    i^*\colon \HHH_T^*(Z) \to \HHH_T^*(Z^T) \simeq \QQ[\fX(T)]^{Z^T},
\end{equation*}
where $i\colon Z^T \to Z$ is the inclusion. This morphism is injective and one has a control on its image. Moreover \(i^*\) is a ring morphism and the rightmost isomorphism identifies multiplication in \(\HHH_T^*(Z^T)\) with coordinate-wise polynomial multiplication on \(\QQ[\fX(T)]^{Z^T}\).
For \(z \in Z^T\), set \(Z_z := Z_z^+\) with notations as in Subsection \ref{subsection:BB-decomposition-and-coho-basis} and \(f_z := i^*(\sigma_{z,Z})\). Then \(f_z = (f_z(p))_{p \in Z^T}\) satisfies: 
\begin{enumerate}[(i)]
    \item The polynomial \(f_z(p)\) is homogeneous of degree \(\mathrm{codim}(Z_z)\). \label{i}
    \item If \(p \notin Z_z\) then \(f_z(p) = 0\). \label{ii}
    \item If \(p \in  Z_z\) is a smooth point, then \(f_z(p)\) is the product of the \(T\)-characters of \(N_{Z_z/Z,p}\). \label{iii}
    \item For \(p,q \in Z^T\), if there exists a \(T\)-stable curve \(C\) joining \(p\) and \(q\), then \(C \simeq \PP^1\) and the character \(\chi_C\) of the action of \(T\) on \(C\) divides \(f_z(p)- f_z(q)\). \label{iv}
\end{enumerate}
\medskip

We apply the previous setting to \(Y\). For \(w\) and \(\beta\) as above, we consider \(w(\sigma_{-\beta}) = [w(Y_{-\beta})]\) as an element in \(\HHH^*_T(Y)\) via the inclusion \(\HHH^*(Y) \subset \HHH^*_T(Y)\). Denote by \(w(f_{-\beta}) \in \HHH^*_T(Y^T)\) the pullback of \(w(\sigma_\beta) \in \HHH^*_T(Y)\). Properties \ref{i} to \ref{iv} above are satisfied by \(w(f_{-\beta})\) replacing \(Y_{-\beta}\) with \(w(Y_{-\beta})\). We have a \(T\)-equivariant version of \eqref{eq:coefficients_a} above given by
\begin{equation}
\label{eq:coefficients_atilde}
   w(f_{-\beta})=\sum_{\pm \gamma \in \Phi_\aleph} \widetilde{a}_{\beta,\gamma} f_\gamma,
\end{equation}
where the coefficients \(\widetilde{a}_{\beta,\gamma}\) satisfy \(\widetilde{a}_{\beta,\gamma} = a_{\beta,\gamma} \mod (\mathbb{Q}[\fX(T)]_+)\).

Recall that we have an order relation on \(\Delta\) given by \(\gamma \leq \gamma'\) if \(\gamma' - \gamma\) is a non-negative linear combination of simple roots. Recall also that \(Y_{-\beta} = X_{-\beta}\) thus \(x_\gamma \in Y_{-\beta} \Leftrightarrow \gamma \leq -\beta\) and \(x_\gamma \in s_\alpha(Y_{-\beta}) \Leftrightarrow s_\alpha(\gamma) \leq -\beta\). By \cite[Proposition 3.32]{BP}, the point \(x_\gamma \in Y_\gamma\) is smooth for any \(\gamma \in \aleph\) and the set \(\Pi_\gamma\) of characters of the normal bundle \(N_{Y_\gamma/Y,x_\gamma}\) of \(Y_{\gamma}\) at \(x_{\gamma}\) is given by 
\begin{equation}
    \label{eq_characters_normal_bundle}
    \Pi_{\gamma} = \{\eta \in \Delta\mid \eta > 0, \eta\neq -\gamma ,\eta + \gamma \in \Pi(V)\},
\end{equation} 
where \(\Pi(V)\) is the set of \(T\)-weights of \(V\) (see Subsection \ref{subsection:coadjoint} for the definition of \(\Delta\) and \(V\)). Note that in particular we have \(f_\gamma(x_\gamma) \neq 0\) for all \(\gamma \in \aleph\). 

We start the explicit computation of the action on middle cohomology of \(w = s_\alpha\) with \(\alpha \in \Phi_{\aleph}\). 


\begin{lemma}
\label{lemma:beta-not-alpha}
If \(\beta \neq \alpha\), then \(s_\alpha(f_{-\beta}) = f_{-\beta}\).
\end{lemma}

\begin{proof}
We prove that \(Y_{-\beta} = s_\alpha(Y_{-\beta})\), which implies \(s_\alpha(f_{-\beta}) = [s_\alpha(Y_{-\beta})] = [Y_{-\beta}] = f_{-\beta}\). Note that by \cite{buch:rigidity} the converse implication also holds. For \(\lambda \in \Pi(V)\), recall the definition of \(V_\lambda = \{ v \in V \ | \ t.v = \lambda(t) v, \forall t \in T\}\). For \(v \in V\), let \(v_\lambda\) be the \(\lambda\)-component of \(v\), explicitely \(v = \sum_{\lambda \in \Pi(V)}v_\lambda\) with  \(v_\lambda \in V_\lambda\). For \(\gamma \in \Delta_\aleph\) with \(\gamma < 0\), we have \(Y_{\gamma} = X_{\gamma} = \{[v] \in X \ | \ v_\lambda = 0 \text{ for all } \lambda \not\leq \gamma\}\). Since \(s_\alpha(v_\lambda) \in V_{s_\alpha(\lambda)}\) and since \(s_\alpha(\{\lambda \in \Pi(V) \ | \ \lambda \leq -\beta\}) = \{\lambda \in \Pi(V) \ | \ \lambda \leq -\beta\}\), we get \(s_\alpha(Y_{-\beta}) = Y_{-\beta}\).
\end{proof}

\begin{lemma}
\label{lemm:alpha}
We have \(s_\alpha(f_{-\alpha}) = f_\alpha\).
\end{lemma}

\begin{proof}
For \(\gamma \in \Phi_\aleph\), we have \(x_\gamma \in s_\alpha(Y_{-\alpha}) \Leftrightarrow s_\alpha(\gamma) \leq -\alpha \Leftrightarrow \gamma = \alpha\). For \(-\gamma \in \Phi_\aleph\), we have \(x_\gamma \in s_\alpha(Y_{-\alpha}) \Leftrightarrow s_\alpha(\gamma) \leq -\alpha \Leftrightarrow (\gamma \neq -\alpha \text{ and } \gamma \not\perp \alpha)\). By \ref{ii}, 
we deduce that \eqref{eq:coefficients_atilde} has the following form
\[s_\alpha(f_{-\alpha}) = \widetilde{a}_{\alpha,\alpha} f_\alpha + \sum_{\beta\in\Phi_\aleph, \beta \not\perp \alpha} \widetilde{a}_{\alpha,-\beta}  f_{-\beta}.\]
Since \(x_{-\alpha}\not\in s_\alpha(Y_{-\alpha})\) and \(x_{-\alpha}\not\in Y_{\alpha}\) by \cite[Remark 3.29.2.b]{BP}, we have \(s_\alpha(f_{-\alpha})(x_{-\alpha}) = 0 = f_\alpha(x_{-\alpha})\). The same holds for any simple root \(\beta\neq \alpha\): \(f_{-\beta}(x_{-\alpha}) = 0\), while \(f_{-\alpha}(x_{-\alpha}) \neq0\), this implies \(\widetilde{a}_{\alpha,-\alpha} = 0\). 
Since \(x_{-\alpha} \in Y_{-\alpha}\) is a smooth point, so is \(x_\alpha \in s_\alpha(Y_{-\alpha})\). By \ref{iii}, \(s_\alpha(f_{-\alpha})(x_\alpha)\) is the image by \(s_\alpha\) of the product of characters of \(N_{Y_{-\alpha}/Y,x_{-\alpha}}\). By  
\eqref{eq_characters_normal_bundle}, we deduce
\[ s_\alpha(f_{-\alpha})(x_\alpha) = \prod_{\eta>0, \eta\neq \alpha, \eta-\alpha\in\Pi(V)} s_\alpha(\eta) = \prod_{\delta>0, \delta \neq -\alpha, \delta + \alpha \in \Pi(V)} \delta = f_\alpha(x_\alpha), \]
thus \(\widetilde{a}_{\alpha,-\alpha} =1\). We are left to compute the coefficients \(\widetilde{a}_{\alpha,-\beta}\) for  \(\beta \in \Phi_\aleph\) with \(\beta \neq \alpha\) and \(\beta \not\perp \alpha\). Note that \(x_{-\alpha-\beta}\notin s_\alpha(Y_{-\alpha})\) since \(x_{-\beta} \notin Y_{-\alpha}\), so \(s_\alpha(f_{-\alpha})(x_{-\alpha-\beta}) = 0\). By the (equivariant) Chevalley formula in \cite{BP}, we have \(f_{-\beta}(x_{-\alpha-\beta})\neq 0\), while \(f_{-\gamma}(x_{-\alpha-\beta}) = 0\) for \(\gamma \in \Phi_\aleph \setminus \{\alpha,\beta\}\) (because \(x_{-\alpha -\beta} \not\in Y_{-\gamma}\)). Finally \(f_\alpha(x_{-\alpha-\beta}) = 0\) by \ref{ii} together with \cite[Proposition 3.23]{BP} telling us that \(x_{-\alpha-\beta}\notin Y_{\alpha}\). Evaluating \eqref{eq:coefficients_atilde} at \(x_{-\alpha-\beta}\) gives \(\widetilde{a}_{\alpha,-\beta} = 0\).
\end{proof}


\begin{proposition}
\label{prop:action-of-W}
    We have the following formula in \(\HHH^*(Y)\):
    \begin{equation*}
s_\alpha(\Gamma_\beta) = \left\{
\begin{array}{cl}
\Gamma_\beta & \text{ if } \beta\perp \alpha, \\
 \Gamma_\beta-\Gamma_\alpha & \text{ if } \beta\not\perp\alpha \text{ and } \beta \neq \alpha \\
-\Gamma_\alpha & \text{ if } \beta = \alpha.
\end{array}\right.
\end{equation*}
\end{proposition}

\begin{proof}
    We have a basis of the middle cohomology of \(Y\) given by \((j^*{\sigma}_{\gamma,X}, \sigma_{-\gamma})_{\gamma\in \Phi_{\aleph}}\). By \autoref{lemma:beta-not-alpha} and \autoref{lemm:alpha}, we know the action of \(s_\alpha\) on \(\sigma_{-\gamma}\): it is the same as the one on the \(f_{-\gamma}\). Since this action comes from an action of monodromy, which is trivial on the ambient cohomology, we have \(s_{\alpha}(j^*\sigma_{\gamma,X}) = j^*\sigma_{\gamma,X}\). The result follows using \autoref{cor-j} and linear algebra.
 \end{proof}


\begin{proposition}
\label{prop:correspondence_of_action}
    The action of \(W_\aleph\) on \(\HHHna(Y)\) identifies with the action of \(W_\aleph\) on \(\oplus_{\alpha \in \Phi_\aleph} \QQ \alpha\) via the identification \(\widetilde{\Gamma_\alpha}\leftrightarrow \alpha\).
\end{proposition}

\begin{proof}
If \(\beta \perp \alpha\), then \(s_\alpha(\widetilde{\Gamma}_\beta) = s_\alpha(\epsilon_\beta\Gamma_\beta) = \epsilon_\beta s_\alpha(\Gamma_\beta)=\epsilon_\beta\Gamma_\beta=\widetilde{\Gamma}_\beta\). If \(\beta \not\perp \alpha\) and \(\alpha \neq \beta\), then \(\epsilon_\alpha = -\epsilon_\beta\) so that
\(s_\alpha(\widetilde{\Gamma}_\beta) = s_\alpha(\epsilon_\beta\Gamma_\beta) = \epsilon_\beta s_\alpha(\Gamma_\beta) = \epsilon_\beta \Gamma_\beta-\epsilon_\beta\Gamma_\alpha = \epsilon_\beta \Gamma_\beta+\epsilon_\alpha\Gamma_\alpha=\widetilde{\Gamma}_\beta+\widetilde{\Gamma}_\alpha\). Finally \(s_\alpha(\widetilde{\Gamma}_\alpha) = -\epsilon_\alpha\Gamma_\alpha = -\widetilde{\Gamma}_\alpha\). This is exactly the desired action.
\end{proof}

\subsection{Lower bound on the monodromy group}

We gather the previous results. Recall that \(A_{W_\aleph}\) is the Artin group associated to \(W_\aleph\) and \(PA_{W_\aleph}\leq A_{W_\aleph}\) is the kernel of the natural projection \(A_{W_\aleph}\to W_\aleph\). By \autoref{prop:correspondence_of_action} we can identify \(\HHHna(Y,\mathbb{Z})\) with \(\Lambda_\aleph\), and \(\HHHna(Y,\mathbb{C})\) with the standard representation of \(W_\aleph\) given by \(\sum_{\alpha\in \Delta_\aleph} \mathbb{C} H_\alpha \subset V_0\simeq V_0^\vee\). 


\begin{theorem}\label{thm:main_theorem_monodromy}
The monodromy representation \(\varrho_S : \pi_1(S,h) \simeq A_{W_\aleph} \to \mathrm{Aut}(\HHHna(Y,\mathbb{Z}))\) factors through the composition of the natural projection \(p\colon A_{W_\aleph}\to W_\aleph\) and the standard representation \(\rho_\aleph\colon W_\aleph \to \OGP(\HHHna(Y,\mathbb{Z})) = \OGP(\Lambda_\aleph)\) realizing \(W_\aleph\) as a reflection group. In particular, we have embeddings of groups
\begin{equation*}
W_\aleph \simeq \Mon(Y_h, B)\leq \Mon(Y_h)\leq W_\aleph\rtimes\Aut(\tau_\aleph).
\end{equation*}
\end{theorem}

\begin{proof}
Combine Corollary \ref{cor:isometry_group_middle_coh}, Proposition \ref{prop_mon_induces_action_W} with its proof, and Proposition \ref{prop:correspondence_of_action}.
\end{proof}

\begin{corollary}\label{cor:monodromy_no_exterior}
    If \(D_\aleph\) has no automorphisms,  
     then \(\Mon(Y_h) \simeq W_\aleph\). This happens in the adjoint cases for types \(C_n\), \(G_2\), \(E_7\), \(E_8\) and the coadjoint (non-adjoint) cases for types \(B_n\), \(G_2\).
\end{corollary}

\begin{corollary}\label{cor:monodromy_simplylaced}
    If \(X\) is adjoint simply-laced, then \(\Mon(Y_h) \simeq W_\aleph\).
\end{corollary}

\begin{proof}
Recall \autoref{exam:reg-ss-adjoint}: in the adjoint simply-laced case, we have  \(W_\aleph = W\) and the universal family of smooth hyperplane sections \(U\subset V\) coincides with the set of regular elements, \emph{i.e.} \(U = G.U_0 = G.S_0\). Moreover the restriction of the GIT quotient map \(V \to V\git G\) to \(U\) is a morphism \(U\to U\git G\) where each fiber is a single closed orbit (of a regular semisimple element). By the Chevalley Theorem \(U \git G\simeq S_0/W = S\). Each semisimple (regular) orbit is simply connected by \cite[Corollary 2.3.4]{coll-mcgov}. Therefore, applying the long exact sequence of homotopy groups for the fibration \(U\to S\) we deduce that \(\pi_1(U,h) \simeq \pi_1(S,h)\), and this isomorphism identifies \(\rho\) with \(\varrho_S\). We get \(\Mon(Y_h) = \Mon(Y_h,U) = \rho(\pi_1(U,h)) = \varrho_S(\pi_1(S,h)) = \Mon(Y_h,S) \simeq W.\)
\end{proof}

\begin{remark}
Similar results about monodromy have been proved in \cite[Lemma 4.3]{belmans2025ade}.
\end{remark}

\section{Reduction of semisimplicity to ambient cohomology}
\label{sec:quantum}

We use the previous result on the monodromy action to prove that, in some cases, semisimplicity of the ambient quantum cohomology implies semisimplicity of the full quantum cohomology.

\subsection{Quantum cohomology}
\label{subsection:qh}

Quantum cohomology is defined for any smooth projective variety  \(Z\) over \(\CC\). To simplify the exposition, we assume that \(Z\) is Fano with Picard rank~\(1\) and refer to \cite{manin} for more details. Let \(i_Z \in \ZZ_{>0}\) be the index of \(Z\) defined by \(-K_Z = i_Z h_Z\) where \(h_Z\) is the ample generator of \(\Pic(Z)\). Fix a homogeneous basis \((\sigma_{i,Z})_{i \in [0,N]}\) of \(\HHH^*(Z) := \HHH^*(Z, \QQ)\) such that \(\sigma_{0,Z} = 1\) and \(\sigma_{1,Z} = h_Z\). Fix dual linear coordinates \((t_i)_{i \in [0,N]}\) and a formal variable \(Q\). Set \(\deg(t_i) = 2(1 - \deg(\sigma_{i,Z}))\) and \(\deg(Q) = 2i_Z\). The underlying space of the big quantum cohomology ring of \(Z\) is the \(\QQ[Q]\formal{{t_0, \cdots , t_N}}\)-module \(\BQH(Z) := \HHH^*(Z) \otimes_\QQ \QQ[Q]\formal{t_0,\cdots,t_u}\) (the fact that we can define \(\BQH(Z)\) over \(\QQ[Q]\) comes from the Fano assumption on \(Z\)).

For cohomology classes \(\sigma,\tau \in \HHH^*(Z)\), the product in \(\BQH(Y)\) is of the form 
\begin{equation*}
\sigma \star \tau := \sigma \cup \tau + \sum_{i = 0}^N \sum_{d\geq 1}\sum_{\underline{k}}\frac{Q^dt_{\uk}}{\uk!} \langle \sigma,\tau,\sigma_{i,Z},\sigma_{\uk,Z}\rangle_{d,|\uk| + 3}^Z  \sigma_{i,Z}^\vee,
\end{equation*}
where \(\uk = (k_0,\cdots,k_N)\) is a sequence of non-negative integers, \(\uk! = \prod_{i = 0}^N k_i!\), \(|\uk| = \sum_{i = 0}^N k_i\) and \(t_\uk = \prod_{i = 0}^N t_i^{k_i}\). The symbol \(\langle \ \bullet \ \rangle_{d,n}^Z\) represents the Gromov-Witten invariant of degree \(d\) with \(n\) insertions (we will often write \(\langle \ \bullet \ \rangle_{d,n}\) or \(\langle \ \bullet \ \rangle_{d}\) for \(\langle \ \bullet \ \rangle_{d,n}^Z\)). Finally \(\sigma_{\uk,Z}\) means the insertion, for all \(i \in [0,N]\), of \(k_i\)-classes \(\sigma_{i,Z}\). Note that the divisor axiom implies that the above formula can be rewritten in the following classical form
\begin{equation*}
\sigma \star \tau = \sigma \cup \tau + \sum_{i = 0}^N \sum_{d\geq 1}\sum_{\underline{k}}\frac{q^dt_{\uk}}{\uk!} \langle \sigma,\tau,\sigma_{i,Z},\sigma_{\uk,Z}\rangle_{d,|\uk| + 3}  \sigma_{i,Z}^\vee,
\end{equation*}
where the sum is over all \(\uk\) with \(k_1 = 0\) and \(q = Qe^{t_1}\). We will use this last presentation.

There are several specialisation of the above quantum product. Setting \(q\) and all \((t_i)_{i \in [0,N]}\) to zero gives back the cup product while setting all \((t_i)_{i \in [0,N]}\) to zero defines the small quantum cohomology \(\QH(Z)\) which is a \(\QQ[q]\)-module with product denoted by \(\sigma \star_0 \tau\). Explicitly, we have
\begin{equation*}
\sigma \star_0 \tau := \sigma \cup \tau + \sum_{i = 0}^N \sum_{d\geq 1}\langle \sigma,\tau,\sigma_{i,Z} \rangle_{d,3}  q^d \sigma_{i,Z}^\vee.
\end{equation*}
Note that \(\QH(Z)=\BQH(Z)/(t_0,\dots,t_N)\). More generally, for a cohomology class \(\eta \in \HHH^*(Z)\), the product in the deformation \(\BQH_\eta(Z)\) along the direction \(\QQ\eta \subset \HHH^*(Z)\) is of the form 
\begin{equation*}
    \sigma \star_\eta \tau := \sigma \star_0 \tau + \sum_{i = 0}^N\sum_{d \geq 1} \sum_{k \geq 0} \langle \sigma ,\tau ,\sigma_{i,Z} , \underbrace{\eta,\cdots,\eta}_{k \textrm{ terms}} \rangle_{d ,k+3} q^d \frac{\mathtt{t}^k}{k!} \sigma_{i,Z}
\end{equation*}
where \(\mathtt{t}\) is the deformation parameter. Note that \(\BQH_\eta(Z)=\BQH(Z)/(\Ker_\eta)\) where \(\Ker_\eta\) is the ideal generated by elements \(f\) in the \(\QQ\)-span of \((t_i)_{i \in [0,N]}\) with \(f(\eta) = 0\).

Finally, it will be useful to consider specification of these rings at \(q=1\), for which we set \(\BQHo(Z) := \BQH(Z)/(q-1)\), \(\QHo(Z) := \QH(Z)/(q-1)\) and  \(\BQHo_\eta(Z) := \BQH_\eta(Z)/(q-1)\).

\medskip

Denote by \(\xi\) and \(\xi_\eta\) the generic points of the affine schemes \(\Spec\left(\QQ[q]\formal{t_0,\dots,t_N}\right)\) and \newline \(\Spec\left(\QQ[q]\formal{t_0,\dots,t_N}/(\Ker_\eta)\right) \simeq \Spec \left( \QQ[q]\formal{\mathtt{t}}\right)\) respectively. Consider the rings \newline
\(\BQH(Z)_\xi := \BQH(Z)\otimes_{\QQ[q]\formal{t_0,\cdots,t_N}} \QQ(q)(\!(t_0,\cdots,t_N)\!)\) and  
\(\BQH_\eta(Z)_\xi := \BQH_\eta(Z)\otimes_{\QQ[q]\formal{\mathtt{t}}} \QQ(q)(\!(\mathtt{t})\!)\).  

\begin{definition}\label{def:generical_semisimplicity}
We will say that \(\BQH(Z)\) (respectively \(\BQH_\eta(Z)\)) is \emph{generically semisimple} if \((\BQH(Z))_\xi\) (resp. \((\BQH_\eta(Z))_{\xi_\eta}\)) is a semisimple algebra.\footnote{Recall that a finite dimensional algebra \(A\) over a field \(\textsf{k}\) is semisimple if it is a product of fields. This is equivalent to the fact that the scheme \(\Spec(A)\) is reduced, that the morphism \(\Spec(A) \to \Spec(\textsf{k})\) is smooth, or that the radical of \(A\) vanishes.}
\end{definition}

\begin{remark}
Note that the generic semisimplicity of \(\BQH_\eta(Z)\) implies that of \(\BQH(Z)\).
\end{remark}

\subsection{Reduction to ambient cohomology}

Let \(j\colon Z \hookrightarrow W\) be a smooth Fano hyperplane section of a smooth Fano variety. Recall the definition of ambient and vanishing cohomology from \autoref{def:vanish-coho}: \(\HHHa(Z) := \im(j^*\mid_{\HHH^*(W)})\) and \(\HHHna(Z) := \Ker(j_*\mid_{\HHH^*(Z)})\). Recall from \cite[Theorem 15.4]{voisin2002theorie} that the vanishing cohomology is an irreducible representation for the monodromy group. The following result is well known. 

\begin{lemma}
\label{rmk:one_primitive_insertion}
    Let \((a_i)_{i \in [1,k]} \in \HHHa(Z)\) and \(b\in \HHHna(Z)\). Then \(\langle a_1,\cdots,a_k,b\rangle_d = 0\).
\end{lemma}

\begin{proof}
 We have a morphism \(u \colon \HHHna(Z) \to \CC\), \(b\mapsto \langle a_1,\cdots,a_k,b\rangle_d\). GW-invariants are invariant under the action of the monodromy group of \(Z\subset W\), so \(u\) should be a monodromy invariant linear function on \(\HHHna(Z)\). However, \(\HHHna(Z)\) is an irreducible representation of the monodromy group; thus \(u=0\), and this holds for any \(k\), any \((a_i)_{i \in [1,k]} \in \HHHa(Z)\) and any \(d\). 
 \end{proof}

Let \(\eta \in \HHHa(Z)\) be an ambient class, let \(\BQH^{\amb}_\eta(Z)\) be the 
\(\QQ[q]\formal{\mathtt{t}}\)-submodule of \(\BQH_\eta(Z)\) generated by \(\HHHa(Z)\) and let \(\QH^{\amb}(Y)\) be its specialisation for \(\eta = 0\).

\begin{corollary}
The submodule \(\BQH^{\amb}_\eta(Z)\) of \(\BQH_\eta(Z)\) is stable for the product \(\star_\eta\).
\end{corollary}

For \(\eta \in \HHHa(Z)\) an ambient class, \(\BQH^{\amb}_\eta(Z)\) is thus a subalgebra of \(\BQH_\eta(Z)\). The following result gives conditions to deduce the semisimplicity of \(\BQH_\eta(Z)\) from that of \(\BQH_\eta^{\amb}(Z)\).

\begin{lemma}
\label{lem_zero_products}
Assume that any product \(b_1\star_\eta b_2\) of vanishing classes is contained in \(\BQH^{\amb}_\eta(Z)\). Then if $\BQH^{\amb}_\eta(Z)$ is generically semisimple, so is $\BQH_\eta(Z)$.
\end{lemma}

\begin{proof}
 Let \(\chi \in \BQH_\eta(Z)\) be a non-zero nilpotent element of square zero \(\chi \star_\eta \chi = 0\). Decompose \(\chi = \chi_a + \chi_v\) in its ambient \(\chi_a\) and vanishing \(\chi_v\) parts. Since \(\BQH_\eta^{\amb}(Z)\) is generically semisimple, we have \(\chi_v \neq 0\). Moreover, since GW-invariants with one vanishing insertion are zero by \autoref{rmk:one_primitive_insertion}, we have \(\chi_a \star_\eta \chi_a \in \BQH_\eta^{\amb}(Z)\), while \(\chi_v \star_\eta \chi_a\) is a \(\QQ[q]\formal{\mathtt{t}}\)-linear combination of vanishing classes. Since any product of two vanishing classes is contained in the ambient cohomology, we have that \(\chi_v \star_\eta \chi_v \in \BQH_\eta^{\amb}(Z)\). Since \(0 = \chi\star_\eta \chi = \chi_a \star_\eta \chi_a + 2\chi_a \star_\eta \chi_v + \chi_v \star_\eta \chi_v\) this implies that \(\chi_a \star_\eta \chi_v = 0\). But then \(0 = \chi \star_\eta \chi \star_\eta \chi_a = \chi_a\star_\eta \chi_a \star_\eta \chi_a = 0\), so \(\chi_a\) is a nilpotent ambient class and therefore \(\chi_a = 0\). We deduce \(\chi = \chi_v\) and also \(\chi_v \star_\eta \chi_v = 0\). However this implies, working modulo \(\mathtt{t}\) and \(q\), that \(\chi_v \cup \chi_v = 0\), which is a contradiction since the intersection form on the vanishing cohomology is non-degenerate.
\end{proof}

\begin{remark}
\label{rem_GW_three_prim}
The condition that all products \(b_1\star_\eta b_2\) of vanishing classes is contained in \(\BQH^{\amb}_\eta(Z)\) is equivalent to the fact that all GW-invariants involving three vanishing classes and any number of copies of \(\eta\) vanish.
\end{remark}

We apply this result to the case of hyperplane sections of an adjoint or coadjoint variety. Let \(Y \subset X\) be a smooth hyperplane section of an adjoint variety or a coadjoint variety \(X\) with \(\Pic(Y) \simeq \ZZ\).

\begin{theorem}
\label{thm_QBH_nonAn}
 Let \(Y\) a smooth hyperplane section of \(X\) for simply-laced \(G\) not of type \(A\). Let \(\eta \in \HHHa(Y)\) be an ambient class and assume that \(G\) is simply laced.
If \(\BQH_\eta^{\mathrm{amb}}(Y)\) is generically semisimple, so is \(\BQH_\eta(Y)\).
\end{theorem}

\begin{proof}
We want to apply \autoref{lem_zero_products} and \autoref{rem_GW_three_prim}. It is enough to prove that, for \(b_1\), \(b_2\), \(b_3\), vanishing classes, we have \(\langle b_1,b_2,b_3,\eta,\cdots,\eta\rangle_d = 0\) for any number of occurrences of \(\eta\). Consider the morphism 
  \begin{equation*}
       u\colon \Sym^3 \HHHna(Y) \to \CC \quad , \quad (b_1,b_2,b_3)\mapsto \langle b_1,b_2,b_3,\eta,\cdots,\eta\rangle_d,
  \end{equation*}
  where we have fixed the number of occurrences of \(\eta\) and the degree \(d\). This function \(u \in \Sym^3 \HHHna(Y)^\vee\) is invariant for the monodromy action. Now, by \autoref{thm:main_theorem_monodromy}, we know that the monodromy group \(\Mon(Y)\) contains \(W_\aleph = W\), which acts on \(\HHHna(Y)\) with its standard representation. 
The algebra of invariants for the standard representation of \(W\) coincides with the algebra of \(G\)-invariants \(\CC[\fg^\vee]^{G}\) by Chevalley Theorem. It is a polynomial algebra generated in degrees \(e_1+1,\cdots,e_r+1\), where the \(e_i\)'s are the exponents of \(G\). Since \(G\) is not of type \(\text{A}_m\) by hypothesis, it contains no invariant of degree three (\(e_1=1\) and all other exponents are \(\geq 3\)). Therefore \(u\) must be zero, and this should happen for any \(d\) and any number of occurrences of \(\eta\). Thus we can apply \autoref{lem_zero_products} and deduce that \(\BQH_\eta(Y)\) is generically semisimple.
\end{proof}

\begin{remark}
Let \(X\) be an adjoint variety and \(Y\) a general hyperplane section. Generic semisimplicity of \(\BQH(Y)\) was proved in \cite{BP} for non-simply laced groups and for groups of type \(\text{A}\). 
\end{remark}
 
If \(X\) is coadjoint not adjoint, then \(G\) is not simply laced, $W_\aleph$ is of type $A_n$ and the above result does not apply (there are always invariants of degree $3$). Note that for type \(\text{B}\) and type \(\text{G}\) the coadjoint varieties are smooth quadrics and the same holds for \(Y\) for which the semisimplicity always holds \cite{perrin2014semisimple}. We are left with two cases in types \(\text{C}\) and \(\text{F}\). We recall the dimension \(\dim(Y)\) and index \(i_Y\) in both cases.

\[\begin{array}{c|cc}
\text{Type}(G) & \dim(Y) & i_Y \\
\hline
\text{C}_n & 2(2n - 3) & 2n - 2 \\
\text{F}_4 & 14 & 10 \\
\end{array}\]

\begin{proposition}
\label{prop:QBH_nonAn}
Assume that \(G\) is not simply laced and $X$ coadjoint.
Let \(a \in \ZZ_{\geq 0}\) be a non trivial divisor of \(i_Y\) and let \(\eta \in \HHH^{2\ell}(Y)^{\amb}\) be an ambient class of degree \(2\ell\), such that \(a\) divides \(\ell - 1\).

If \(\BQH_\eta^{\mathrm{amb}}(Y)\) is generically semisimple, so is \(\BQH_\eta(Y)\).
\end{proposition}

\begin{proof}
As in the previous proof it is enough to prove that, for \(b_1\), \(b_2\), \(b_3\) vanishing classes, we have  
\[\langle b_1,b_2,b_3,\underbrace{\eta,\cdots,\eta}_{k \textrm{ terms}}\rangle_{d,k+3} = 0 \]
 for any \(d\) and any \(k\). This invariant vanishes unless \(\sum_{i = 1}^3\deg(b_i) + 2k\ell = 2\dim(Y) + 2k+ 2di_Y\). Since \(\dim(Y) = 2N\) and \(\deg(b_i) = N\) for \(i \in [1,3]\), we get the desired vanishing unless \(N = 2di_Y + 2(1 - \ell )k\) and thus unless \(a\) divides \(N\). This is not possible. 
\end{proof}

We expect that \autoref{thm_QBH_nonAn} and \autoref{prop:QBH_nonAn} will lead to generic semisimplicity of \(\BQH(Y)\) as conjectured in \cite{BP}. In what follows we will focus on the non simply laced case where the above technique does not fully apply and on the coadjoint case where the generic sesmisimplicity is not known. In Section \ref{sec:small_quantum_cohomology} we compute a presentation of \(\QH(Y)\) in that case and deduce a presentation of \(\QH^{\amb}(Y)\). In Section \ref{sec:computations_for_IG(2,6)},  we rather consider a deformation in a vanishing cohomology direction \(\eta \in \HHHna(Y)\) and deduce generic semisimplicity in that case.

\section{Presentation of the small quantum cohomology 
}\label{sec:small_quantum_cohomology}

For \(Z\) Fano  of Picard rank one, according to \cite[Conjecture 1.3]{KS}, under the assumption that \(\BQH(Z)\) is generically semisimple, the algebraic structure of \(\QH(Z)\), in particular its radical, should give information on the structure of \(\text{D}^b(Z)\) the bounded derived category of \(Z\). It is therefore desirable to compute a presentation of \(\QH(Z)\). We will fix the quantum deformation parameter as \(q = 1\). Recall that for any variety \(Z\) we have that \(\QHo(Z) := \QH(Z)/(q-1)\) is a \(\ZZ/2i_Z\ZZ\)-graded ring, where \(i_Z\) is the index of the variety, whose graded pieces are \(\QHo(Z)^k \simeq \bigoplus_j \HHH^{k+2ji_Y}(Z)\) as vector spaces.

In this section we give a presentation of \(\QH(Y)\) for \(Y\) a hyperplane section of a coadjoint variety \(X\) for a non simply laced group \(G\) \emph{i.e.} for types \(\text{C}_n\) and \(\text{F}_4\).   We will be mostly interested in the structure of the radical of \(\QHo(Y)\). Since we only work with the small quantum cohomology we will often replace the product \(\star_0\) by the standard product notation, explicitely we write \(\sigma\sigma'\) for \(\sigma \star_0 \sigma'\) and \(\sigma^k\) for the $k$-th power of \(\sigma\) in the \(\star_0\) product.

We will often use the fact that quantum multiplication by the hyperplane class is of maximal rank so the dimensions of the graded pieces of \(\QHo(Y)\) play an important role in what follows. We listed these dimensions in \autoref{table-Cn} and \autoref{table-F4} .

\begin{table}[ht]
\[\begin{array}{c|ccccccccc}
 k & 0 & 1 & 2 & 3 & \cdots & 2n - 6 & 2n -5 & 2n - 4 & 2n - 3 \\
 \hline
  \\[-10pt]
\dim(\QHo^{2k}(Y)) & n & n - 1 & n & n - 1 & \cdots & n & n - 1 & n & 2(n-1) 
\end{array}\]
     \centering
   \caption{\label{table-Cn} Dimensions of the graded pieces of \(\QHo(Y)\) in the \(\text{C}_n\) case.}
\end{table}

\begin{table}[ht]
\[\begin{array}{c|cccccccccc}
 k & 0 & 1 & 2 & 3 & 4 & 5 & 6 & 7 & 8 & 9 \\
 \hline
 \\[-10pt]
\dim(\QHo^{2k}(Y)) & 3 & 2 & 2 & 2 & 3 & 2 & 2 & 4 & 2 & 2 
\end{array}\] 
     \centering
   \caption{\label{table-F4} Dimensions of the graded pieces of \(\QHo(Y)\) in the \(\text{F}_4\) case.}
\end{table}




\subsection{Reminders on \(\QH(X)\) and \(\QH(Y)\)}

In this subsection we recall some results from \cite{BP} useful for the computation of the presentation of \(\QH(Y)\). Our main focus will be to compare results in \(X\) with results in \(Y\) especially because \(\QH(X)\) is well understood from \cite{maxim}. Recall the definition of the Schubert classes in \(X\) and \(Y\) from Subsection \ref{subsection:BB-decomposition-and-coho-basis}. The following two results are valid for \(X\) adjoint or coadjoint.

\begin{proposition}[{\cite[Cor. 4.15]{BP}}]
  \label{cor-j}
  Let \(\alpha \in \aleph\).
  \begin{enumerate}
    \item For \(\alpha\) positive, non simple, we have \(j^*\sigma_{\alpha,X} = \sigma_\alpha\).
    \item For \(\alpha\) negative, we have \(j^*\sigma_{\alpha,X} = \sum_{\beta \neq \alpha} a_{\alpha}^{\beta}\sigma_{\beta}\) and \(j_*\sigma_\alpha = \sigma_{\alpha,X}\).
    \item For \(\alpha\) simple, we have \({j^*\sigma_{\alpha,X} = \sigma_\alpha + \sigma_{-\alpha} + \sum_{\beta\in \Phi_\aleph,\ \alpha + \beta \in \Delta_\cdot} \sigma_{-\beta}}\) and \(j_*\sigma_\alpha = \sigma_{-\alpha,X}\).
  \end{enumerate}
  \end{proposition}

\begin{proposition}[{\cite[Proposition 5.13]{BP}}]
  \label{prop:comp-GW}
  Let \(\sigma_Y,\sigma'_Y \in \HHH^*(Y,\ZZ)\) and \(\gamma_X \in \HHH^*(X,\ZZ)\) with \(\deg(\gamma_X) < i_Y = i_X - 1\). 
  Then we have \(\langle \sigma_{Y} , \sigma'_{Y} , j^*\gamma_X \rangle_1^Y = \langle j_*\sigma_{Y} , j_*\sigma'_{Y} , \gamma_X \rangle_1^X\).
\end{proposition}

%
%

We now assume that \(X\) is coadjoint not adjoint. Recall the quantum Chevalley formula for \(X\).

\begin{theorem}[Quantum Chevalley formula for \(X\), \cite{fulton-woodward,ChPe}]
\label{thm-qch-hom-qmin}
Let \(h_X \in \HHH^2(X)\) be the hyperplane class and let \(\alpha \in \aleph\). We have
\begin{equation*}
    h_X \star_0 \sigma_{\alpha,X} = h_X \cup \sigma_{\alpha,X} + q \delta_{\alpha \leq \theta - \Theta} \sigma_{\alpha + \Theta}.
\end{equation*}

\end{theorem}

\begin{theorem}[Quantum Chevalley formula for \(Y\),  {\cite[Theorem 5.38]{BP}}]
  \label{thm-qch-qm}
    Let \(h \in \HHH^2(Y)\) be the hyperplane class and let \(\alpha \in \aleph\).
  \begin{enumerate}
  \item If \(\alpha > 0\) is not simple, then \(h \star_0 \sigma_\alpha = h \cup \sigma_\alpha\).
  \item If \(\alpha\) is simple, then \(h \star_0 \sigma_\alpha = h \star_0 \sigma_{-\alpha}\).
  \item If \(\alpha < 0\), then \(h \star_0 \sigma_\alpha = h \cup \sigma_\alpha + \delta_{\alpha \leq \theta - \Theta} q j^*\sigma_{\alpha + \Theta,X}\).
  \end{enumerate}
\end{theorem}

\begin{corollary}[{\cite[Corollary 5.20]{BP}}] 
\label{5.20}
We have \(h \star_0 \Gamma_\alpha = 0\) for all \(\alpha \in \Phi_\aleph\).
\end{corollary}

Consider the endomorphisms \(E_X = (h_X \star -)\) on \(\QH(X)\) and \(E_Y = (h \star -)\) on \(\QH(Y)\). For \(k \in [0,c_1(X) - 1]\), we set \(\QH^k(X) = \oplus_{i \geq 0}\HHH^{k+2i c_1(X)}(X)\) and for \(k \in [0,c_1(Y) - 1]\), we set \(\QH^k(Y) = \oplus_{i \geq 0} \HHH^{k + 2ic_1(Y)}(X)\). Note that the map \(j_* : \HHH^{\dim Y}(Y) \to \HHH^{\dim Y + 2}(X)\) is surjective and define a right inverse \((j_*)^{-1} : \HHH^{\dim Y + 2}(X) \to \HHH^{\dim Y}(Y)\) by \((j_*)^{-1}(\sigma_{-\alpha,X}) = \sigma_{-\alpha}\) for \(\alpha \in \Phi_\aleph\). Define a \(\QQQ[q]\)-linear map \(\mathbb{j} : \QH(X) \to \QH(Y)\) via \(\mathbb{j}(q) = q\) (by abuse of notation \(q\) is the quantum parameter for both \(X\) and \(Y\)) and 
\begin{equation*}
    \mathbb{j}\vert_{\HHH^{2k}(X)} = \left\{\begin{array}{ll}
j^* & \textrm{ for \(k \in [0,\dim Y/2]\)} \\
(j_*)^{-1} & \textrm{ for \(k \in [\dim Y/2+1,\dim X]\)} \\
\end{array}\right.
\end{equation*}

\begin{lemma}
  The map \(\mathbb{j} : \QH(X) \to \QH(Y)\) is an isomorphism of \(\QQ\)-vector spaces.
\end{lemma}

Define the map \(\bar E_X : \QH(X) \to \QH(X)\) by

\begin{equation*}
    \bar E_X \vert_{\HHH^{2k}(X)} =
\left\{\begin{array}{ll}
E_X & \textrm{ for \(k \neq \dim Y/2\)} \\
E_X^2 & \textrm{ for \(k = \dim Y/2\)}. \\
  \end{array}\right.
\end{equation*}
  
  \begin{corollary}[{\cite[Corollary 5.42]{BP}}]
  \label{5.42}
We have \(\Ker E_Y = \mathbb{j}(\Ker \bar E_X)\).
\end{corollary}

\begin{proposition}[{\cite[Proposition 5.44]{BP}}]
\label{5.44}
  Let \(X\) be coadjoint not adjoint.
  \begin{enumerate}
  \item The radical of \(\QH(X)\) is contained in \(\Ker E_X\).
  \item The non-zero eigenvalues of \(E_X\) have multiplicity one.
\item The minimal polynomial \(\mu_{E_X}\) of \(E_X\) is of the form \(\mu_{E_X}(T) = T P(T^{c_1(X)})\) with \(P\) having non-zero simple roots.
  \end{enumerate}
\end{proposition}

\begin{corollary}[{\cite[Corollary 5.45]{BP}}]
\label{5.45}
  The minimal polynomial \(\mu_{E_Y}\) of \(E_Y\) is \(\mu_{E_Y} = T P(T^{c_1(Y)})\) and the non-zero eigenvalues of \(E_Y\) have multiplicity one. In particular, the radical of \(\QH(Y)\) is contained in \(\Ker E_Y\). 
\end{corollary}

By \cite{maxim}, the rank of \(\Ker E_X\) is equal to \(s + 1\) where \(s\) is the number of short roots \(s\) of \(G\); moreover, the Dynkin diagram of short roots of \(G\) describes completely the radical of \(\QH(X)\) and \(\Ker E_X\). Indeed, the radical of \(\QH(X)\) is isomorphic to \(\QQ[\tau_X]/(\tau_X^2)\) if \(X = F_4/P_4\) and to \(\QQ[\tau_X]/(\tau_X^{n-1})\) if \(X = \IGr(2,2n)\). Here \(\tau_X\) is of degree \(8 + 2i_X = 30\) in type \(\text{F}_4\) and degree \(4 + 2i_X = 4n + 2\) in type \(\text{C}_n\). Denote also by \(\tau_X\) the corresponding element modulo \(q - 1\). By \autoref{5.44}, we get the following presentation:

\begin{proposition}
\label{prop_QHoX}
  Let \(X\) be coadjoint not adjoint. Then \begin{equation*}
      \QHo(X)\simeq \QQ[h_X, \tau_X]/(h_XP(h_X^{i_X}),\tau_X^s, h_X\tau_X)
  \end{equation*}
   where \(P\) is a polynomial of degree \(s\) with non-zero constant coefficient.
\end{proposition}

\subsection{Presentations of \(\QH(Y)\) and \(\QH^{\amb}(Y)\)}
\label{subsection:presentation}

Let \(h \in \HHH^2(Y)\) be the hyperplane class. Set \(\tau:=\mathbb{j}\tau_X\), we have \(\deg(\tau) = 28\) in type \(\text{F}_4\) and \(\deg(\tau) = 4n\) in type \(\text{C}_n\). Recall that \((\Gamma_\alpha)_{\alpha\in \Phi_\aleph}\) is a basis of \(\HHHna(Y)\). 

\begin{theorem}
\label{thm_pres_smallQHbar}
Let \(X\) be coadjoint not adjoint. Then \begin{align*}
    \QHo(Y)\simeq \QQ[h, \tau, (\Gamma_\alpha)_{\alpha\in \Phi_\aleph}]/\!\left(hP(h^{i_Y}),\tau^s, h\tau, \{h\Gamma_\alpha\}_{\alpha\in \Phi_\aleph},\{\Gamma^2_\alpha - \Gamma_\beta^2\}_{\alpha\neq \beta\in \Phi_\aleph}\right.,\\ 
    \nonumber\left.\{\Gamma^2_\alpha - 2\Gamma_\alpha \Gamma_\beta\}_{\alpha\neq \beta\in \Phi_\aleph}, \{\tau \Gamma_\alpha\}_{\alpha\in \Phi_\aleph}, \{\Gamma^2_\alpha-\tau^{s-1}\}_{\alpha\in \Phi_\aleph}\right).
\end{align*}
\end{theorem}

We will prove this result in several steps. 

\begin{lemma}
\label{lem_tau_ker}
 We have \(\tau \in \Ker E_Y\).
\end{lemma}

\begin{proof}
Note that we have a decomposition \(\tau_X = q a_2^X + b_2^X\), where \(a_2^X \in \HHH^{\deg(\tau_X)-2i_X}(X)\) and \(b_2^X \in \HHH^{\deg(\tau_X)}(X)\). Moreover we have \(\dim(X) \not\in \{ \deg(\tau_X) , \deg(\tau_X) - 2 i_X\}\), 
which implies that \(\overline{E}_X \tau_X = E_X \tau_X=0\). By \autoref{5.42}, we have \(\tau = \mathbb{j}\tau_X \in \Ker E_Y\). 
\end{proof}

We prove that \(\Ker E_Y \cap \QH^{\deg{\tau}}(Y) = \QQ\tau\). 
If \((\alpha_i)_{i \in [1,r]}\) are the simple roots of \(G\) with notation as in \cite{bourbaki} and if \(\alpha \in \Delta\) is of the form \(\alpha = \sum_{i = 1}^r m_i \alpha_i\), then we will replace \(\sigma_\alpha\) by \(\sigma_{m_1\cdots m_r}\). 

\begin{lemma}
\label{lem_tau_expr}
We have \(\dim(\Ker E_Y \cap \QH^{\deg{\tau}}(Y)) = 1\). In particular, up to scalar, we have:
\begin{enumerate}
    \item If \(X = F_4/P_4\), then \(\tau = \sigma_{-\theta} + q(-\sigma_{1111} + \sigma_{0121}).\)
    \item If $X=\IGr(2,2n)$, then 
    \begin{align*}
        \tau = \left( -\sigma_{-11110\cdots 0} + \sigma_{-011110\cdots 0} + \cdots + (-1)^{n-3}\sigma_{-0\cdots 01111} + (-1)^{n-2}\sigma_{-0\cdots 00121}\right) + 
    \\
    \nonumber
    + q \left(\sigma_{1112\cdots 21} - \sigma_{012\cdots 21}\right).
    \end{align*}
\end{enumerate}
\end{lemma}

\begin{proof}
Let \(\tau \in \Ker E_Y\) of degree \(28\) in the \(\text{F}_4\)-case and \(4n\) in the \(\text{C}_n\)-case. Write \(\tau = \tau' + q \tau''\), with \(\tau' \in \HHH^{\deg(\tau)}(Y)\) and \(\tau'' \in \HHH^{\deg(\tau) - 2i_Y}(Y)\). We have \(0 = \tau \star_0 h = \tau' \cup h + q( \tau'' \cup h + \sum_\gamma \langle \tau',h,\sigma_\gamma \rangle_1 \sigma^\vee_\gamma)\) (for degree reasons there is no \(q^2\) term), thus \(\tau'\cup h = 0\) and \(\tau'' \cup h + \sum_\gamma \langle \tau',h,\sigma_\gamma \rangle_1 \sigma^\vee_\gamma = 0\).
  
From Chevalley formula in \(X\) and \autoref{cor-j}, there exists a unique (modulo scalar) element \(\tau' \in \HHH^{\deg(\tau)}(Y)\) such that \(\tau' \cup h = 0\). This element is \(\sigma_{-\theta}= [\pt]\) in the \(\text{F}_4\) case and \(-\sigma_{-11110\cdots 0} + \sigma_{-011110\cdots 0} + \cdots +(-1)^{n-3}\sigma_{-0\cdots 01111} + (-1)^{n-2}\sigma_{-0\cdots 00121}\) in the \(\text{C}_n\) case.
  
  In the \(\text{F}_4\)-case, from \autoref{thm-qch-hom-qmin}, we get \(\tau'  \star_0 h = q\sigma_{1110}\), so we need to find \(\tau''\) such that \(\tau''\cup h = -\sigma_{1110}\), and this is given by \(-\sigma_{1111} + \sigma_{0121}\) again by the Chevalley formula. In the \(\text{C}_n\)-case, we get \(\tau'\star_0 h = -q\sigma_{11112\cdots 21}\), so we need to find \(\tau''\) such that \(\tau'' \cup h = \sigma_{11112\cdots 21}\), and this is given by \(\sigma_{1112\cdots 21} - \sigma_{012\cdots 21}\).
\end{proof}

\begin{proposition}
\label{prop_nonzeropower}
 We have \(\tau^{s - 1} \neq 0\) in \(\QH(Y)\) 
 \end{proposition}

\begin{proof}
  In the \(\text{F}_4\)-case \(s = 2\), so there is nothing to prove. Assume that \(X = \IGr(2,2n)\). As already mentioned, 
  \(\tau_X\) has degree \(4n +2 = 4 + 2 i_X\) and \(\tau\) has degree \(4n = 4 + 2 i_Y\). Set \(\tau_2^X := \tau_X\) and \(\tau_2 := \tau\). By \autoref{prop_QHoX}, \(\tau_{2i}^X := (\tau_2^X)^i\) is a non-zero element in \(\Ker E_X \subset \QH(X)\) for \(i = 1,\cdots,s-1\). Denote by \(\tau_{2i} := \mathbb{j}\tau_{2i}^X\).   Since \(\mathbb{j}\) is an isomorphism of vector spaces, if we show that \(\tau_{2i} = \tau_2^i\) for \(i = 1,\cdots,s-1\), then since \(\tau_{2i} \neq 0\), we have proved our claim. In fact, it is sufficient to proceed by induction, showing that \(\tau_{2i+2} = \tau_{2i} \star_0 \tau_2\). Otherwise stated, our statement is proven if we show that 
 \begin{equation}
 \label{eq_prop_nonzeropower}
     \mathbb{j}(\tau_{2i}^X \star_0 \tau_{2}^X)=\mathbb{j}\tau_{2i}^X \star_0 \mathbb{j}\tau_2^X.
 \end{equation}
 Fix \(i \in [1,s-2]\). For degree reasons, we always have a decomposition \(\tau_{2i}^X = q^{i-1}(qa_{2i}^X + b_{2i}^X)\) with \(\deg(a_{2i}^X) = 4i \leq \dim(Y)\) and \(\deg(b_{2i}^X) = 4i + 2 i_X \geq \dim(Y)+2\) (here we are essentially using the fact that \(2i_Y>\dim(Y)\)). This also means that \(\tau_{2i} = q^{i-1}(qa_{2i} + b_{2i})\), with \(a_{2i} = j^* a_{2i}^X\) and \(b_{2i}^X = j_* b_{2i}\). A similar decomposition holds for the RHS of \eqref{eq_prop_nonzeropower}: \(\mathbb{j}\tau_{2i}^X \star_0 \mathbb{j}\tau_2^X = q^{i}(qc_{2i+2} + d_{2i+2})\).
 
 Since cup product by \(h\) has maximal rank, and looking at the Betti numbers of \(Y\), there can exist at most one element in \(\Ker E_Y\) of degree \(\deg(\tau_{2k})\) for \(k = 1,\cdots,s-1\). Moreover both \(\mathbb{j}\tau_{2i+2}^X = \mathbb{j}(\tau_{2i}^X \star_0 \tau_2^X)\) and \(\mathbb{j} \tau_{2i}^X \star_0 \mathbb{j}\tau_2^X\) belong to \(\Ker E_Y\) and have the same degree, so they are multiple one of the other. To prove that they are equal it is enough to show that \(b_{2i+2} = d_{2i+2}\). 
 
 Let us start from the LHS of \eqref{eq_prop_nonzeropower}, \emph{i.e.} from \(b_{2i+2}\). We have 
 \[\begin{array}{ll}
b_{2i+2} & = \displaystyle{\mathbb{j}(a_{2i}^X\cup b_2^X + b_{2i}^X\cup a_2^X + q \sum_\lambda \langle b_{2i}^X,b_2^X,\sigma_\lambda^X \rangle_1(\sigma_\lambda^X)^\vee )} \\
  & = \displaystyle{\mathbb{j}(a_{2i}^X\cup j_* b_2 + j_* b_{2i}\cup a_2^X + q \sum_\lambda \langle j_* b_{2i},j_*b_2,\sigma_\lambda^X \rangle_1(\sigma_\lambda^X)^\vee ).}
  \end{array}\]
 By the projection formula for \(j\) 
 and  \autoref{prop:comp-GW}, together with the definition of \(\mathbb{j}\), we have 
 \[\begin{array}{ll}
b_{2i+2} & = \displaystyle{\mathbb{j}\left(j_*(j^*a_{2i}^X\cup b_2 + b_{2i}\cup j^*a_2^X) + q \sum_\lambda \langle  b_{2i},b_2,j^*\sigma_\lambda^X \rangle_1(\sigma_\lambda^X)^\vee \right)} \\
 & = \displaystyle{a_{2i}\cup b_2 + b_{2i}\cup a_2 + q \sum_\lambda \langle  b_{2i},b_2,j^*\sigma_\lambda^X \rangle_1 \mathbb{j}\left((\sigma_\lambda^X)^\vee\right).}
   \end{array}\]
   The last expression is equal to \(d_{2i+2}\) provided we show that \(\mathbb{j}\left((\sigma_\lambda^X)^\vee\right) = (j^* \sigma_\lambda^X)^\vee\). Tracking degrees one realizes that \(\mathbb{j}\left((\sigma_\lambda^X)^\vee\right) = (j_*)^{-1}((\sigma_\lambda^X)^\vee)\). Applying \(j_*\) we are reduced to show that \(j_* (j^*( \sigma_\lambda^X)^\vee)=(\sigma_\lambda^X)^\vee\). Consider a class \(\sigma_\mu^X\) of the right degree. Since \(j_*[\pt]=[\pt]\), projection formula implies the following equalities, \begin{equation*}
     \sigma_\mu^X \cup j_*((j^* \sigma_\lambda^X)^\vee)=j_*(j^*\sigma_\mu^X \cup (j^*\sigma_\lambda^X)^\vee)=\delta_{\mu,\lambda}j_*[\pt]=\delta_{\mu,\lambda}[\pt],
 \end{equation*} 
 where \(\delta_{\bullet,\bullet}\) is Kronecker's delta, proving the claim.
\end{proof}

Looking at the Betti numbers of \(Y\) one checks, as in \autoref{lem_tau_expr}, that there exists, up to scalar, a unique  element \(\tau_0\) in \(\Ker E_Y\) with \(\deg(\tau_0) = 2i_Y\) which is easy to compute explicitly. Indeed, up to scalar we have \(\tau_0 = q + \tau_0'\), with \(\tau_0' \in \HHH^{2i_Y}(Y)\). We may compute \(\tau_0 \star_0 h = \tau_0'\cup h + qh(1+\langle \tau_0',h,h^\vee \rangle_1)\), so \(\tau_0'\) is defined by the condition that \(\tau_0' \cup h = 0\) and \(\langle \tau_0',h,h^\vee \rangle_1 = -1\). Arguing as in \autoref{lem_tau_expr} again, one can check that in the \(\text{F}_4\) case \(\tau_0 = q - \sigma_{-1111} + \sigma_{-0121}\), and in the \(\text{C}_n\) case \(\tau_0 = q-\sigma_{-110\cdots 0} + \sigma_{-110\cdots 0} - \cdots +(-1)^{n-1}\sigma_{-0\cdots 011}\).

\begin{proposition}
\label{prop_mult_radical}
  Inside \(\QHo(Y)\), \(\Ker E_Y\) is generated multiplicatively by \(\tau_0\), \(\tau\), and \((\Gamma_\alpha)_{\alpha\in \Phi_\aleph}\).
\end{proposition}

\begin{proof}
The elements \(\tau_0\), \((\Gamma_\alpha)_{\alpha\in \Phi_\aleph}\), \(\tau, \tau^2, \cdots, \tau^{s-1}\) are linearly independent (they have pairwise distinct degrees). One can check that \(i_Y\deg(P) + 2s = \rank \QHo(Y) = \rank \HHH^*(Y)$. This shows that \(\rank(\Ker E_Y) = 2s\), and we have already found \(2s\) linearly independent elements.
\end{proof}

\begin{proposition}
\label{prop_rel_QHo}
For all \(\alpha, \beta\in \Phi_\aleph\), we have the following relations in \(\QHo(Y)\).
  \begin{enumerate}
\item \(\Gamma^2_\alpha - \Gamma_\beta^2 = \Gamma^2_\alpha - 2\Gamma_\alpha \star_0 \Gamma_\beta = 0\) 
\item \(\tau \star_0 \Gamma_\alpha = 0\).
\item \(\tau^s=0\).
\item There exists a non-zero rational number \(d\) such that \(\Gamma^2_\alpha-d\tau^{s-1} = 0\) for all \(\alpha\in \Phi_\aleph\).
  \end{enumerate}
\end{proposition}

\begin{proof}
Recall that \(\QHo(Y)\) is \(\ZZ/2i_Y\ZZ\)-graded.

1. In degree \(2 \dim(Y)\), the element \(\tau^{s-1}\) is the unique, up to scalar, element in \(\Ker E_Y\). Thus \(\Gamma_\alpha^2\) and \(\Gamma_\alpha\Gamma_\beta\) are multiples of \(\tau^{s-1}\) for any \(\alpha,\beta \in \Phi_\aleph\). Furthermore \(\tau^{s-1} = \lambda[\pt] + qu\), for some \(\lambda \in \QQ\) and \(u \in \HHH^{2\dim(Y) - 2i_Y}(Y)\). Note that since \(\Gamma_\alpha \cup \Gamma_\alpha \neq 0\), we must have \(\lambda \neq 0\) and 
\[\Gamma_\alpha^2 = \frac{\Gamma_\alpha \cup \Gamma_\alpha}{\lambda} \tau^{s-1} \text{and similarly } \Gamma_\alpha \star_0 \Gamma_\beta = \frac{\Gamma_\alpha \cup \Gamma_\beta}{\lambda} \tau^{s-1} .\] 
The result follows from \autoref{thm_inter_product_middle_cohom}. Item 4 also follows from this.

2. In degree \(\deg(\Gamma_\alpha) + \deg(\tau)\), there is no nonzero element in \(\Ker E_Y\), so \(\tau\Gamma_\alpha\) must vanish.

3. By items 2 and 4, we have \(\tau^s = \tau^{s-1} \star_0 \tau = \frac{1}{d}\Gamma_\alpha^2 \star_0 \tau = 0\).
\end{proof}

\begin{proof}[Proof of  \autoref{thm_pres_smallQHbar}]
  All the relations on the RHS are satisfied by \autoref{5.45},  \autoref{5.20},  \autoref{lem_tau_ker}, and \autoref{prop_rel_QHo}. Since \(\tau_0 = 1+a\), with \(a\) in the image of \(E_Y\), by \autoref{prop_mult_radical} the elements \(h\), \(\tau\) and \((\Gamma_\alpha)_{\alpha\in \Phi_\aleph}\) generate \(\QHo(Y)\). So we have a well-defined surjective map from the RHS to \(\QHo(Y)\). A straightforward computation shows that \(\QHo(Y)\) and the RHS have the same rank, so they must be isomorphic.
\end{proof}

\begin{corollary}
\label{cor_ambient_cohom}
We have a ring isomorphism 
\begin{equation*}
    \QHoa(Y)\simeq \QQ[h, \tau ]/(hP(h^{i_Y}),\tau^s, h\tau).
\end{equation*} 
\end{corollary}

\begin{corollary}
\label{cor_nilp_even_degree}
When \(X = \IGr(2,2n)\), all nilpotents in \(\QHa(Y)\) have even degree. Moreover the set of degree \(2n - 3\) nilpotents coincides with the vanishing cohomology.
\end{corollary}

\begin{proof}
  By \autoref{cor_ambient_cohom}, the radical in \(\QHa(Y)\) is generated by \(\tau\) modulo \(q - 1\), so nilpotent elements have degree \(2n , 2n+2 , \dots , 4n-6\) modulo \(\deg(q) = 2n-2\). In particular there are no nilpotents in degree \(2n - 3\) in \(\QH^{\amb}(Y)\). This, together with the fact that by \autoref{thm_pres_smallQHbar} the vanishing cohomology is contained in the radical of \(\QH(Y)\), implies the last claim.
\end{proof}

We end this section with a different presentation of \(\QH^{\amb}(Y)\). We focus on hyperplanes sections of \(X = \IGr(2,6)\) this is the only case we will consider in what follows. Let \(h_X,\sigma_X\) be the Chern classes of the dual tautological bundle \(\cU^\vee\) of rank \(2\) on \(X\) and set \(h := j^*h_X\) and \(\sigma := j^*\sigma_X\).

\begin{proposition}
\label{prop_pres_small_quantum_Y}
Let \(Y\) be a general hyperplane section of  \(X = \IGr(2,6)\). Then we have 
\begin{equation*}
\QH^{\amb}(Y)=\QQ[h,\sigma,q]/(R_4(q),R_5(q)),
\end{equation*} 
with \(R_4(q) = 3\sigma^2 + h^4 - 4\sigma \star_0 h^2\) and \(R_5(q) = -5h^5 + 14\sigma  \star_0 h^3 + 9qh\).
\end{proposition}

\begin{proof}
Note that \(\sigma_X = \sigma_{011,X}\) and \(\sigma = \sigma_{011}\). 
 By \cite[Proposition 3.2]{ig22n}, we have a presentation \(\HHH(X) = \QQ[h_X,\sigma_X]/(R_4^X,R_6^X)\), with \(R_4^X = 3\sigma_X^2 + h^4 - 4\sigma_X h_X^2\) and \(R_6^X = -5h_X^6 + 14\sigma_X h_X^4\). We get a surjective map \(j^* : \QQ[h_X,\sigma_X]/(R^X_{4},R^X_{6}) \to \HHHa(Y)\). From the Betti numbers, there are relations \(R_4\) and \(R_5\) in degree \(4\) and \(5\) between \(h\) and \(\sigma\) in \(\HHHa(Y)\) and by \cite[Lemma 4.40]{BP}, we must have \(R_4 = 3\sigma_X^2 + h^4 - 4\sigma_Xh_X^2\) and \(R_5 = - 5h^5 + 14\sigma_X h^3\) where \(R_4\) is obtained from \(R_4^X\) replacing \(h_X,\sigma_X\) by \(h,\sigma\) and \(R_5\) is obtained in the same way from \(R_5^X\) with \(R_6^X = h_XR_5^X\). This induces a presentation \(\HHH^{\amb}(Y) = \QQ[h,\sigma]/(R_4,R_5)\).

We now deform the above relations \(R_4,R_5\) in \(\QH^{\amb}(Y)\). To finish the proof, we only need to prove that the following relations hold in \(\QH(Y)\):
 \begin{equation*}
 3\sigma^2+h^4-4\sigma \star_0 h^2=0 \text{ and } -5h^5+14\sigma \star_0 h^3=-9qh.
 \end{equation*}
To compute these products in \(\QH(Y)\) it is sufficient to use the quantum Chevalley formula (\autoref{thm-qch-qm}) and the product \(\sigma\star_0 \sigma\). The latter only involves the Gromov-Witten invariant \(\langle \sigma,\sigma,\pt\rangle_1\). This is the number of lines in \(Y\) passing through general translates of \(\sigma,\sigma,\pt\). Let \(V_6\) be a \(6\)-dimensional symplectic vector space such that \(X = \IGr(2,V_6)\). Fix a general point \([U_2] \in Y\) and two general translates of \(\sigma\), which are given by two general \(5\)-dimensional subspaces \(U_5,U_5'\) of \(V_6\) (the class of \(\sigma\) is the class of the zero locus of a general section of \(\cU^\vee\) on \(Y\)). Lines in \(Y\) are contained in lines in \(\Gr(2,V_6)\), which are parametrized by \(F(1,3;V_6)\). Consider a line \(l\), given by a flag \(F_1\subset F_3\), passing through \([U_2]\) and intersecting \(\Gr(2,U_5) \cap Y\) and \(\Gr(2,U_5')\cap Y\). The first condition implies that \(F_1\subset U_2\subset F_3\). Let \([V_2] \in l\cap \Gr(2,V_5)\cap Y\); then \(V_2 \neq U_2\) which means that \(V_2\cap U_2 = F_1 = U_5\cap U_2\). Similarly for \(U_5'\), which gives that \(U_5 \cap U_2 = F_1 = U_5'\cap U_2\), which is not possible by genericity of \(U_2,V_5,V_5'\). Therefore \(\langle \sigma,\sigma,\pt\rangle_1=0\). Then the quantum Chevalley formula gives the claimed relations.
\end{proof}

\begin{corollary}
The only fat point in \(\QHo^{\amb}(Y)\) is the one defined by the ideal \(\cM = (h,\sigma)\).
\end{corollary}

\begin{proof}
  This follows from \autoref{cor_ambient_cohom} since \(\tau = \sigma^3 - 2h^2 q + 3\sigma q\).
\end{proof}

\begin{remark}
\label{rem:pres-HHH-coadjoint}
A similar argument gives presentations of \(\HHHa(Y)\) in type \(\text{C}_n\) in general and type \(\text{F}_4\). In \(\text{C}_n\), keeping \(h\) and \(\sigma\) as before, we get a presentation \(\HHHa(Y) = \QQ[h,\sigma]/(R_{2n-2},R_{2n-1})\), with \(R_{2n-2},R_{2n-1}\) relations in degrees \(2n-2,2n-1\). In \(\text{F}_4\), set \(\sigma := \sigma_{0121} = j^*\sigma_{0121,X}\) and keep \(h\) for the hyperplane class, we get a presentation \(\HHH^{\amb}(Y) = \QQ[h,\sigma]/(S_{8},S_{11})\), with \(S_8 = 2 h^8 -6 h^4\sigma + 3 \sigma^2\), \(S_{11} = -11 h^{11} + 26 h^7\sigma\).
\end{remark}

\begin{proposition}
\label{prop_pres_small_quantum_Y}
Let \(Y\) be a general hyperplane section of  \(X = \text{F}_4/\text{P}_4\). Then we have 
\begin{equation*}
\QH^{\amb}(Y)=\QQ[h,\sigma,q]/(S_8(q),S_{11}(q)),
\end{equation*} 
with \(S_8(q) = 2 h^8 -6 h^4\sigma+3 \sigma^2\) and \(S_{11}(q) = -11 h^{11} +26 h^7\sigma+3qh\).
\end{proposition}

\begin{proof}
  The fact that the quantum cohomology has a presentation of the form $\QQ[h,\sigma,q]/(S_8(q),S_{11}(q))$, where $S_{8}(q)$ and $S_{11}(q)$ are deformation of $S_{8}$ and $S_{11}$ is a consequence of the fact that the quantum cohomology is a deformation of the classical cohomology and \autoref{rem:pres-HHH-coadjoint}. Since $q$ has degree $10$, in order to compute these relations in $\QH(Y)$ it is sufficient to use the quantum Chevalley formula (\autoref{thm-qch-qm}).
\end{proof}

\begin{remark}
\label{rem_3_points}
The arguments in \autoref{prop_pres_small_quantum_Y}. also apply for \(X=\IGr(2,6)\) leading to a well known presentation (see for example \cite[Corollary 4.2]{ig22n}) that we explicit for later use
  \begin{equation*}
      \QH(X)=\QQ[h,\sigma_X,q]/(R_4^X(q),R_6^X(q)),
  \end{equation*} 
  where \(R_4^X(q)=3\sigma^2+h^4-4\sigma \star_0 h^2\) and \(R_6^X(q)=-5h^6+14\sigma  \star_0 h^4+9qh\). Using the quantum Chevalley formula for \(X\), one can express each Schubert class in \(\QH(X)\) in terms of \(h,\sigma,q\). Using the above presentation, one can recover any \(3\)-points Gromov-Witten invariant in \(X\). Similarly, the quantum Chevalley formula for \(Y\) and \autoref{prop_pres_small_quantum_Y} recover any \(3\)-points Gromov-Witten invariant involving ambient classes in \(Y\). We implemented this strategy on a computer program, through which we were able to compute \(3\)-points GW-invariants in \(X\) and \(Y\) (in the latter case, only involving ambient classes) that we will need in Section \ref{sec:computations_for_IG(2,6)}.
\end{remark}

\section{Dubrovin's conjecture for hyperplane sections of \(\IGr(2,6)\)}
\label{sec:computations_for_IG(2,6)}

The results of the previous sections imply that to prove generic semisimplicity of \(\BQH(Y)\), it is sufficient to obtain semisemplicity of a deformation \(\BQH_\eta(Y)^{\amb}\) of \(\QH^{\amb}(Y)\) along an ambient class \(\eta\). However, we could not find such an ambient class giving the desired semisimplicity (at least, for small deformations). We still hope to find an appropriate direction in future work but we will follow in the rest of this work a different approach and deform along vanishing directions.

\medskip 

In this section we complete this program for a general hyperplane section \(Y\subset X = \IGr(2,6)\). We start with general simplification results for computing \(\BQH(Y)\) in Subsection \ref{subsec:red-of-gw} and specialise to the special case of \(\IGr(2,6)\) in later subsections.

\subsection{Reduction of GW-invariants}
\label{subsec:red-of-gw}

Let \(Z\) be a smooth Fano variety of Picard rank one and \(h\) be an ample generator of \(\Pic(Z)\). 
Fix \(d\) and \(n\). For \((\alpha_i)_{i \in [1,n]}\) and \((\alpha'_i)_{i \in [1,n]}\) cohomology classes, we write 
\begin{equation*}
    \langle \alpha_1, \cdots, \alpha_n \rangle_d \equiv \langle \alpha_1', \cdots , \alpha_n' \rangle_d + o(d,n)
\end{equation*}
if one can recover \(\langle \alpha_1, \cdots, \alpha_n \rangle_d\) from \(\langle \alpha_1', \cdots , \alpha_n' \rangle_d\) and \(n'\)-points GW-invariants of degree \(d'\) with either \(d'<d\) and \(n' \leq n\) or \(d' \leq d\) and \(n' < n\). If \((\beta_i)_{i \in [1,n]}\) is a sequence of classes and \(\uk = \{k_1 < \cdots < k_r\} \subset [1,n]\), we set \(|\uk| := r\) and write \(\beta_\uk\) for the sequence \(\beta_{k_1},\cdots,\beta_{k_r}\). Inclusion of sequences induces a partial order \(\uk' \leq \uk\). The following result is certainly a well known consequence of the divisor axiom.

\begin{proposition} 
Let \(d,n \in \ZZ_{\geq0}\) and \(\alpha_1,\alpha_2,\alpha_3,\beta_\uk\) cohomology classes with \(|\uk| = n\). We have
\begin{equation*}
    \langle \alpha_1\cup h,\alpha_2,\alpha_3,\beta_{\underline{k}} \rangle_d \equiv \langle \alpha_1,\alpha_2\cup h,\alpha_3,\beta_{\underline{k}} \rangle_d +o(d,n+3).
\end{equation*}
\end{proposition}

\begin{proof}
The LHS is the coefficient of \(\frac{q^d t_{\underline{k}}}{\underline{k}!}\alpha_3^\vee\) in the quantum product \((\alpha_1\cup h )\star \alpha_2\). We have
\begin{equation*}
    \alpha_1 \star h =\alpha_1\cup h + \sum_\gamma \sum_{1\leq d'\leq d}  \sum_{\underline{k}'\leq \underline{k}} \frac{q^{d'}t_{\underline{k}'}}{\underline{k}'!}\langle \alpha_1,h,\gamma,\beta_{\underline{k}'}\rangle_{d'} \gamma^\vee +\mbox{ h.d.t.}
\end{equation*}
here \emph{h.d.t.} means \emph{higher degree terms} and refers to coefficients of \(\frac{q^{d'} t_{\underline{k}'}}{\underline{k}!}\) with either \(d' > d\) or \(\underline{k}' \not \leq \underline{k}\), which will play no role in our argument. Similarly
\begin{equation*}
    \alpha_2 \star h =\alpha_2\cup h + \sum_\gamma \sum_{1\leq d'\leq d}  \sum_{\underline{k}'\leq \underline{k}} \frac{q^{d'}t_{\underline{k}'}}{\underline{k}'!}\langle \alpha_2,h,\gamma,\beta_{\underline{k}'}\rangle_{d'} \gamma^\vee +\mbox{ h.d.t.}
\end{equation*}
We get
\begin{equation*}
 \begin{array}{rl}
 (\alpha_1\cup h )\star \alpha_2 = & \displaystyle{(\alpha_1\star h )\star \alpha_2 - \sum_\gamma \sum_{1\leq d'\leq d}  \sum_{\underline{k}'\leq \underline{k}} \frac{q^{d'}t_{\underline{k}'}}{\underline{k}'!}\langle \alpha_1,h,\gamma,\beta_{\underline{k}'}\rangle_{d'} (\gamma^\vee \star \alpha_2) +\mbox{ h.d.t.} } \\
= & \displaystyle{ \alpha_1\star (h \star \alpha_2) - \sum_\gamma \sum_{1\leq d'\leq d}  \sum_{\underline{k}'\leq \underline{k}} \frac{q^{d'}t_{\underline{k}'}}{\underline{k}'!}\langle \alpha_1,h,\gamma,\beta_{\underline{k}'}\rangle_{d'} (\gamma^\vee \star \alpha_2) +\mbox{ h.d.t.}}\\
= & \displaystyle{  \alpha_1\star (\alpha_2\cup h) +\sum_\gamma \sum_{1\leq d'\leq d}  \sum_{\underline{k}'\leq \underline{k}} \frac{q^{d'}t_{\underline{k}'}}{\underline{k}'!}\langle \alpha_2,h,\gamma,\beta_{\underline{k}'}\rangle_{d'} (\alpha_1\star\gamma^\vee) }\\
 &   \displaystyle{- \sum_\gamma \sum_{1\leq d'\leq d}  \sum_{\underline{k}'\leq \underline{k}} \frac{q^{d'}t_{\underline{k}'}}{\underline{k}'!}\langle \alpha_1,h,\gamma,\beta_{\underline{k}'}\rangle_{d'} (\gamma^\vee \star \alpha_2) +\mbox{ h.d.t.} }\\
\end{array}
\end{equation*}
where again \emph{h.d.t.} will not influence our reasoning. Indeed, in the above expression, we are interested in the coefficient of \(\frac{q^d t_{\underline{k}}}{\underline{k}!}\alpha_3^\vee\). Using the last expression we have found for \((\alpha_1\cup h )\star \alpha_2\), this coefficient is given by the sum of three factors:
\begin{align*}
\langle \alpha_1,\alpha_2\cup h ,\alpha_3,\beta_{\underline{k}}\rangle_d +  \sum_\gamma \sum_{1\leq d'\leq d}  \sum_{\underline{k}'\leq \underline{k}} \frac{\underline{k}!}{\underline{k}'!(\underline{k}-\underline{k}')!}\langle \alpha_2,h,\gamma,\beta_{\underline{k}'}\rangle_{d'} \langle \alpha_1,\gamma^\vee, \alpha_3,\beta_{\underline{k}-\underline{k}'}\rangle_{d-d'} +\\
\nonumber- \sum_\gamma \sum_{1\leq d'\leq d}  \sum_{\underline{k}'\leq \underline{k}}\frac{\underline{k}!}{\underline{k}'!(\underline{k}-\underline{k}')!}\langle \alpha_1,h,\gamma,\beta_{\underline{k}'}\rangle_{d'} \langle \gamma^\vee,\alpha_2, \alpha_3,\beta_{\underline{k}-\underline{k}'}\rangle_{d-d'}.
\end{align*}
In the second and third term above, the GW-invariants \(\langle \alpha_2,h,\gamma,\beta_{\underline{k}'}\rangle_{d'}\) and \(\langle \alpha_1,h,\gamma,\beta_{\underline{k}'}\rangle_{d'}\) are \(o(d,n+3)\) by the divisor axiom, and \(\langle \gamma^\vee,\alpha_2, \alpha_3,\beta_{\underline{k}-\underline{k}'}\rangle_{d-d'}\) and \(\langle \gamma^\vee,\alpha_2, \alpha_3,\beta_{\underline{k}-\underline{k}'}\rangle_{d-d'}\) are \(o(d,n+3)\) since \(d'\geq 1\) hence \(d - d' < d\). We deduce that the coefficient above, which coincides with \(\langle \alpha_1\cup h,\alpha_2,\alpha_3,\beta_{\underline{k}} \rangle_d\), can be recovered from \(\langle \alpha_1,\alpha_2\cup h ,\alpha_3,\beta_{\underline{k}}\rangle_d\) and \(n'\)-points GW-invariants of degree \(d'\) with either \(d' < d\) and \(n' \leq n+3\) or \(d'\leq d\) and \(n'<n+3\).
\end{proof}

From now on we assume that \(X = \IGr(2,6)\) and \(Y \subset X\) is a general hyperplane section. Set \(\Gamma_i := \Gamma_{\alpha_i} \in \HHHna(Y,\mathbb{Z})\), for \(i=1,2\). We have seen in \autoref{prop_pres_small_quantum_Y} that \(\HHHa(Y)\) is generated by the classes \(h\)  and \(\sigma\) and the presentation shows that the classes \(1,\sigma,\Gamma_1,\Gamma_2\) are the only classes that are not \(\cup\)-multiples of \(h\) in \(\HHH(Y)\).

\begin{proposition}
\label{prop_only_some_invariants}
The following invariants determine \(\BQH(Y)\).
\begin{itemize}
    \item \(\langle \sigma,\dots,\sigma,\Gamma_{i_1},\dots,\Gamma_{i_k} \rangle_d\);
    \item \(\langle \sigma,\dots,\sigma \rangle_d\);
    \item \(\langle \alpha,\sigma,\dots,\sigma \rangle_d\), where \(\alpha = \sigma\cup h^j\) but \(\alpha\) is not a \(\cup\)-power of \(h\).
\end{itemize}
\end{proposition}

\begin{proof}
We use an inductive argument. Assume that we know how to compute all \(n'\)-points degree \(d'\) GW-invariants for \(n'<n\) and \(d' \leq d\) or \(n'\leq n\) and \(d'<d\). Consider any GM-invariant \(\langle \alpha_1,\dots,\alpha_n\rangle_d\). If any of the \(\alpha_i\)'s for \(i\geq 2\) is a \(\cup\)-multiple of \(h\), say \(\alpha_2 = h\cup \beta_2\), we can use the previous proposition and deduce that \(    \langle \alpha_1,\dots,\alpha_n\rangle_d\equiv \langle \alpha_1\cup h,\beta_2,\alpha_3,\dots,\alpha_n\rangle_d +o(d,n)\).

We can repeat this procedure until \(    \langle \alpha_1,\dots,\alpha_n\rangle_d\equiv \langle \alpha_1\cup h^j,\beta_2,\dots,\beta_n\rangle_d +o(d,n)
\)
with \(\beta_2,\dots,\beta_n\) not \(\cup\)-multiples of \(h\). Since we can find a basis of the cohomology of \(Y\) in which the only elements which are not \(\cup\)-multiples of \(h\) are \(1,\sigma,\Gamma_1,\dots,\Gamma_s\), we have shown that all \(n\)-point degree \(d\) GW-invariants can be recovered from \(n'\)-points degree \(d'\) GW-invariants for \(n'<n\) and \(d'\leq d\) or \(n'\leq n\) and \(d'<d\), and from the GW-invariants \(\langle \sigma,\dots,\sigma,\Gamma_{i_1},\dots,\Gamma_{i_k} \rangle_d\), \(\langle \Gamma_{i_1},\dots,\Gamma_{i_n} \rangle_d\), \(\langle \sigma,\dots,\sigma \rangle_d\), \(\langle \alpha,\sigma,\dots,\sigma, \Gamma_{i_1},\dots,\Gamma_{i_j}\rangle_d\). Notice that for degree reasons, there is no non-zero GW-invariants with only vanishing insertions. Moreover if \(\alpha = \beta\cup h\) and \(j\neq 0\), by the previous proposition and the fact that \(h\cup \Gamma_i=0\) we have \(    \langle \alpha,\sigma,\dots,\sigma,\Gamma_{i_1},\dots,\Gamma_{i_j} \rangle_d \equiv \langle \beta,\sigma,\dots,\sigma ,\Gamma_{i_1}\cup h,\dots,\Gamma_{i_j}\rangle_d +o(d,n+3) \equiv o(d,n+3)\).
Therefore the result follows if we exclude that \(\alpha\) is a \(\cup\)-power of \(h\). If \(\alpha = h^{\cup j}\) then the previous proposition implies that \(\langle \alpha,\sigma,\dots,\sigma \rangle_d \equiv \langle h,\sigma\cup h^{\cup (j-1)},\sigma,\dots,\sigma \rangle_d +o(d,n+3) \equiv o(d,n+3)\),
where in the last step we have used the divisor axiom.
\end{proof}

\subsection{Deformation of \(\QH(Y)\) along vanishing classes}

Recall that \(Y\subset X=\IGr(2,6)\) is a general hyperplane section. We start by computing explicit expressions for the degree two and three invariants \(\phi_2,\phi_3\) of the monodromy action \(\Mon(Y)\to \Aut(\HHH^*(Y,\mathbb{Z}))\) in terms of the classes \(\Gamma_1 := \Gamma_{\alpha_1}\) and \(\Gamma_2 = \Gamma_{\alpha_2}\). Denote by \(t_i\) the dual coordinate to \(\Gamma_i\) for \(i = 1,2\).


\begin{lemma}
We have \(\phi_2 = \Gamma_1^2-\Gamma_1\Gamma_2 +\Gamma_2^2\) and \(\phi_3=-2\Gamma_1^3-2\Gamma_2^3+3\Gamma_1^2\Gamma_2+3\Gamma_1\Gamma_2^2\).  
In the dual basis  \(\phi_2\) and \(\phi_3\) read (modulo scalar) as: \(s:={t_1^2} + \tone\ttwo + {t_2^2}\) and \(t := \tone \ttwo (\tone+\ttwo)\).
\end{lemma}

\begin{proof}
The two reflections associated to the simple roots \(\alpha_i\), \(i=1,2\) operate as follows (see \autoref{prop:action-of-W}): \(s_1(\Gamma_1)=-\Gamma_1\), \(s_1(\Gamma_2)=\Gamma_2-\Gamma_1\), and \(s_2(\Gamma_1)=\Gamma_1-\Gamma_2\), \(s_2(\Gamma_2)=-\Gamma_2\). Then one checks that the two reflections both fix \(\phi_2\) and \(\phi_3\).  
\end{proof}

\begin{lemma}
\label{lem_GW_prim_insertions_simplified}
Let \(\alpha_1,\cdots,\alpha_r\) be ambient classes. Then for any \(d\):
\begin{enumerate}[(1)]
    \item \(\langle \alpha_1,\cdots,\alpha_r,\Gamma_i \rangle_d=0\) for \(i=1,2\);
    \item \(\langle \alpha_1,\cdots, \alpha_r,\Gamma_1,\Gamma_1 \rangle_d=2\langle \alpha_1,\cdots,\alpha_r,\Gamma_1,\Gamma_2 \rangle_d=\langle \alpha_1,\cdots,\alpha_r,\Gamma_2,\Gamma_2 \rangle_d\);
    \item \(\langle \alpha_1,\cdots,\alpha_r,\Gamma_1,\Gamma_1,\Gamma_1 \rangle_d=\langle \alpha_1,\cdots,\alpha_r,\Gamma_2,\Gamma_2, \Gamma_2 \rangle_d=0\);
    \item \(\langle \alpha_1,\cdots,\alpha_r,\Gamma_1,\Gamma_1,\Gamma_2 \rangle_d=\langle \alpha_1,\cdots,\alpha_r,\Gamma_1,\Gamma_2, \Gamma_2 \rangle_d\).
\end{enumerate}
\end{lemma}

\begin{proof}
This follows from \autoref{prop:action-of-W} and the fact that \(\langle \alpha_1,\cdots,\alpha_r, \bullet,\cdots,\bullet \rangle_d \in \Sym(\HHHna(Y))\) is invariant under the monodromy action (because ambient classes are).
\end{proof}

Associativity of quantum product and the divisor axiom will be sufficient to compute all \(4\)-points invariants apart from \(\langle \sigma,\Gamma_1,\Gamma_1,\Gamma_2 \rangle_1\), for which we will need geometric considerations (see Subsection \ref{subsection:inv-via-geom}). All GW-invariants involving three points and ambient classes can be computed from the explicit presentation of the small quantum cohomology given in Subsection \ref{subsection:presentation} (see also \autoref{rem_3_points}). The outcome of the computations is given by the following lemma. We will use the notation \(\gamma_{100} := j^* \sigma_{100,X}\) and \(\gamma_{010} := j^* \sigma_{010,X}\). Recall also from \autoref{prop_pres_small_quantum_Y}, the notation \(\sigma = \sigma_{011} = j^*\sigma_{011,X}\).

\begin{lemma}\label{lem:non-geometric_GW}
    We have:
    \begin{itemize}
        \item \(\langle \sigma_{-110},\Gamma_1,\Gamma_1 \rangle_1 = 2 \) and \(\langle \sigma_{-011},\Gamma_1,\Gamma_1 \rangle_1= -2\),
        \item \(\langle \gamma_{100}, \sigma_{110} ,\Gamma_1,\Gamma_1 \rangle_1=\langle \gamma_{100}, \sigma_{011} ,\Gamma_1,\Gamma_1 \rangle_1=\langle \gamma_{010}, \sigma_{110} ,\Gamma_1,\Gamma_1 \rangle_1=2\) \\ and \(\langle \gamma_{010}, \sigma_{011} ,\Gamma_1,\Gamma_1 \rangle_1=0\),
        \item \(\langle h\cup h,\Gamma_1,\Gamma_1 ,\Gamma_2\rangle_1=0\),
        \item \(\langle \pt,\Gamma_1,\Gamma_1,\Gamma_2 \rangle_2 = \langle \sigma,\Gamma_1,\Gamma_1,\Gamma_2 \rangle_1\),
        \item  \(\langle \gamma_{100},\pt,\Gamma_1,\Gamma_1 \rangle_2=0\) and \(\langle \gamma_{010},\pt,\Gamma_1,\Gamma_1 \rangle_2 = 2\).
    \end{itemize}
\end{lemma}
A detailed proof can be found in \autoref{appendix:computation_of_GW_invariants}.

\subsubsection{Invariants through geometry}
\label{subsection:inv-via-geom}

In this section we compute the GW-invariant \(\langle \sigma,\Gamma_1,\Gamma_1,\Gamma_2 \rangle_1\) via geometry. 

\begin{lemma}
\label{lem_difficult_pullback_invariant}
We have \(\langle \sigma,\gamma_{100},\gamma_{100},\gamma_{010} \rangle_1=8\).
\end{lemma}

\begin{proof}
We will apply the strategy of the proof of \autoref{prop_only_some_invariants} and show that this invariant can be computed only from \(3\)-points GW-invariants. In this proof, all computations will be made in \(A:=\BQH_{\sigma}(Y)/(t_\sigma^2)\) and we will use the Schubert basis except in degree \(6\) for which we use the basis \((\gamma_{100},\gamma_{010},\Gamma_1,\Gamma_2)\). In this basis we have \(\gamma_{100} = 2\gamma_{100}^\vee + \gamma_{010}^\vee\) and \(\gamma_{010} = \gamma_{100}^\vee + 2\gamma_{010}^\vee\). 
The above invariant appears as coefficient of \(qt_\sigma \gamma_{100}^\vee\) in the product \(\gamma_{100} \star \gamma_{010}\).
In \(A\) we have \(\gamma_{010} = h\star \sigma\) (this uses \autoref{rem_3_points}), thus we get \(\gamma_{100}\star \gamma_{010}=\gamma_{100}\star (h \star \sigma)\). Using associativity, the divisor axiom and three-points GW-invariants we deduce that \(\langle \sigma,\gamma_{100},\gamma_{100},\gamma_{010}\rangle_1 = \langle \sigma,\gamma_{100},3\sigma_{-110}+\sigma_{-011},\sigma\rangle_1 + 2\). 
The GW-invariant on the RHS appears as the coefficient of \(qt_\sigma \sigma^\vee\) in the quantum product \(\gamma_{100}\star(3\sigma_{-110}+\sigma_{-011})\). Since in \(A\) we have \(\gamma_{100} = h\star h \star h - 2\sigma\star h - qt_\sigma\), again associativity, the divisor axiom and three-points GW-invariants yield \(\langle \sigma,\gamma_{100},3\sigma_{-110}+\sigma_{-011},\sigma\rangle_1=-8\langle \sigma,l,\sigma,\sigma\rangle_1 +6\). Choosing representatives of \(\sigma,l,\sigma,\sigma\) one easily shows that there are no lines passing through such representatives, i.e. \(\langle \sigma,l,\sigma,\sigma\rangle_1=0\), finally yielding the result.
\end{proof}

\begin{lemma}
\label{lem_geometric_invariant}
We have \(\langle \sigma,\sigma_{-100},\sigma_{-100},\sigma_{-100} \rangle_1 = 1\), \(\langle \sigma,\sigma_{-100},\sigma_{-100},\sigma_{-010} \rangle_1=0\), \\ \(\langle \sigma,\sigma_{-100},\sigma_{-010},\sigma_{-010} \rangle_1=0\), and \(\langle \sigma,\sigma_{-010},\sigma_{-010},\sigma_{-010} \rangle_1=0\).
\end{lemma}

\begin{proof}
The hyperplane section \(Y\subset \IGr(2,V_6)\subset \Gr(2,V_6)\) is defined as a double linear section by two generic symplectic forms \(\omega,\omega_1\) on \(V_6\). Up to the action of \(\GL_6\), one can assume that \(\omega=\sum_{i}x_i\wedge y_i\) while \(\omega_1=\sum_{i}\lambda_i x_i\wedge y_i\), for pairwise distinct coefficients \(\lambda_1,\lambda_2,\lambda_3\). The pencil \(\langle \omega,\omega_1\rangle\) defines a decomposition of \(V_6 = A\oplus B \oplus C\), where coordinates on \(A\) (resp. \(B\) and \(C\)) are \(x_1,y_1\) (resp. \(x_2,y_2\) and \(x_3,y_3\)). Choose a flag of the following form: 
\[A_1\subset A_1+B_1\subset A_1+B_1+C_1\subset A_1+B_1+C\subset A_1+B+C \subset V,\]
where \(A_1 \subset A\), \(B_1\subset B\), \(C_1\subset C\) are one-dimensional. From \cite[Section 3]{BP}, we deduce that a representative of \(\sigma_{-100}\) is given by \(X_{\sigma_{-100}}\), which can be embedded in \(Y\) as
\[Y_{A_1} := \{ [V_2]\in Y\subset \Gr(2,V)\mid A_1\subset V_2  \subset A_1 + B + C\subset V_6\} = X_{\sigma_{-100}}\]
for any \(A_1\subset A\). A representative of \(\sigma_{-010}\) is given by \(X_{\sigma_{-010}}\), which can be embedded in \(Y\) as
\begin{equation*}
    Y_{A_1,B_1}:=\{ [V_2]\in  \Gr(2,V)\mid V_2\cap (A_1+B_1)\neq 0,V_2 \subset A_1+B_1+C\subset V_6\} = X_{\sigma_{-010}}
\end{equation*}
for any \(A_1\subset A\) and \(B_1\subset B\). Here the subscript denotes the dimension of the vector space. 

Notice that we can act with the automorphism group \(\PGL(A)\times \PGL(B)\times \PGL(C)\) of \(Y\), which moves \(A_1\subset A\), \(B_1\subset B\), \(C_1\subset C\) but leaves cohomology classes invariant (since the group is connected): thus, for instance, we have \([Y_{A_1'}] = [Y_{A_1}]\) for \(A_1\neq A_1'\subset A\). The role of \(A\), \(B\), \(C\) can be permuted, giving other subvarieties of \(Y\) which do not represent the same cohomology class. In order to see this, recall the intersection form for \(Y\) in the basis given by \(\{\sigma_\alpha\}_{\alpha\in \Phi_\aleph \sqcup (-\Phi_\aleph)}\):
\[\begin{pmatrix}
    C_\aleph-2 I_\aleph&I_\aleph\\ I_\aleph &0
\end{pmatrix}\]
where \(C_\aleph\) is the Cartan matrix associated to the root system \(\Phi_\aleph\) and \(I_\aleph\) the identity matrix.
Consider the subvarieties \(Y_{B_1}\) and \(Y_{C_1}\), and \(Y_{B_1,C_1}\) and \(Y_{A_1,C_1}\). We have \(j_*[Y_{A_1}]=j_*[Y_{B_1}]=j_*[Y_{C_1}]\) and \(j_*[Y_{A_1,B_1}]=j_*[Y_{A_1,C_1}]=j_*[Y_{B_1,C_1}]\), since as subvarieties of \(X\) they define the same cohomology class (just move the complete isotropic flag with the action of the symplectic group). So, for instance, \([Y_{B_1}] - \sigma_{-100} = [Y_{B_1}] - [Y_{A_1}] = a\Gamma_1 + b\Gamma_2\) for some coefficients \(a,b\). These coefficients can be determined by the following computation:
\begin{itemize}
    \item \([Y_{A_1}] \cup [Y_{B_1}]=1\) since \(Y_{A_1}\) and \(Y_{B_1}\) intersect transversally at the point \(A_1+B_1 \in Y\).
    \item \([Y_{A_1,B_1'}]\cup [Y_{B_1}] = 0\) since \(Y_{A_1,B_1'} \cap Y_{B_1} = \emptyset\), where we have chosen \(B_1' \neq B_1\).
\end{itemize}
Using the intersection form above, we get 
\(a=1,b=0\), \emph{i.e.} \([Y_{B_1}] = \sigma_{100}\). One can do a similar reasoning with the other classes and the following computations:
\begin{itemize}
\item \([Y_{A_1,B_1}] \cup [Y_{A_1',C_1}] = 1\) since \(Y_{A_1,B_1} \cap Y_{A_1',C_1} = \{[B_1+C_1]\}\).
    \item \([Y_{A_1}] \cup [Y_{A_1',C_1}] = 0\), \([Y_{A_1}] \cup [Y_{C_1}] = 1\), \([Y_{B_1,C_1}] \cup [Y_{B_1'}] = 0\), \([Y_{A_1}] \cup [Y_{A_1',C_1}] = 0\), \([Y_{A_1,B_1}] \cup [Y_{A_1',C_1}] = 1\) by symmetry and the above computations.
\end{itemize}
All these comuptations lead to the following description of classes:
\[[Y_{A_1}]=\sigma_{-100}, \quad [Y_{B_1}]=\sigma_{100}, \quad [Y_{C_1}]=\sigma_{100}-\sigma_{010}+\sigma_{-010};\]
\[[Y_{A_1,B_1}]=\sigma_{-010},\quad [Y_{A_1,C_1}]=\sigma_{010},\quad [Y_{B_1,C_1}]=\sigma_{010}-\sigma_{100}+\sigma_{-100}.\]
Representatives for \(\sigma\) and \(\gamma_{100}\) are given by the varieties \(Y_{D_5}\) and \(Y_{E_2}\)
\begin{equation*}
    Y_{D_5}:=\{[V_2] \in Y \subset \Gr(2,V)\mid V_2\subset D_5 \} \text{ and } Z_{E_2}:=\{[V_2] \in Y \subset X \mid V_2 \cap E_2 \neq 0 \},
\end{equation*}
where \(D_5 \subset V_6\) is any \(5\)-dimensional subspace and \(E_2 \subset V_6\) is isotropic for \(\omega\) but not necessarily for \(\omega_1\).
We are now in position to compute some useful GW-invariants.
\begin{itemize}
\item \(\langle \sigma,\sigma_{-100},\sigma_{-100},\sigma_{100} \rangle_1 = \langle \sigma,\sigma_{-100},\sigma_{-100},\sigma_{100}-\sigma_{010}+\sigma_{-010} \rangle_1 = 0\): for the first GW-invariant, choose generic representatives of the four classes \(Y_{D_5}\), \(Y_{A_1}\), \(Y_{A_1'}\), \(Y_{B_1}\). We need to compute the number of lines \(l_{L_1\subset L_3}\) passing through these four subvarieties. Clearly \(L_3= A_1 + A_1' + B_1 = A + B_1\), which forces \(L_1 = B_1\); however \(L_1\) should be contained in \(D_5\), which is in general position, hence there are no lines passing through the four representatives. The second GW-invariant is computed in the same way.
\item \(\langle \sigma,\sigma_{-100},\sigma_{-100},\gamma_{100} \rangle_1 = 1\): choose generic representatives of the four classes \(Y_{D_5}\), \(Y_{A_1}\), \(Y_{A_1'}\), \(Z_{E_2}\). We need to compute the number of lines \(l_{L_1\subset L_3}\) passing through these four subvarieties. Set \(W_3 := D_5 \cap (B + C)\) and  \(W_1 := E_2 \cap (A + W_3)\), by our genericity assumption, we have \(\dim W_i = i\) for \(i \in \{1,3\}\). We have \(A \subset L_3 \subset A + W_3\), thus \(L_3 = W_1 + A\) and \(L_1 = L_3 \cap W_3\) \emph{i.e.} there is a unique line passing through the four representatives.
\item \(\langle \sigma,\sigma_{-100},\sigma_{-100},\gamma_{010} \rangle_1 = 1\): choose generic representatives of the four classes \(Y_{D_5}\), \(Y_{A_1}\), \(Y_{A_1'}\), \(Z_{E_3,E_5} := \{ [V_2] \in Y \subset \Gr(2,V) \mid V_2 \cap E_3 \neq 0,V_2 \subset E_5 \}\), where \(E_3\) is isotropic for \(\omega\) but not necessarily for \(\omega_1\) and \(E_3\subset E_5\). We need to compute the number of lines \(l_{L_1\subset L_3}\) passing through these four subvarieties. First notice that \(A_1,A_1' \subset L_3\) must be orthogonal to \(L_1\), so \(L_1\subset W_2 := (B+C) \cap D_5\cap E_5\), or equivalently \([L_1]\in \PP(W_2)\). Then \(L_3 = L_1 + A\). However, \(L_3\cap E_3\neq 0\). Working modulo \(A\), one checks that there is a unique point in \(\PP(W_2)\) satisfying this condition, thus there is a unique line passing through the four representatives.
\item \(\langle \sigma,\sigma_{-010},\sigma_{-010},\sigma_{010} \rangle_1 = \langle \sigma,\sigma_{-010},\sigma_{-010},\sigma_{010}-\sigma_{100}+\sigma_{-100} \rangle_1 = 0\): for the first GW-invariant, choose generic representatives of the four classes \(Y_{D_5}\), \(Y_{A_1,B_1}\), \(Y_{A_1',B_1'}\), \(Y_{A_1'',C_1}\). In this case $L_1$ should be contained in the triple intersection 
\[L_1\subset (A_1+B_1+C)\cap (A_1'+B_1'+C)\cap (A_1''+B+C_1)=C_1.\]
However, we also have \(L_1 \subset D_5\), which is impossible since \(D_5\cap C_1 = 0\) for \(D_5\) generic. So there is no line passing through the four representatives. The second GW-invariant is computed in the same way.
\item \(\langle \sigma,\sigma_{100}-\sigma_{010}+\sigma_{-010},\sigma_{-010},\sigma_{-010} \rangle_1 = 0\): choose generic representatives of the four classes \(Y_{D_5}\), \(Y_{C_1}\), \(Y_{A_1,B_1}\), \(Y_{A_1',B_1'}\). In this case \(L_1\) should be contained in the triple intersection 
\[L_1\subset (A_1+B_1+C)\cap (A_1'+B_1'+C)\cap D_5=C\cap D_5=C_1',\]
which in general will be different from \(C_1\). But then \(L_3\) should contain a line in \(A_1+B_1\), a line in \(A_1'+B_1'\), \(C_1\) and \(C_1'\). This is impossible, since generically these subspaces span a subspace of dimension \(4\), hence there is no line passing through the four representatives. 
\item \(\langle \sigma,\sigma_{-010},\sigma_{-010},\gamma_{100} \rangle_1 = 0\): choose generic representatives of the four classes \(Y_{D_5}\), \(Y_{A_1,B_1}\), \(Y_{A_1',B_1'}\), \(Z_{E_2}\). We need to compute the number of lines \(l_{L_1\subset L_3}\) passing through these four subvarieties. We have \(L_1\subset (A_1+B_1+C)\cap (A_1'+B_1'+C)\cap D_5=C\cap D_5\), and since the latter is \(1\)-dimensional, we get \(L_1 = C\cap D_5\). Moreover \(E_1 := L_1^{\perp_\omega} \cap E_2\) has to be contained in \(L_3\). However, \(L_3\) also contains a \(1\)-dimensional subspace of \(A_1+B_1\) and of \(A_1'+B_1'\), and this is not simultaneously possible, so there is no line passing through the four representatives.
\end{itemize}
Recall that \(\gamma_{100} = \sigma_{100}+\sigma_{-100}+\sigma_{-010}\) and \(\gamma_{010}=\sigma_{010}+\sigma_{-100}+\sigma_{-010}\). Using this and the last three bullets above, we get \(\langle \sigma,\sigma_{-100},\sigma_{-010},\sigma_{-010} \rangle_1 = \langle \sigma,\sigma_{-010},\sigma_{-010},\sigma_{-010} \rangle_1 = 0\). Using the first three bullets above, we get \(\langle \sigma,\sigma_{-100},\sigma_{-100},\sigma_{-100} \rangle_1 = 1\) and $\langle \sigma,\sigma_{-100},\sigma_{-100},\sigma_{-010} \rangle_1=0$.
\end{proof}

Using the relations \(\Gamma_1 = \gamma_{100}-2\sigma_{-100}-\sigma_{-010}\) and \(\Gamma_2 = \gamma_{010} - \sigma_{-100} - 2\sigma_{-010}\) together with the previously computed GW-invariants, one gets the following. 

\begin{proposition}
\label{prop_geometric_invariant}
 We have $\langle \sigma,\Gamma_1,\Gamma_1,\Gamma_2 \rangle_1 = 6$.
\end{proposition}

\begin{proof}
We have \(\Gamma_1=\gamma_{100}-2\sigma_{-100}-\sigma_{-010}\) and \(\Gamma_2=\gamma_{010}-\sigma_{-100}-2\sigma_{-010}\), thus we get 
\begin{equation}
\label{eq_comp_invariant_geom}
\begin{array}{rl}
\langle \sigma,\Gamma_1,\Gamma_1,\Gamma_2 \rangle_1 = &  
\langle \sigma,\gamma_{100},\Gamma_1,\Gamma_2 \rangle_1 -\langle \sigma,2\sigma_{-100}+\sigma_{-010},\Gamma_1,\Gamma_2 \rangle_1 \\
 = & \langle \sigma,\gamma_{100},\Gamma_1,\Gamma_2 \rangle_1 - \langle \sigma,2\sigma_{-100}+\sigma_{-010},\gamma_{100},\Gamma_2 \rangle_1  \\ 
 &  + \langle \sigma,2\sigma_{-100} + \sigma_{-010},2\sigma_{-100}+\sigma_{-010},\Gamma_2 \rangle_1 \\
   = & \langle \sigma,\gamma_{100},\Gamma_1,\Gamma_2 \rangle_1 - \langle \sigma,2\sigma_{-100}+\sigma_{-010},\gamma_{100},\Gamma_2 \rangle_1 \\
   & +  \langle \sigma,2\sigma_{-100}+\sigma_{-010},2\sigma_{-100}+\sigma_{-010},\gamma_{010} \rangle_1 +\\
    &   -\langle \sigma,2\sigma_{-100}+\sigma_{-010},2\sigma_{-100}+\sigma_{-010},\sigma_{-100}+2\sigma_{-010} \rangle_1.
\end{array}
\end{equation}
By \autoref{lem_GW_prim_insertions_simplified} and \autoref{lem:non-geometric_GW} we deduce that  \(\langle \sigma,\gamma_{100},\Gamma_1,\Gamma_2 \rangle_1 = 1\) and 
\begin{equation*}
    \langle \sigma,2\sigma_{-100}+\sigma_{-010},\gamma_{100},\Gamma_2\rangle_1 = -\langle \sigma,\Gamma_1,\gamma_{100},\Gamma_2\rangle_1 + \langle \sigma,\gamma_{100},\gamma_{100},\Gamma_2 \rangle_1=-1 + 0 = -1.
\end{equation*}
We determine the last two terms in (\ref{eq_comp_invariant_geom}). Start with \(a := \langle \sigma,2\sigma_{-100}+\sigma_{-010},2\sigma_{-100}+\sigma_{-010},\gamma_{010} \rangle_1\).
\begin{equation*}
\begin{array}{rl}
a = & \langle \sigma,2\sigma_{-100}+\sigma_{-010},\gamma_{100},\gamma_{010} \rangle_1 + \langle \sigma,2\sigma_{-100}+\sigma_{-010},-\Gamma_1,\gamma_{010} \rangle_1\\
    = &\langle \sigma,\gamma_{100},\gamma_{100},\gamma_{010} \rangle_1 - \langle \sigma,\Gamma_{1},\gamma_{100},\gamma_{010} \rangle_1 - \langle \sigma,\gamma_{100},\Gamma_1,\gamma_{010} \rangle_1 + \langle \sigma,\Gamma_1,\Gamma_1,\gamma_{010} \rangle_1 \\
    =& \langle \sigma,\gamma_{100},\gamma_{100},\gamma_{010} \rangle_1 - 0-0+0=\langle \sigma,\gamma_{100},\gamma_{100},\gamma_{010} \rangle_1=8,
\end{array}
\end{equation*}
where in the last line we have used \autoref{lem:non-geometric_GW} and \autoref{lem_difficult_pullback_invariant}. We deduce that
\begin{equation*}
\begin{array}{rl} 
\langle \sigma,\Gamma_1,\Gamma_1,\Gamma_2 \rangle_1 = & 1+ 1 + 8 - 4\langle \sigma,\sigma_{-100},\sigma_{-100},\sigma_{-100} \rangle_1 - 12\langle \sigma,\sigma_{-100},\sigma_{-100},\sigma_{-010} \rangle_1 \\
&    -9\langle \sigma,\sigma_{-100},\sigma_{-010},\sigma_{-010} \rangle_1 - 2\langle \sigma,\sigma_{-010},\sigma_{-010},\sigma_{-010} \rangle_1 \\
 = & 6,\\
\end{array}
\end{equation*}
where in the last equality we have used \autoref{lem_geometric_invariant}.
\end{proof}

\subsubsection{Multiplication by \(h\)}
Recall the \(\Mon(Y)\)-invariant polynomials \(s:= {t_1^2} + \tone \ttwo + {t_2^2}\) and \(t:= \tone \ttwo(\tone+\ttwo)\). Define
\begin{equation*}
    r_1:=s-\frac{3{t_1^2}}{2} \quad , \quad r_2:= s-\frac{3{t_2^2}}{2}.
\end{equation*}

The quantum potential \(F\) determining the quantum multiplication in the big quantum cohomology is also \(\Mon(Y)\)-invariant, thus it depends on \(s\) and \(t\) rather than \(\tone, \ttwo\). The following is a direct consequence of \autoref{prop_geometric_invariant} and \autoref{lem_three_point_primitive}.
\begin{proposition}
\(F\) depends non-trivially on \(s\) and \(t\).
\end{proposition}
Consider the operator \(E_Y:\BQH(Y)\to \BQH(Y)\), \(\alpha \mapsto h\star \alpha\) given by the multiplication by \(h\).
In fact, we will only consider deformations along vanishing directions of the small quantum cohomology ring, \emph{i.e.} along \(\Gamma_1,\Gamma_2\), whose dual variables are \(\tone,\ttwo\). Moreover, we are interested in deformations up to order two in \(\tone,\ttwo\), therefore we look at the ring
\begin{equation*}
    \BQH_{\Gamma_1,\Gamma_2}^{(2)}(Y):=\BQH_{\Gamma_1,\Gamma_2}(Y) / (t_1^3,t_1^2\ttwo,\tone t_2^2,t_2^3),
\end{equation*}
and the corresponding operator 
\begin{equation*}
    \hat{E}_Y: \BQH_{\Gamma_1,\Gamma_2}^{(2)}(Y) \to \BQH_{\Gamma_1,\Gamma_2}^{(2)}(Y) , \quad \alpha \mapsto h\star \alpha.
\end{equation*}
The ring \(\BQH_{\Gamma_1,\Gamma_2}^{(2)}(Y)\) is the \emph{smallest} deformation in which \(\langle \sigma,\Gamma_1,\Gamma_1,\Gamma_2 \rangle_1\) appears.

Denote by \(O(3)\) any element in the ideal \((t_1^3,t_1^2\ttwo,\tone t_2^2,t_2^3)\). We compute the action of \(\hat{E}_Y\) on the additive basis \(\mathfrak{B} := (1,h,\sigma_{110},\sigma,\gamma_{100},\gamma_{010}, \Gamma_1,\Gamma_2,\sigma_{-110},\sigma_{-011}, l, \pt)\)
of the classical cohomology \(\HHH^*(Y)\). We find the following 

\begin{theorem}
\label{thm:mult-h}
The matrix of \(\hat{E}_Y\) in the basis \(\mathfrak{B}\) is 
\begin{equation*}
\left(\begin{smallmatrix}
0 & 0 & 0 & 0 & 2q & q+2q^2s & 6q^2 \frac{\partial t}{\partial \tone} & 6q^2\frac{\partial t}{\partial \ttwo} & 0 & 0 & 0 & 0 \\
1 & 0 & 0 & 0 & 0 & 0 & 0 & 0 & q & 2q^2s & 0 & 0 \\
 0 & 1+ qs & 0 & 0 & 0 & 0 &  q\frac{\partial s}{\partial \tone} & q\frac{\partial s}{\partial \ttwo} & 0 & 0 & q & 0 \\
 0 & 1- qs & 0 & 0 & 0 & 0 &  - q\frac{\partial s}{\partial \tone} & -  q\frac{\partial s}{\partial \ttwo} & 0 & 0 & 2q^2s & 0 \\
 0 & 0  & 1+ \frac{qs}{3} &  2\frac{qs}{3} & 0 & 0 & 0 & 0 & 0 & 0 & 0 & q- \frac{2q^2s}{3} \\
 0 & 0  & 1+ \frac{qs}{3} & 1-\frac{qs}{3} & 0 & 0 & 0 & 0 & 0 & 0 & 0 & \frac{4q^2s}{3} \\
  0 & 0 & 2r_1 q & -2r_1 q & 0 & 0 & 0 & 0 & -q\tone & q\tone & 0 & -4r_1 q^2\\
 0 & 0  & 2r_2 q & -2r_2 q & 0 & 0 & 0 & 0 & -q\ttwo & q\ttwo & 0 & -4r_2 q^2\ \\
  0 & 0 & 0 & 0 & 3+ qs & 3+qs & - 3q\frac{\partial t}{\partial \tone} & - 3 q\frac{\partial t}{\partial \ttwo} & 0 & 0 & 0 & 0 \\
  0 & 0 & 0 & 0 & 1+ qs & 2 &  3q\frac{\partial t}{\partial \tone} & 3q\frac{\partial t}{\partial \ttwo} & 0 & 0 & 0 & 0 \\
  0 & 0 & 0 & 0 & 0 & 0 & 0 & 0 & 1+qs & 1-qs & 0 & 0 \\
  0 & 0 & 0 & 0 & 0 & 0 & 0 & 0 & 0 & 0 & 1 & 0 \\
\end{smallmatrix}\right).
\end{equation*}
\end{theorem}

Detailed computations can be found in \autoref{appendix:multiplication_by_h}.




\subsubsection{Multiplication by the Euler field}
 
We impose \(\ttwo=0\) and consider the multiplication by the quantum Euler vector field  which is a deformation of the classical Euler field (see \cite[(4.33)]{kontsevich1994gromov} for a definition)
\begin{equation*}
    \eu := 4h-2\tone \Gamma_1
\end{equation*} 
up to order two in \(\tone\). More precisely, we work in the ring
\begin{equation*}
\BQH_{\Gamma_1}^{(2)}(Y):=\BQH_{\Gamma_1,\Gamma_2}^{(2)}(Y) / (\ttwo).
\end{equation*}
The essential computation was made in \autoref{thm:mult-h}. We obtain the following result.

\begin{theorem}
The matrix of the multiplication by \(\eu\) in \(\BQH_{\Gamma_1}^{(2)}(Y)\) in the basis \(\mathfrak{B}\) is
\begin{equation}
\label{mat:eu-t}
    \left(\begin{smallmatrix}
 0 & 0 & 0 & 0 & 8q & 4q+ 4q^2 t_1^2 & 0 &   12q^2 t_1^2 & 0 & 0 & 0 & 0 \\
 4 & 0 & 0 & 0 & 0 & 0  & 0 & 0 & 4q & 4q^2 t_1^2 & 0 & 0 \\
  0 & 4 & 0 & 0 & 0 & 0 &  4q\tone &  2q\tone & 0 & 0 & 4q & 0 \\
  0 & 4 & 0 & 0 & 0 & 0 &  - 4q\tone &  - 2q\tone & 0 & 0 &  4q^2 t_1^2 & 0 \\
  0 & 0  & 4 & 0 & 0 & 0 & 0 & 0 & 0 & 0 & 0 & 4q- \frac{4q^2 t_1^2}{3} \\
  0 & 0  & 4 & 4 & 0 & 0 & 0 & 0 & 0 & 0 & 0 &  \frac{8q^2 t_1^2}{3} \\
   -2\tone & 0 & 0 & 0 & 0 & 0 & 0 & 0 & -2q\tone & 2q\tone & 0 & 4 q^2 t_1^2 \\
  0 & 0  & 0 & 0 & 0 & 0 & 0 & 0 & 0 & 0 & 0 &  -8 q^2 t_1^2 \\
   0 & 0 & 0 & 0 & 12 & 12 & 0 & 0 & 0 & 0 & 0 & 0 \\
   0 & 0 & 0 & 0 & 4 & 8 & 0 & 0 & 0 & 0 & 0 & 0 \\
   0 & 0 & 0 & 0 & 0 & 0 & 0 & 0 & 4 & 4 & 0 & 0 \\
   0 & 0 & 0 & 0 & 0 & 0 & 4\tone & 2\tone & 0 & 0 & 4 & 0 \\
\end{smallmatrix}\right).
\end{equation}
\end{theorem}
Detailed computations can be found in \autoref{appendix:multiplication_by_E}.

\begin{corollary}\label{cor:semi-semplicity_for_IG(2,6)}
The big quantum cohomology \(\BQH(Y)\) of a smooth hyperplane section \(Y\subset \IGr(2,6)\) is generically semisimple.
\end{corollary}
\begin{proof}
Set \(q =1\). Let \(M^{(2)}\) be the matrix \eqref{mat:eu-t} (with \(q = 1\)) and \(M\) be the matrix of multiplication by \(\eu\) in \(\BQH_\eta(Y)\) in the basis \(\mathfrak{B}\). Note that \(M\) is a deformation of \(M^{(2)}\) with coefficients in \(\QQ\formal{\mathtt{t}}\). We prove that \(M\) has simple non-zero eigenvalues in the field \(K := \CC\langle\!\langle \mathtt{t} \rangle\!\rangle\) of Puiseux series in \(\mathtt{t}\) with complex coefficients. Let \(\nu : \K^\times \to \QQ\) be the \(\mathtt{t}\)-valuation, let \((\lambda_i)_{i \in [1,12]}\) be the eigenvalues of \(M\) and let \(\chi_M = x^{12} + \sum_{i = 0}^{11} a_i x^i \in \QQ\formal{\mathtt{t}}[x]\) be the characteristic polynomial of \(M\). 

We first claim that \(\nu(\lambda_i) \geq 0\) for all \(i \in [1,12]\). Indeed, set \(\nu_m = \min_{i \in [1,12]}(\nu(\lambda_i))\) and \(I := \{ i \in [1,12] \ | \ \nu(\lambda_i) = \nu_m\}\). Then \(|I| \nu_m = \nu(a_{12 - |I|}) \geq 0\), 
proving the claim. Setting \(\mathtt{t} =0\), we get \(\chi_M\vert_{t = 0} =  x^4 (x^4 - 6912) (x^4 + 256)\) so that \(8\) of the eigenvalues are pairwise distinct with trivial valuation, while the last \(4\) have a positive valuation. Up to reordering the eigenvalues, we may assume that \(\nu(\lambda_i)  = 0\) for \(i \in [5,12]\) and \(0 < v_1 \leq v_2 \leq v_3 \leq v_4\) with \(v_i = \nu(\lambda_i)\) for \(i \in [1,4]\).

We now claim that \(\nu(a_0) = 6\), \(\nu(a_1) \geq 5\), \(\nu(a_2) \geq 4\)  and \(\nu(a_3) \geq 3\). We start with \(i = 0\), we have \(a_0 = \det(M)\). Let \(F_k\) be the Fitting ideal of order \(k \in [0,3]\) of \(M^{(2)}\) in \(\QQ\formal{\mathtt{t}}\). Then \(F_0= (\mathtt{t}^6)\), \(F_1 = (\mathtt{t}^4)\), \(F_2 = (\mathtt{t}^2)\), \(F_3 = (\mathtt{t})\) and \(F_4 = (1)\). Write \(M = M^{(2)} + \mathtt{t}^3 M'\) with \(M'\) having coefficients in  \(\QQ\formal{\mathtt{t}}\). Expanding the determinant of \(M\), we see that \(a_0 = \det(M^{(2)}) + \sum_{i = 1}^{12}\mathtt{t}^{3i} b_{i}\) where \(b_i \in F_i\). This proves the claim (note that \(\det(M^{(2)})\) spans \(F_0\), actually \(\det(M^{(2)}) = 150994944 \mathtt{t}^6 +  301989888 \mathtt{t}^8 + 150994944 \mathtt{t}^{10}\)).
The cases for \(i > 0\) are similar: the coefficient \(a_i\) is the sum of principal minors of order \(12 - i\) of \(M\). Using this, the fact that the sums of principal minors of size \(11\), \(10\) and \(9\) of \(M^{(2)}\) vanish, and similar arguments, the claim follows.

Our next claim is that \(\nu(\lambda_i) = \frac{3}{2}\) for \(i \in [1,4]\). Assume that \(v_3 < v_4\). By the previous claims, we have \(v_1 + v_2 + v_3 + v_4 = \nu(a_0) = 6\) and \(v_1 + v_2 + v_3 = \nu(a_1) \geq 5\), thus \(v_4 \leq 1\) which is impossible since \(6 \leq 4 v_4\). The same arguments imply that \(v_1 = v_2 = v_3 = v_4\) and the claim follows. 

We now finish the proof. For \(i \in [1,4]\), write \(\lambda_i = \mathtt{t}^{\frac{3}{2}}\mu_i\) with \(\mu_i \in K\) and \(\nu(\mu_i) = 0\) and let \(\mu_i^{(0)} = \mu_i\vert_{\mathtt{t} = 0}\). By the previous claims, we have \(\mu_1^{(0)}\mu_2^{(0)}\mu_3^{(0)}\mu_4^{(0)} = 150994944\) while the other symmetric functions in \((\mu_i^{(0)})_{i \in [1,4]}\) vanish since \(\nu(a_i) > (4 - i) \frac{3}{2}\) for \(i \in [1,3]\). In particular \(\{\mu_1^{(0)},\mu_2^{(0)},\mu_3^{(0)},\mu_4^{(0)}\}\) is the set of fourth roots of \(150994944\) which are all distinct so that all eigenvalues of \(M\) are distinct.

This implies that \(\BQH_\eta(Y) \otimes_\QQ K \simeq K^{12}\) as an algebra and \(\BQH_\eta(Y)\) is semisimple. 
\end{proof}

\begin{corollary}The \emph{coarse atoms} of a smooth hyperplane section of $\IGr(2,6)$ are one dimensional.
\end{corollary}

\begin{proof}
  Coarse atoms are generalized eigenspaces for the action of the Euler field on the big quantum cohomology for a certain choice of the quantum parameters (we refer to \cite{kkpy}). By the proof of the previous theorem, we have shown that the multiplication by the Euler vector field has only simple eigenvalues, hence the claim.
\end{proof}

\begin{corollary}\label{cor:Dubrovin_for_IG(2,6)}Dubrovin's conjecture holds for smooth hyperplane sections \(Y\) of $\IGr(2,6)$.
\end{corollary}
\begin{proof}
Use \autoref{cor:semi-semplicity_for_IG(2,6)} and the fact that \(D^b(Y)\) has a full exceptional collection by \cite{kuznetsov2008exceptional}.
\end{proof}

\appendix

\section{Computations of GW-invariants}\label{appendix:computation_of_GW_invariants}

We prove \autoref{lem:non-geometric_GW}. We thus work with \(Y \subset \IGr(2,6)\) a smooth general hyperplane section.

\begin{lemma}
\label{lem_three_point_primitive}
We have \(\langle \sigma_{-110},\Gamma_1,\Gamma_1 \rangle_1 = 2\) and \(\langle \sigma_{-011},\Gamma_1,\Gamma_1 \rangle_1 =-2\).
\end{lemma}

\begin{proof}
We work in the small quantum cohomology \(\QH(Y)\). The only element of degree six in the kernel of the multiplication by \(h\) is \(\tau = \pt + q(\sigma-\sigma_{110})\). Since \(\Gamma_1\star_0 h = 0\), we deduce that \(\Gamma_1 \star_0 \Gamma_1\) is in the kernel of the multiplication by \(h\) as well, hence it is a multiple of \(\tau\). Since \(\Gamma_1 \cup \Gamma_1  = -2  \pt\) and \(\sigma_{110}^\vee = \sigma_{-110}\), \(\sigma_{011}^\vee = \sigma_{-011}\), we deduce that \(-2 \tau = \Gamma_1\star_0 \Gamma_1 = -2 \pt + q(\langle \sigma_{-011},\Gamma_1,\Gamma_1 \rangle_1 \sigma + \langle \sigma_{-110},\Gamma_1,\Gamma_1 \rangle_1 \sigma_{110})\), hence the result.
\end{proof}

Recall the definition of \(\gamma_{100} := j^* \sigma_{100,X}\) and \(\gamma_{010} := j^* \sigma_{010,X}\).

\begin{lemma}
\label{lem_4_GW_4points}
We have \(\langle \gamma_{100}, \sigma_{110} ,\Gamma_1,\Gamma_1 \rangle_1 = \langle \gamma_{100}, \sigma_{011} ,\Gamma_1,\Gamma_1 \rangle_1 = \langle \gamma_{010}, \sigma_{110} ,\Gamma_1,\Gamma_1 \rangle_1 = 2\)  and \(\langle \gamma_{010}, \sigma_{011} ,\Gamma_1,\Gamma_1 \rangle_1 = 0\).
\end{lemma}

\begin{proof}
These GW-invariants appear as the coefficients of \(q\tone \Gamma_1^\vee\) in the quantum products \(\gamma_{100}\star \sigma_{110}\), \(\gamma_{100}\star \sigma_{011}\), \(\gamma_{010}\star \sigma_{110}\) and \(\gamma_{010}\star \sigma_{011}\). We compute these products. Since we are interested in the coefficient of \(q\tone \Gamma_1^\vee\), we will work in the deformation of the small quantum cohomology along \(\tone\). Since \(\deg(\tone)=-2\), for degree reasons and using \autoref{lem_GW_prim_insertions_simplified} (1), in any product of ambient classes of degree at most \(4\), there is no degree \(1\) term in \(\tone\). We get:
\begin{equation*}
\begin{array}{ll} 
\gamma_{100}=(\sigma_{110}-\sigma_{011})\star h +O({t_1^2}), &    \sigma_{110}\star \sigma_{110}=  4\sigma_{-110}+2\sigma_{-011}+2q+O({t_1^2}), \\
\gamma_{010}=\sigma_{011}\star h +O({t_1^2}), & 
    \sigma_{110}\star \sigma_{011}=  2\sigma_{-110}+\sigma_{-011}+q+O({t_1^2}), \\ 
 &   \sigma_{011}\star \sigma_{011}= \sigma_{-110}+\sigma_{-011}+O({t_1^2}). \\
\end{array}
\end{equation*}
By the associativity of the quantum product we deduce that
\begin{equation*}
\gamma_{100}\star \sigma_{110}= ((\sigma_{110}-\sigma_{011})\star h)\star \sigma_{110}+O({t_1^2}) =
    h\star (2\sigma_{-110}+\sigma_{-011}+q)+O({t_1^2}).
\end{equation*}
By the divisor axiom and \autoref{lem_three_point_primitive}, the coefficient of $q\tone \Gamma_1^\vee$ of the last expression is given by \(\langle h, 2\sigma_{-110}+\sigma_{-011},\Gamma_1,\Gamma_1 \rangle_1=\langle 2\sigma_{-110}+\sigma_{-011},\Gamma_1,\Gamma_1 \rangle_1 = 2\) thus
\(\langle \gamma_{100}, \sigma_{110} ,\Gamma_1,\Gamma_1 \rangle_1 = 2\). The other invariants can be computed similarly using the following formulas
\[\begin{array}{l}
    \gamma_{100}\star \sigma_{011}= h\star (\sigma_{-110}+q)+O({t_1^2}), \\
    \gamma_{010}\star \sigma_{110}= h\star (2\sigma_{-110}+\sigma_{-011}+q)+O({t_1^2}), \\
    \gamma_{010}\star \sigma_{011}= h\star (\sigma_{-110}+\sigma_{-011})+O({t_1^2}),
\end{array}\]
the divisor axiom and \autoref{lem_three_point_primitive}.
\end{proof}

\begin{lemma}
\label{lem_hcuph_invar_prim}We have
$\langle h\cup h,\Gamma_1,\Gamma_1 ,\Gamma_2\rangle_1=0$.
\end{lemma}

\begin{proof}
The above GW-invariant appears as the coefficient of \(q\tone \Gamma_2^\vee\) in the quantum product \((h\cup h)\star \Gamma_1\). As in the previous proof, we will work in the deformation of the small quantum cohomology along \(\tone\). For degree reasons we have \(h\star h=h \cup h + O({t_1^2})\). Moreover we have
\begin{equation*}
\begin{array}{ll}
h\star \Gamma_1 & = q\tone (\langle h,\Gamma_1,\Gamma_1,\sigma_{-110}\rangle_1 \sigma_{110} + \langle h,\Gamma_1,\Gamma_1,\sigma_{-011}\rangle_1 \sigma_{011})+O({t_1^2})  \\
& =q\tone(2\sigma_{110}-2\sigma_{011})+O({t_1^2}), 
\end{array}
\end{equation*}
where we used the divisor axiom and \autoref{lem_three_point_primitive}. By associativity of \(\star\), we get
\begin{equation*}
    (h\cup h)\star \Gamma_1=h\star (h\star \Gamma_1)+O({t_1^2})=h\star(q\tone(2\sigma_{110}-2\sigma_{011})) +O({t_1^2}).
\end{equation*}
Since the cup product of ambient classes is ambient, we deduce that the coefficient of \(q\tone \Gamma_2^\vee\) in the last expression is equal to zero, hence the statement of the lemma follows.
\end{proof}

\begin{lemma}
\label{lem_interesting_invariants_are_opposite}We have an equality
\(\langle \pt,\Gamma_1,\Gamma_1,\Gamma_2 \rangle_2=\langle \sigma,\Gamma_1,\Gamma_1,\Gamma_2 \rangle_1\).
\end{lemma}

\begin{proof}
The former GW-invariant appears as the coefficient of \(q^2\tone \Gamma_2^\vee\) in the quantum product \(\pt\star \Gamma_1\). We will work in the deformation of the small quantum cohomology along \(\tone\). By the first point of \autoref{lem_GW_prim_insertions_simplified}, if \(\alpha,\beta,\gamma\) are ambient classes, the coefficient of \(q^d \tone \gamma^\vee\) in the product \(\alpha\star \beta\) vanishes for any \(d\). We thus have \(l\star h = \pt +q\sigma_{110} + O({t_1^2})\) and 
\[l\star \Gamma_1 = 
    q\tone \langle l,\Gamma_1,\Gamma_1,[Y]\rangle_1 \pt + q^2\tone \langle l,\Gamma_1,\Gamma_1,\sigma_{-110}\rangle_2 \sigma_{110} + q^2\tone \langle l,\Gamma_1,\Gamma_1,\sigma_{-011}\rangle_2 \sigma_{011} +O({t_1^2}).\]
We have \(\langle l,\Gamma_1,\Gamma_1,[Y]\rangle_1 = 0\) since \([Y]\) appears as an insertion. By associativity, we deduce
\[\begin{array}{ll}    \pt\star \Gamma_1& = h\star (l\star \Gamma_1) - q\sigma_{110}\star \Gamma_1  +O({t_1^2}) \\
  &  =q^2\tone \langle l,\Gamma_1,\Gamma_1,\sigma_{-110}\rangle_2 \sigma_{110}\star h + q^2\tone \langle l,\Gamma_1,\Gamma_1,\sigma_{-011}\rangle_2 \sigma_{011}\star h - q\sigma_{110}\star \Gamma_1 +O({t_1^2}).
\end{array}\]
Since the cup product of ambient classes is ambient, we deduce that the first two terms in the last expression do not contribute to the coefficient of \(q^2\tone \Gamma_2^\vee\). Therefore this coefficient, which is nothing else than \(\langle \pt,\Gamma_1,\Gamma_1,\Gamma_2 \rangle_2\), is equal to
\begin{equation*}
    -\langle \sigma_{110},\Gamma_1,\Gamma_1 ,\Gamma_2\rangle_1 = \langle -h\cup h + \sigma,\Gamma_1,\Gamma_1 ,\Gamma_2\rangle_1,
\end{equation*}
which by \autoref{lem_hcuph_invar_prim} is equal to \(\langle \sigma,\Gamma_1,\Gamma_1 ,\Gamma_2\rangle_1\).
\end{proof}

\begin{lemma}
\label{lem_4_pt_quadric}
We have \(\langle \gamma_{100},\pt,\Gamma_1,\Gamma_1 \rangle_2 = 0\) and \(\langle \gamma_{010},\pt,\Gamma_1,\Gamma_1 \rangle_2 = 2\).
\end{lemma}

\begin{proof}
The above GW-invariants appear as the coefficients of \(q^2\tone \Gamma_1^\vee\) in the quantum products \(\gamma_{100} \star \pt\), \(\gamma_{010}\star \pt\). We will work in the deformation of the small quantum cohomology along \(\tone\). As in \autoref{lem_4_GW_4points} we have \(\gamma_{100} = (\sigma_{110}-\sigma_{011})\star h +O({t_1^2})\) and \(\gamma_{010} = \sigma_{011}\star h +O({t_1^2})\). Moreover,  for degree reasons and using \autoref{lem_GW_prim_insertions_simplified}.(1), we have
\[\sigma_{110}\star \pt = 2q\sigma_{-110}+q\sigma_{-011}+q^2+O({t_1^2}) \quad \text{ and } \quad \sigma_{011}\star \pt = q\sigma_{-110}+q^2+O({t_1^2}).\]
Putting these together and by the associativity of the quantum product we deduce that
\begin{equation*}
    \gamma_{100}\star \pt= h\star ((\sigma_{110}-\sigma_{011})\star \pt)+O({t_1^2})=qh\star(\sigma_{-110}+\sigma_{-011})+O({t_1^2}).
\end{equation*}
Now, by the divisor axiom and \autoref{lem_three_point_primitive}, the coefficient of \(q^2\tone \Gamma_1^\vee\) of the last expression vanishes, and therefore \(\langle \gamma_{100}, \pt ,\Gamma_1,\Gamma_1 \rangle_2 = 0\). Similarly, we have
\begin{equation*}
    \gamma_{010}\star \pt = h\star (\sigma_{011}\star \pt)+O({t_1^2})=qh\star(\sigma_{-110}+q)+O({t_1^2}),
\end{equation*}
from which we deduce that \(\langle \gamma_{010},\pt,\Gamma_1,\Gamma_1 \rangle_2=2\).
\end{proof}

We finish this appendix with a lemma that we shall use in \autoref{appendix:multiplication_by_h}.

\begin{lemma}
\label{lem_4pts_avec_l}
We have \(\langle \sigma_{-110},l,\Gamma_1,\Gamma_1 \rangle_2 = 0\) and \(\langle \sigma_{-011},l,\Gamma_1,\Gamma_1 \rangle_2 = 2\).
\end{lemma}

\begin{proof}
Note that, by the intersection form on vanishing classes, we have \(2\Gamma_1^\vee + \Gamma_2^\vee = -
\Gamma_1\) and \(\Gamma_1^\vee + 2\Gamma_2^\vee = - \Gamma_2\). The GW-invariants in the statement of the lemma appear as the coefficients of \(q^2\tone \Gamma_1^\vee\) in the quantum products \(\sigma_{-110}\star l\) and \(\sigma_{-011}\star l\). We work in the deformation of the small quantum cohomology along \(\tone\). By the divisor axiom, \autoref{lem_three_point_primitive} and \autoref{lem_GW_prim_insertions_simplified}, we have:
\[\begin{array}{ll}
    \sigma_{-011}\star h & = l + q\tone \langle h, \sigma_{-011}, \Gamma_1,\Gamma_1 \rangle_1 \Gamma_1^\vee+q\tone \langle h, \sigma_{-011}, \Gamma_1,\Gamma_2 \rangle_1 \Gamma_2^\vee+O({t_1^2}) \\
& =l+q\tone(-2 \Gamma_1^\vee- \Gamma_2^\vee) + O(\tone)= l + q\tone\Gamma_1 + O({t_1^2}).
\end{array}\]
Moreover, we have \(\sigma_{-011}\star \sigma_{-110}= q\sigma_{-110}+O({t_1^2})\) and \(\sigma_{-011}\star \sigma_{-011}= q^2 + O({t_1^2})\). By the associativity of the quantum product we deduce that we have
\begin{equation*}
\sigma_{-110}\star l = (\sigma_{-110}\star \sigma_{-011})\star h - q\tone \sigma_{-110}\star \Gamma_1+O({t_1^2}) = qh\star \sigma_{-110} - q\tone \sigma_{-110}\star \Gamma_1 +O({t_1^2}).
\end{equation*}
By the divisor axiom, the coefficient of \(q^2\tone \Gamma_1^\vee\) of the last expression is equal to
\begin{equation*}
    \langle h,\sigma_{-110},\Gamma_1,\Gamma_1\rangle_1-\langle \sigma_{-110},\Gamma_1,\Gamma_1 \rangle_1=0;
\end{equation*}
therefore \(\langle \sigma_{-110},l,\Gamma_1,\Gamma_1 \rangle_2 = 0\). Similarly, we have 
\[\sigma_{-011}\star l = 
 (\sigma_{-011}\star \sigma_{-011})\star h -q\tone \sigma_{-011}\star \Gamma_1+O({t_1^2}) = q^2h - q\tone \sigma_{-011}\star \Gamma_1 +O({t_1^2}),\]
from which we deduce that \(\langle \sigma_{-011},l,\Gamma_1,\Gamma_1 \rangle_2=2\).
\end{proof}

\section{Multiplication by \(h\)}\label{appendix:multiplication_by_h}

In this appendix we give the details of the proof of \autoref{thm:mult-h}.

\begin{lemma}We have 
\(\hat{E}_Y(1) = h\) and \(\hat{E}_Y(h) = (1+ qs)\sigma_{110}+(1-qs)\sigma_{011}\).
\end{lemma}

\begin{proof}
The first is trivial, and the second is a direct computation:
\begin{equation*}
    h\star h = \sigma_{110}+\sigma_{011} +\frac{qs}{2}(\langle  h,h,\Gamma_1,\Gamma_1,\sigma_{-110} \rangle_1 \sigma_{110}+\langle  h,h,\Gamma_1,\Gamma_1,\sigma_{-011} \rangle_1 \sigma_{011})+O(3).
\end{equation*}
 The above computation is a consequence of the definition of the quantum product and of \autoref{lem_GW_prim_insertions_simplified}. Then the divisor axiom and \autoref{lem_three_point_primitive} conclude the proof.
\end{proof}

\begin{lemma}
We have \(\hat{E}_Y(\sigma_{110}) = (1 + \frac{qs}{3})\gamma_{100}+(1 + \frac{qs}{3})\gamma_{010} + 2qr_1 \Gamma_1 + 2qr_2 \Gamma_2\) and \\ \(\hat{E}_Y(\sigma_{011})=\frac{2qs}{3}\gamma_{100}+(1-\frac{qs}{3})\gamma_{010} - 2q r_1 \Gamma_1- 2q r_2 \Gamma_2\).

\end{lemma}

\begin{proof}
Recall Poincaré duality on \(\HHH^{\dim Y}(Y)\) obtained from the intersection matrix of \(Y\):
\begin{equation*}
    3\gamma_{100}^\vee = 2\gamma_{100}-\gamma_{010}, \quad 3\gamma_{010}^\vee =-\gamma_{100} + 2\gamma_{010}
, \quad -3\Gamma_{1}^\vee = 2\Gamma_{1}-\Gamma_{2}, \quad -3\Gamma_{2}^\vee = -\Gamma_{1} + 2\Gamma_{2}.
\end{equation*}
The definition of quantum product and \autoref{lem_GW_prim_insertions_simplified} yield
\begin{equation*}
\begingroup
 \begin{array}{ll} 
h\star \sigma_{110} = & \displaystyle{\gamma_{100}+\gamma_{010} + \frac{qs}{2}\langle h,\sigma_{110},\Gamma_1,\Gamma_1,\gamma_{100} \rangle_1 \gamma_{100}^\vee + \frac{qs}{2}\langle h,\sigma_{110},\Gamma_1,\Gamma_1,\gamma_{010} \rangle_1 \gamma_{010}^\vee} \\[1.4ex]
    & \displaystyle{+  \frac{q\ttwo(2\tone+\ttwo)}{2}\langle h,\sigma_{110},\Gamma_1,\Gamma_2,\Gamma_1 \rangle_1 \Gamma_1^\vee + \frac{q\tone(\tone+2\ttwo)}{2}\langle h,\sigma_{110},\Gamma_1,\Gamma_1,\Gamma_2 \rangle_1 \Gamma_2^\vee +O(3).}
    \end{array}
    \endgroup
\end{equation*}
By the divisor axiom, \autoref{lem_4_GW_4points}, \autoref{lem_hcuph_invar_prim} and \autoref{prop_geometric_invariant}, we get  

\begin{equation*}
 h\star \sigma_{110} = \gamma_{100} + \gamma_{010} + qs(\gamma_{100}^\vee + \gamma_{010}^\vee) - 3q \ttwo(2\tone+\ttwo) \Gamma_1^\vee - 3q \tone(\tone+2\ttwo)\Gamma_2^\vee +O(3).
\end{equation*}
The first statement follows from this and the previous formulae for Poincaré duals. For the second statement, the same arguments yield
\begin{equation*}
\begingroup
 \begin{array}{rl} 
    h\star \sigma_{011} = & \displaystyle{\gamma_{010} + \frac{qs}{2}\langle h,\sigma_{011},\Gamma_1,\Gamma_1,\gamma_{100} \rangle_1 \gamma_{100}^\vee + \frac{qs}{2}\langle h,\sigma_{011},\Gamma_1,\Gamma_1,\gamma_{010} \rangle_1 \gamma_{010}^\vee} \\[1.4ex]
    &  \displaystyle{+  \frac{q\ttwo(2\tone+\ttwo)}{2}\langle h,\sigma_{011},\Gamma_1,\Gamma_2,\Gamma_1 \rangle_1 \Gamma_1^\vee + \frac{q\tone(\tone+2\ttwo)}{2}\langle h,\sigma_{011},\Gamma_1,\Gamma_1,\Gamma_2 \rangle_1 \Gamma_2^\vee +O(3)} \\[1.4ex]
  =   &  \displaystyle{ \gamma_{010} + qs\gamma_{100}^\vee + 3q \ttwo(2\tone+\ttwo) \Gamma_1^\vee + 3q \tone(\tone+2\ttwo) \Gamma_2^\vee + O(3),}
     \end{array}
    \endgroup
\end{equation*}
which implies the desired formula.
\end{proof}

\begin{lemma} We have 
$\hat{E}_Y(\gamma_{100})=(3+ qs)\sigma_{-110}+(1+ qs)\sigma_{-011}+2q$ and\\ $\hat{E}_Y(\gamma_{010})=(3+qs)\sigma_{-110}+2\sigma_{-011}+q+ 2q^2s$.
\end{lemma}

\begin{proof}
The definition of quantum product and \autoref{lem_GW_prim_insertions_simplified} yield
\begin{equation*}
\begingroup
 \begin{array}{rl} 
    h\star \gamma_{100} = & \displaystyle{3\sigma_{-110}+\sigma_{-011}+2q+\frac{qs}{2}\langle h,\gamma_{100},\Gamma_1,\Gamma_1,\sigma_{110} \rangle_1 \sigma_{-110}} \\[1.4ex]
 &  \displaystyle{+\frac{qs}{2}\langle h,\gamma_{100},\Gamma_1,\Gamma_1,\sigma_{011} \rangle_1 \sigma_{-011} + \frac{q^2s}{2}\langle h,\gamma_{100},\Gamma_1,\Gamma_1,\pt \rangle_2 +O(3).}
    \end{array}
    \endgroup
\end{equation*}
By the divisor axiom, \autoref{lem_4_GW_4points} and \autoref{lem_4_pt_quadric}, we get 
\begin{equation*}
    h\star \gamma_{100} = 3\sigma_{-110}+\sigma_{-011}+2q + qs\sigma_{-110}+ qs\sigma_{-011} +O(3),
\end{equation*}
which implies the first statement. The same arguments yield
\begin{equation*}
\begingroup
 \begin{array}{rl} 
    h\star \gamma_{010} = & \displaystyle{3\sigma_{-110}+2\sigma_{-011}+q+\frac{qs}{2}\langle h,\gamma_{010},\Gamma_1,\Gamma_1,\sigma_{110} \rangle_1 \sigma_{-110} }\\[1.4ex]
    & \displaystyle{+\frac{qs}{2}\langle h,\gamma_{010},\Gamma_1,\Gamma_1,\sigma_{011} \rangle_1 \sigma_{-011} + \frac{q^2s}{2}\langle h,\gamma_{010},\Gamma_1,\Gamma_1,\pt \rangle_2 +O(3).} \\[1.4ex]
    = & 3\sigma_{-110}+2\sigma_{-011}+q+ qs\sigma_{-110}+ 2q^2s+O(3),
    \end{array}
    \endgroup
\end{equation*}
which implies the second statement.
\end{proof}

\begin{lemma} We have
\(\hat{E}_Y(\Gamma_1) = q\frac{\partial s}{\partial \tone}(\sigma_{110} - \sigma_{011})+\frac{\partial t}{\partial \tone}(6q^2+ 3q(-\sigma_{-110}+\sigma_{-011}))\) and\\ \(\hat{E}_Y(\Gamma_2)=q\frac{\partial s}{\partial \ttwo}(\sigma_{110} -\sigma_{011})+\frac{\partial t}{\partial \ttwo}(6q^2+ 3q(-\sigma_{-110}+\sigma_{-011}))\).
\end{lemma}

\begin{proof}
The definition of quantum product and \autoref{lem_GW_prim_insertions_simplified} yield
\begin{equation*}
\begingroup
 \begin{array}{rl} 
     h\star \Gamma_{1} = & \displaystyle{\frac{q(2\tone+\ttwo)}{2}\langle h,\Gamma_1,\Gamma_1,\sigma_{-110}\rangle_1\sigma_{110} + \frac{q(2\tone+\ttwo)}{2}\langle h,\Gamma_1,\Gamma_1,\sigma_{-011}\rangle_1\sigma_{011}} \\[1.4ex]
    & \displaystyle{ +\frac{q\ttwo(2\tone+\ttwo)}{2}\langle h,\Gamma_{1},\Gamma_1,\Gamma_2,\sigma_{110} \rangle_1 \sigma_{-110} + \frac{q\ttwo(2\tone+\ttwo)}{2}\langle h,\Gamma_{1},\Gamma_1,\Gamma_2,\sigma_{011} \rangle_1 \sigma_{-011}} \\[1.4ex]
& \displaystyle{+\frac{q^2\ttwo(2\tone+\ttwo)}{2}\langle h,\Gamma_{1},\Gamma_1,\Gamma_2,\pt \rangle_2 +O(3).}
   \end{array}
    \endgroup
\end{equation*}
By the divisor axiom, \autoref{lem_three_point_primitive}, \autoref{lem_hcuph_invar_prim}, \autoref{lem_interesting_invariants_are_opposite} and \autoref{prop_geometric_invariant}, we get
\begin{equation*}
h\star \Gamma_{1} = q(2\tone+\ttwo)(\sigma_{110} - \sigma_{011})+({t_2^2}+2\tone\ttwo)( 6 q^2+3q(-\sigma_{-110}+\sigma_{-011})) +O(3),
\end{equation*}
which implies the first statement. The same arguments yield
\begin{equation*}
\begingroup
 \begin{array}{rl} 
     h\star \Gamma_{2} = & \displaystyle{\frac{q(\tone+2\ttwo)}{2}\langle h,\Gamma_1,\Gamma_1,\sigma_{-110}\rangle_1\sigma_{110} + \frac{q(\tone+2\ttwo)}{2}\langle h,\Gamma_1,\Gamma_1,\sigma_{-011}\rangle_1\sigma_{011}} \\[1.4ex]
    &  \displaystyle{+\frac{q\tone(\tone+2\ttwo)}{2}\langle h,\Gamma_{1},\Gamma_1,\Gamma_2,\sigma_{110} \rangle_1 \sigma_{-110} + \frac{q\tone(\tone+2\ttwo)}{2}\langle h,\Gamma_{1},\Gamma_1,\Gamma_2,\sigma_{011} \rangle_1 \sigma_{-011} } \\[1.4ex]
 &  \displaystyle{+\frac{q^2\tone(\tone+2\ttwo)}{2}\langle h,\Gamma_{1},\Gamma_1,\Gamma_2,\pt \rangle_2 +O(3)} \\[1.4ex]
  = & q(\tone+2\ttwo)(\sigma_{110}-  \sigma_{011})+({t_1^2}+2\tone\ttwo)(6q^2+ 3q(-\sigma_{-110}+\sigma_{-011})) +O(3),
   \end{array}
    \endgroup
\end{equation*}
which implies the second statement.
\end{proof}

\begin{lemma} We have
$\hat{E}_Y(\sigma_{-110})=(1+ qs)l+qh-q\tone\Gamma_1-q\ttwo\Gamma_2$ and\\ $\hat{E}_Y(\sigma_{-011})=(1- qs)l+ 2q^2sh+q\tone\Gamma_1+q\ttwo\Gamma_2$.
\end{lemma}

\begin{proof}
The definition of quantum product and \autoref{lem_GW_prim_insertions_simplified} yield
\begin{equation*}
\begingroup
 \begin{array}{rl} 
    h\star \sigma_{-110}= & \displaystyle{l+qh+\frac{qs}{2}\langle h,\sigma_{-110},\Gamma_1,\Gamma_1,h \rangle_1 l+\frac{q^2 s}{2}\langle h,\sigma_{-110},\Gamma_1,\Gamma_1,l \rangle_2 h} \\[1.4ex]
    & \displaystyle{+ \frac{q(2\tone+\ttwo)}{2}\langle h,\sigma_{-110},\Gamma_1,\Gamma_1 \rangle_1\Gamma_1^\vee + \frac{q(\tone+2\ttwo)}{2}\langle h,\sigma_{-110},\Gamma_1,\Gamma_1 \rangle_1\Gamma_2^\vee +O(3).}
   \end{array}
    \endgroup
\end{equation*}
By the divisor axiom, \autoref{lem_three_point_primitive}, \autoref{lem_4pts_avec_l}, we get 
\begin{equation*}
    h\star \sigma_{-110}=  l+qh+ qs l+ q(2\tone+\ttwo)\Gamma_1^\vee +q(\tone+2\ttwo)\Gamma_2^\vee +O(3),
\end{equation*}
which implies the first statement. The same arguments yield
\begin{equation*}
\begingroup
 \begin{array}{rl} 
    h\star \sigma_{-011}= & \displaystyle{l+\frac{qs}{2}\langle h,\sigma_{-011},\Gamma_1,\Gamma_1,h \rangle_1 l+\frac{q^2 s}{2}\langle h,\sigma_{-011},\Gamma_1,\Gamma_1,l \rangle_2 h}\\
    & \displaystyle{+ \frac{q(2\tone+\ttwo)}{2}\langle h,\sigma_{-011},\Gamma_1,\Gamma_1 \rangle_1\Gamma_1^\vee + \frac{q(\tone+2\ttwo)}{2}\langle h,\sigma_{-011},\Gamma_1,\Gamma_1 \rangle_1\Gamma_2^\vee +O(3).} \\[1.4ex]
   = & l- qs l+ 2q^2sh-q(2\tone+\ttwo)\Gamma_1^\vee - q(\tone+2\ttwo)\Gamma_2^\vee +O(3),
   \end{array}
    \endgroup
\end{equation*}
which implies the second statement.
\end{proof}

\begin{lemma}We have 
$\hat{E}_Y(l)=l\star h = \pt+q\sigma_{110}+ 2q^2s\sigma_{011}$ and\\ $\hat{E}_Y(\pt)=(q-\frac{2q^2s}{3})\gamma_{100}+\frac{4q^2s}{3}\gamma_{010}-4r_1q^2\Gamma_1-4r_1q^2 \Gamma_2$.
\end{lemma}

\begin{proof}
The definition of quantum product and \autoref{lem_GW_prim_insertions_simplified} yield
\begin{equation*}
    h\star l = \pt +q\sigma_{110} +\frac{q^2 s}{2}\langle h,l,\Gamma_1,\Gamma_1,\sigma_{-110}\rangle_2 \sigma_{110}+\frac{q^2 s}{2}\langle h,l,\Gamma_1,\Gamma_1,\sigma_{-011}\rangle_2 \sigma_{011} +O(3).
\end{equation*}
By the divisor axiom and \autoref{lem_4pts_avec_l} one obtains
\begin{equation*}
    l\star h = \pt+q\sigma_{110}+ 2q^2s\sigma_{011}+O(3).
\end{equation*}
The definition of quantum product and \autoref{lem_GW_prim_insertions_simplified} yield
\begin{equation*}
\begingroup
 \begin{array}{rl} 
     h\star \pt = & \displaystyle{q\gamma_{100} + \frac{q^2s}{2} \langle h,\pt,\Gamma_1,\Gamma_1,\gamma_{100}\rangle_2\gamma_{100}^\vee +\frac{q^2s}{2} \langle h,\pt,\Gamma_1,\Gamma_1,\gamma_{010}\rangle_2\gamma_{010}^\vee } \\[1.4ex]
 &  \displaystyle{  + \frac{q^2\ttwo(2\tone+\ttwo)}{2}\langle h , \pt, \Gamma_1,\Gamma_1,\Gamma_2 \rangle_2 \Gamma_1^\vee  + \frac{q^2\tone(\tone+2\ttwo)}{2}\langle h , \pt, \Gamma_1,\Gamma_1,\Gamma_2 \rangle_2 \Gamma_2^\vee +O(3).}
  \end{array}
    \endgroup
\end{equation*}
By the divisor axiom, \autoref{lem_4_pt_quadric}, \autoref{lem_interesting_invariants_are_opposite} and \autoref{prop_geometric_invariant}, we get
\begin{equation*}
    h\star \pt = q\gamma_{100}  +2q^2s \gamma_{010}^\vee + 6q^2\ttwo(2\tone+\ttwo) \Gamma_1^\vee  + 6q^2\tone(\tone+2\ttwo) \Gamma_2^\vee  +O(3),
\end{equation*}
which implies the second statement.
\end{proof}

\section{Multiplication by \(\eu\)}\label{appendix:multiplication_by_E}

In this appendix we work in \(\BQH_{\Gamma_1}^{(2)}(Y)\) with \(Y \subset \IGr(2,6)\) a smooth hyperplane section.

\begin{lemma}
We have \(\eu \star 1 = 4h-2\tone \Gamma_1\) and \(\eu\star h= 4\sigma_{110}+4\sigma_{011}\). 
\end{lemma}
\begin{proof}
The first equality holds by definition of \(\eu\). For \(\ttwo=0\) we have \(s = t_1^2\), 
thus we have
\(\eu\star h=4(1+ q{t_1^2})\sigma_{110}+4(1-q{t_1^2})\sigma_{011}-2{t_1}(2q{t_1}( \sigma_{110}-\sigma_{011}))=4\sigma_{110}+4\sigma_{011}\)
\end{proof}

\begin{lemma}
We have \(\eu \star \sigma_{110}= 4\gamma_{100}+4\gamma_{010}\) and \(\eu \star \sigma=4\gamma_{010}\). 
\end{lemma}
\begin{proof}
Recall that for \({t_2}=0\), we have \(2r_1= -t_1^2\) and \(r_2={t_1^2}\). This gives 
\begin{equation*}
    h\star \sigma_{110}=(1+ \frac{q{t_1^2}}{3})\gamma_{100}+(1+\frac{q{t_1^2}}{3})\gamma_{010}-q{t_1^2}\Gamma_1+2q{t_1^2}\Gamma_2.
\end{equation*}
Moreover, using that \(\langle \sigma_{110} , \Gamma_1,\Gamma_1,\Gamma_1 \rangle_1 = 0\) by \autoref{lem_GW_prim_insertions_simplified}, we have 
\begin{equation*}
\begingroup
 \begin{array}{ll} 
   t_1\Gamma_1\star \sigma_{110}& = q{t_1^2}\langle\Gamma_1,\sigma_{110},\Gamma_1,\Gamma_2\rangle_1\Gamma_2^\vee+q{t_1^2}\langle\Gamma_1,\sigma_{110},\Gamma_1,\gamma_{100}\rangle_1\gamma_{100}^\vee+ q{t_1^2}\langle\Gamma_1,\sigma_{110},\Gamma_1,\gamma_{010}\rangle_1\gamma_{010}^\vee   \\[1.4ex]
   & =   \displaystyle{2q {t_1^2}(2\Gamma_2- \Gamma_1)+ \frac{2q{t_1^2}}{3}(\gamma_{100}+\gamma_{010}),}
    \end{array}
    \endgroup
\end{equation*}
this gives the desired equality. Similarly as before, we have 
\begin{equation*}
   h\star \sigma = \frac{2q{t_1^2}}{3}\gamma_{100}+(1-\frac{q{t_1^2}}{3})\gamma_{010}+q{t_1^2}\Gamma_1-2q{t_1^2}\Gamma_2
\text{ and } 
t_1 \Gamma_1 \star \sigma=2q{t_1^2}(\Gamma_1 - 2\Gamma_2)+\frac{2q{t_1^2}}{3}(2\gamma_{100}-\gamma_{010})
\end{equation*}

which gives the desired equality. 
\end{proof}

\begin{lemma}
We have \(\begin{aligned}[t]
\eu \star \gamma_{100}& = 12\sigma_{-110}+4\sigma_{-011}+8q \text{ and }\\
 \eu \star \gamma_{010}& =12\sigma_{-110}+8\sigma_{-011}+4q+4q^2{t_1^2}. 
 \end{aligned}\)
\end{lemma}

\begin{proof}
We have \(h\star \gamma_{100}=(3+q{t_1^2})\sigma_{-110}+(1+ q{t_1^2})\sigma_{-011}+2q\)
and
\[\begin{array}{rl}
    t_1\Gamma_1\star\gamma_{100}= & q{t_1^2}\langle\Gamma_1,\gamma_{100},\Gamma_1,\sigma_{110}\rangle_1\sigma_{-110} + q{t_1^2}\langle\Gamma_{1},\gamma_{100},\Gamma_1,\sigma\rangle_1\sigma_{-011} + q^2{t_1^2}\langle\Gamma_1,\gamma_{100},\Gamma_1,\pt\rangle_2 \\
   = &  2q{t_1}\sigma_{-110} + 2q{t_1}\sigma_{-011},
\end{array}\]
giving the first equality. We also have \(h \star \gamma_{010}=(3+ q{t_1^2})\sigma_{-110}+2\sigma_{-011}+q+ 2q^2{t_1^2}\) 
\[\begin{array}{rl}
    t_1\Gamma_1\star \gamma_{010}= & q{t_1^2}\langle\Gamma_1,\gamma_{010},\Gamma_1,\sigma_{110}\rangle_1\sigma_{-110}  +q{t_1^2}\langle\Gamma_{1},\gamma_{010},\Gamma_1,\sigma\rangle_1\sigma_{-011}+q^2{t_1^2}\langle\Gamma_1,\gamma_{010},\Gamma_1,\pt\rangle_2 \\
 =& 2q{t_1}\sigma_{-110} + 2q^2{t_1}.
\end{array}\]
proving the second equality.
\end{proof}

\begin{lemma}
We have \(\begin{aligned}[t]
\eu \star \Gamma_1&= 4{t_1}( \pt+q( \sigma_{110}+ - \sigma)) \text{ and } \\
\eu \star \Gamma_2& = 2{t_1}\pt+2q{t_1}( \sigma_{110}- \sigma)+12q^2{t_1^2}.
\end{aligned}\)
\end{lemma}

\begin{proof}
By \autoref{lem:non-geometric_GW}
 we have \(t_1\Gamma_1\star \Gamma_1= - 2 t_1 \pt +2qt_1(\sigma_{110}-\sigma)\) and 
\[\begin{array}{rl}
t_1\Gamma_1\star \Gamma_2 = & - t_1 \pt +q t_1( \sigma_{110}- \sigma)
+ q{t_1^2}\langle\Gamma_1,\Gamma_1,\Gamma_2,\sigma_{110}\rangle_1 \sigma_{-110} \\
    &  +q{t_1^2}\langle\Gamma_1,\Gamma_1,\Gamma_2,\sigma\rangle_1 \sigma_{-011} + q^2 {t_1^2}\langle\Gamma_1,\Gamma_1,\Gamma_2,\pt\rangle_2,\\
     = & t_1(-\pt +q(\sigma_{110} - \sigma)) + 6q{t_1^2} (-\sigma_{-110}+\sigma_{-011})+ 6q^2 {t_1}.
\end{array}\]
For \(t_2 = 0\), we have \(\frac{\partial s}{\partial {t_1}}=2{t_1}\), \(\frac{\partial s}{\partial {t_2}}={t_1}\),  \(\frac{\partial t}{\partial {t_1}}=0\) and \(\frac{\partial t}{\partial {t_2}}={t_1^2}\). We get 
 \(h\star \Gamma_1 = 2q{t_1}(\sigma_{110}-\sigma)\) and \(h\star\Gamma_2=q{t_1}(\sigma_{110}- \sigma)+ 6q^2{t_1^2} + 3q{t_1^2}(-\sigma_{-110}+\sigma_{-011})\). In conclusion, 
\begin{equation*}
    \eu \star \Gamma_1=8q{t_1}(\sigma_{110}- \sigma)+ 4{t_1}\pt + 4q{t_1}(-\sigma_{110}+\sigma)=4{t_1}(\pt+q( \sigma_{110}-\sigma))
\end{equation*}
\[\begin{array}{rl} 
\eu \star \Gamma_2 = &4q{t_1}( \sigma_{110}- \sigma)+ 24q^2{t_1^2}+12q{t_1^2}(-\sigma_{-110}+\sigma_{-011})\\
& + 2{t_1}\pt + 2q{t_1}(-\sigma_{110}+\sigma)- 12q{t_1^2}(-\sigma_{-110}+\sigma_{011})-12q^2{t_1^2} \\
= &  2{t_1}\pt+2q{t_1}(\sigma_{110}- \sigma)+ 12q^2{t_1^2},
\end{array}\]
giving both formulas.
\end{proof}

\begin{lemma}
We have \(\begin{aligned}[t]
\eu \star \sigma_{-110} & = 4l+4qh-2q{t_1}\Gamma_1 \text{ and } \\
\eu \star \sigma_{-011} & = 4l+ 4q^2{t_1^2}h+2q{t_1}\Gamma_1.
\end{aligned}\) 
\end{lemma}

\begin{proof}
We have \(h\star \sigma_{-110} = (1+q{t_1^2})l+qh-q{t_1}\Gamma_1\) and 
\[\begin{array}{rl} 
t_1\Gamma_1\star\sigma_{-110} = & qt_1\langle\Gamma_1,\sigma_{-110},\Gamma_1 \rangle_1\Gamma_1^\vee + qt_1\langle\Gamma_1,\sigma_{-110},\Gamma_2\rangle_1\Gamma_2^\vee\\
 &+q{t_1^2}\langle\Gamma_1,\sigma_{-110},\Gamma_1,h\rangle_1 l+ q^2{t_1^2}\langle\Gamma_1,\sigma_{-110},\Gamma_1,l\rangle_2 h \\
 = & -qt_1 \Gamma_1 +2q{t_1^2}l, \\
\end{array}\]
giving the first formula. We also have \(h\star \sigma_{-011}=(1- q{t_1^2})l +2q^2{t_1^2}h+q{t_1}\Gamma_1\) and 
\[\begin{array}{rl} 
t_1\Gamma_1\star \sigma_{-011}= & q t_1\langle\Gamma_1,\sigma_{-011},\Gamma_1\rangle_1\Gamma_1^\vee+q t_1\langle\Gamma_1,\sigma_{-011},\Gamma_2\rangle_1\Gamma_2^\vee\\
 & +q{t_1^2}\langle\Gamma_1,\sigma_{-011},\Gamma_1,h\rangle_1l +q^2{t_1^2}\langle\Gamma_1,\sigma_{-011},\Gamma_1,l\rangle_2 h \\
 = & qt_1\Gamma_1- 2q{t_1^2}l+ 2q^2{t_1^2}h
 \end{array}\]
 giving the second formula.
\end{proof}

\begin{lemma}
We have \(\begin{aligned}[t]
\eu \star \  l & = 4\pt +4q\sigma_{110} + 4q^2{t_1^2}\sigma \text{ and } \\
\eu \star \pt & = 4(q-\frac{q^2{t_1^2}}{3})\gamma_{100} + \frac{8q^2{t_1^2}}{3}\gamma_{010}+4q^2{t_1^2}\Gamma_1-8q^2{t_1^2}\Gamma_2.
\end{aligned}\)
\end{lemma}

\begin{proof}
We have \(h\star l = \pt+q\sigma_{110}+2q^2{t_1^2}\sigma\) and 
\[\begin{array}{rl} 
  t_1\Gamma_1\star l = & q{t_1^2}\langle\Gamma_1,l,\Gamma_1,[Y]\rangle_1 \pt+q^2{t_1^2}\langle\Gamma_1,l,\Gamma_1,\sigma_{-110}\rangle_2\sigma_{110} + q^2{t_1^2}\langle\Gamma_1,l,\Gamma_1,\sigma_{-011}\rangle_2\sigma_{011} \\
= &  2q^2{t_1}\sigma_{011},
 \end{array}\]
giving the first formula. We have \(h\star \pt=(q- \frac{2q^2{t_1^2}}{3})\gamma_{100} + \frac{4q^2{t_1^2}}{3}\gamma_{010}+2q^2{t_1^2}\Gamma_1-4q^2{t_1^2}\Gamma_2\)
 and 
\[\begin{array}{rl} 
t_1\Gamma_1\star \pt   = &  q^2{t_1^2}\langle\Gamma_1,\pt,\Gamma_1,\Gamma_2\rangle_2\Gamma_2^\vee+q^2{t_1^2}\langle\Gamma_1,\pt\,\Gamma_1,\gamma_{100}\rangle_2\gamma_{100}^\vee+q^2{t_1^2}\langle\Gamma_1,\pt\,\Gamma_1,\gamma_{010}\rangle_2\gamma_{010}^\vee \\
  = &   2q^2{t_1^2}(\Gamma_1 - 2\Gamma_2) + \frac{2q^2{t_1^2}}{3}(-\gamma_{100} + 2\gamma_{010}),
 \end{array}\]
 giving the second formula. 
\end{proof}

\bibliography{References}

\newcommand{\etalchar}[1]{$^{#1}$}
\begin{thebibliography}{CMMPS19}

\bibitem[ACG{\etalchar{+}}26]{atomic.ci.flag}
Boris Alexeev, Leonardo~F. Cavenaghi, Giovane Galindo, Bogdan Georgiev, Ludmil
  Katzarkov, and Pedro Antonio~Muniz Martins.
\newblock On the atomic decomposition of complete intersection in flag
  varieties.
\newblock {\em arXiv:2609.21914}, 2026.

\bibitem[BCP25]{buch:rigidity}
Anders~S. Buch, Pierre-Emmanuel Chaput, and Nicolas Perrin.
\newblock Equivariant rigidity of {Richardson} varieties.
\newblock {\em Pac. J. Math.}, 338(2):209--229, 2025.

\bibitem[Bea86]{beauville-monodromie}
Arnaud Beauville.
\newblock Le groupe de monodromie des familles universelles d'hypersurfaces et
  d'intersections compl\`etes.
\newblock In {\em Complex analysis and algebraic geometry ({G}\"ottingen,
  1985)}, volume 1194 of {\em Lecture Notes in Math.}, pages 8--18. Springer,
  Berlin, 1986.

\bibitem[Bea95]{beauville-qh}
Arnaud Beauville.
\newblock Quantum cohomology of complete intersections.
\newblock {\em Mat. Fiz. Anal. Geom.}, 2(3-4):384--398, 1995.

\bibitem[BGL{\etalchar{+}}25]{belmans2025ade}
Pieter Belmans, Sergey Galkin, Naichung~Conan Leung, Changzheng Li, Markus
  Reineke, and Rui Xiong.
\newblock {ADE} diagrams, {Hodge--Tate} hyperplane sections and semisimple
  quantum cohomology.
\newblock {\em arXiv:2509.01101}, 2025.

\bibitem[BGMP26]{verra}
Vladimiro Benedetti, J{\'e}r{\'e}my Gu{\'e}r{\'e}, Laurent Manivel, and Nicolas
  Perrin.
\newblock Quantum cohomology and birational geometry of {Verra} fourfolds.
\newblock {\em arXiv:2605.30450}, 2026.

\bibitem[BM24]{manben}
Vladimiro Benedetti and Laurent Manivel.
\newblock On the automorphisms of hyperplane sections of generalized
  {G}rassmannians.
\newblock {\em Transform. Groups}, 29(3):1115--1159, 2024.

\bibitem[BMP26]{GM}
Vladimiro Benedetti, Laurent Manivel, and Nicolas Perrin.
\newblock Quantum cohomology and irrationality of {Gushel-Mukai} fourfolds.
\newblock {\em arXiv:2603.17487}, 2026.

\bibitem[Bou02]{bourbaki}
N.~Bourbaki.
\newblock {\em Lie groups and {L}ie algebras. {C}hapters 4--6}.
\newblock Elements of Mathematics (Berlin). Springer-Verlag, Berlin, 2002.
\newblock Translated from the 1968 French original by Andrew Pressley.

\bibitem[BP22]{BP}
Vladimiro Benedetti and Nicolas Perrin.
\newblock Cohomology of hyperplane sections of (co) adjoint varieties.
\newblock {\em arXiv:2207.02089}, 2022.

\bibitem[Bri05]{Brion}
M.~Brion.
\newblock Lectures on the geometry of flag varieties.
\newblock In {\em Topics in cohomological studies of algebraic varieties},
  Trends Math., pages 33--85. Birkh\"auser, Basel, 2005.

\bibitem[CDG24]{CDG}
Giordano Cotti, Boris~A. Dubrovin, and Davide Guzzetti.
\newblock {\em Helix structures in quantum cohomology of {Fano} varieties},
  volume 2356 of {\em Lect. Notes Math.}
\newblock Cham: Springer, 2024.

\bibitem[CM93]{coll-mcgov}
David~H. Collingwood and William~M. McGovern.
\newblock {\em Nilpotent orbits in semisimple {L}ie algebras}.
\newblock Van Nostrand Reinhold Mathematics Series. Van Nostrand Reinhold Co.,
  New York, 1993.

\bibitem[CMMPS19]{ig22n}
J.~A. Cruz~Morales, A.~Mellit, N.~Perrin, and M.~Smirnov.
\newblock On quantum cohomology of {G}rassmannians of isotropic lines,
  unfoldings of {$A_n$}-singularities, and {L}efschetz exceptional collections.
\newblock {\em Ann. Inst. Fourier (Grenoble)}, 69(3):955--991, 2019.
\newblock With an appendix by Alexander Kuznetsov.

\bibitem[CMP10]{semisimplicity}
P.~E. Chaput, L.~Manivel, and N.~Perrin.
\newblock Quantum cohomology of minuscule homogeneous spaces. {III}.
  {Semi}-simplicity and consequences.
\newblock {\em Can. J. Math.}, 62(6):1246--1263, 2010.

\bibitem[CP11]{ChPe}
P.~E. Chaput and N.~Perrin.
\newblock On the quantum cohomology of adjoint varieties.
\newblock {\em Proc. Lond. Math. Soc. (3)}, 103(2):294--330, 2011.

\bibitem[Dub98]{dubrovin1998geometry}
Boris Dubrovin.
\newblock Geometry and analytic theory of {F}robenius manifolds.
\newblock In {\em Proceedings of the International Congress of Mathematicians
  1998}, pages 315--326. European Mathematical Society-EMS-Publishing House
  GmbH, 1998.

\bibitem[FW04]{fulton-woodward}
W.~Fulton and C.~Woodward.
\newblock On the quantum product of {S}chubert classes.
\newblock {\em J. Algebraic Geom.}, 13(4):641--661, 2004.

\bibitem[GMS15]{galkin2015dubrovin}
Sergey Galkin, Anton Mellit, and Maxim Smirnov.
\newblock Dubrovin's conjecture for {IG(2, 6)}.
\newblock {\em International Mathematics Research Notices},
  2015(18):8847--8859, 2015.

\bibitem[Gu{\'e}26]{guere}
J{\'e}r{\'e}my Gu{\'e}r{\'e}.
\newblock On the irrationality of cubic fourfolds.
\newblock {\em arXiv:2603.04518}, 2026.

\bibitem[Hum92]{humphreys1992reflection}
James~E. Humphreys.
\newblock {\em Reflection groups and Coxeter groups}.
\newblock Number~29. Cambridge university press, 1992.

\bibitem[IK23]{iritani-koto}
Hiroshi Iritani and Yuki Koto.
\newblock Quantum cohomology of projective bundles.
\newblock {\em arXiv:2307.03696}, 2023.

\bibitem[KKPY25]{kkpy}
L.~Katzarkov, M.~Kontsevich, T.~Pantev, and T.Y. Yu.
\newblock {B}irational {I}nvariants from {H}odge {S}tructures and {Q}uantum
  {M}ultiplication.
\newblock {\em arXiv:2508.05105}, 2025.

\bibitem[KM94]{kontsevich1994gromov}
Maxim Kontsevich and Yu~Manin.
\newblock {Gromov-Witten} classes, quantum cohomology, and enumerative
  geometry.
\newblock {\em Communications in Mathematical Physics}, 164(3):525--562, 1994.

\bibitem[KS20]{KS1}
Alexander Kuznetsov and Maxim Smirnov.
\newblock On residual categories for {Grassmannians}.
\newblock {\em Proc. Lond. Math. Soc. (3)}, 120(5):617--641, 2020.

\bibitem[KS21]{KS}
A.~Kuznetsov and M.~Smirnov.
\newblock Residual categories for (co)adjoint {G}rassmannians in classical
  types.
\newblock {\em Compos. Math.}, 157(6):1172--1206, 2021.

\bibitem[Kuz08]{kuznetsov2008exceptional}
Alexander Kuznetsov.
\newblock Exceptional collections for {G}rassmannians of isotropic lines.
\newblock {\em Proceedings of the London Mathematical Society}, 97(1):155--182,
  2008.

\bibitem[Man99]{manin}
Yuri~I. Manin.
\newblock {\em Frobenius manifolds, quantum cohomology, and moduli spaces},
  volume~47 of {\em American Mathematical Society Colloquium Publications}.
\newblock American Mathematical Society, Providence, RI, 1999.

\bibitem[Par14]{paris2014k}
Luis Paris.
\newblock {$ K (\pi, 1)$} conjecture for {A}rtin groups.
\newblock In {\em Annales de la Facult{\'e} des sciences de Toulouse:
  Math{\'e}matiques}, volume~23, pages 361--415, 2014.

\bibitem[Per14]{perrin2014semisimple}
Nicolas Perrin.
\newblock Semisimple quantum cohomology of some {F}ano varieties.
\newblock {\em arXiv:1405.5914}, 2014.

\bibitem[PS25]{maxim}
Nicolas Perrin and Maxim~N. Smirnov.
\newblock On the big quantum cohomology of coadjoint varieties.
\newblock {\em Selecta Mathematica}, 31(2):30, 2025.

\bibitem[Tev05]{tevelev}
E.~A. Tevelev.
\newblock {\em Projective duality and homogeneous spaces}, volume 133 of {\em
  Encyclopaedia of Mathematical Sciences}.
\newblock Springer-Verlag, Berlin, 2005.
\newblock Invariant Theory and Algebraic Transformation Groups, IV.

\bibitem[VdL83]{van1983homotopy}
Harm Van~der Lek.
\newblock {\em The homotopy type of complex hyperplane complements}.
\newblock PhD thesis, Katholieke Universiteit te Nijmegen, 1983.

\bibitem[Voi07]{voisin2002theorie}
Claire Voisin.
\newblock {\em Hodge theory and complex algebraic geometry. {I} \& {II}},
  volume 76 \& 77 of {\em Cambridge Studies in Advanced Mathematics}.
\newblock Cambridge University Press, Cambridge, english edition, 2007.
\newblock Translated from the French by Leila Schneps.

\end{thebibliography}
\bibliographystyle{alpha}

\end{document}